\documentclass[10pt,english]{amsart}
\usepackage{helvet}
\usepackage[T1]{fontenc}
\usepackage[latin1]{inputenc}
\usepackage{graphicx}
\advance\hoffset by -20mm
\usepackage{pifont}
\usepackage{textcomp}
\usepackage{amssymb}

\makeatletter
\numberwithin{equation}{section} 
\numberwithin{figure}{section} 
 \theoremstyle{plain}
 \newtheorem{thm}{Theorem}[subsection]
 \newtheorem*{thm*}{Theorem}
 \theoremstyle{definition}
 \newtheorem{defn}[thm]{Definition}
 \theoremstyle{plain}
 \newtheorem{lem}[thm]{Lemma}
 \newtheorem{cor}[thm]{Corollary}
 \theoremstyle{remark}
 \newtheorem{rem}[thm]{Remark}
 \theoremstyle{definition}
 \newtheorem{example}[thm]{Example}
 \theoremstyle{plain}
 \newtheorem{prop}[thm]{Proposition}

\usepackage{latexsym,amsfonts, xypic,amssymb}
\usepackage{pstricks,pst-node,pst-tree,fancybox}
\xyoption{all}

\newcommand{\p}{\mathbb{P}}
\renewcommand{\k}{\mathrm{k}}
\newcommand{\kk}{\mathrm{\overline{k}}}

\newcommand{\tr}{\!\triangleright\!}

\newcommand{\Pn}{\mathbb{P}^2} 
\newcommand{\leftexp}[2]{{\vphantom{#2}}^{#1}{\! #2}}
\newcommand{\T}[1]{\leftexp{t}{#1}}
\newcommand{\Aff}{\mathbf{Aff}}
\newcommand{\Aut}{\mathrm{Aut}}
\newcommand{\PGL}{\mathrm{PGL}}
\newcommand{\Base}{\mathcal{B}}
\newcommand{\dJo}{\mathbf{Jon}}
\newcommand{\A}{\mathbb{A}}
\theoremstyle{remark}

\theoremstyle{definition}
\newtheorem{enavant}[thm]{}

\AtBeginDocument{

}

\usepackage{babel}
\makeatother

\author{J\'er\'emy Blanc} 
\address{J\'er\'emy Blanc, Universit\'e de Gen\`eve, Section de math\'ematiques, 2-4 rue du Li\`evre
Case postale 64, 1211 Gen\`eve 4, Suisse} 
\email{Jeremy.Blanc@unige.ch} 
\author{Adrien Dubouloz} 

\address{Adrien Dubouloz, Institut de Math\'ematiques de Bourgogne, Universit\'e de Bourgogne, 9 avenue Alain Savary - BP 47870, 21078 Dijon cedex, France} 
\email{Adrien.Dubouloz@u-bourgogne.fr}
\thanks{This research has been partially supported by FABER Grant 07-512-AA-010-S-179} 
\begin{document}

\title{Automorphisms of $\mathbb{A}^{1}$-fibered affine surfaces}

\begin{abstract}{ We develop technics of birational geometry to study automorphisms of affine surfaces admitting many distinct rational fibrations, with a particular focus on the interactions between automorphisms and these fibrations.  In particular, we associate to each surface $S$ of this type a graph encoding equivalence classes of rational fibrations from which it is possible to decide for instance if the automorphism group of $S$ is generated by automorphisms preserving these fibrations. }
\end{abstract}
\maketitle

\indent\newline



\section*{Introduction}

Motivated by the example of the affine plane $\mathbb{A}^{2}$ on
which algebraic automorphisms act transitively, it is a natural problem
to determine which affine surfaces are homogeneous under the action
of their automorphism group. It turned out that it is more interesting
to consider affine surfaces that are only almost homogeneous in the
sense that the orbit of a general point has a finite complement. Indeed,
in his pioneer work, M.H. Gizatullin \cite{Gi} obtained a geometric
characterization of such surfaces in term of the structure of the
boundary divisors in minimal projective completions of these. Namely,
he established that up to finitely many exceptional cases, such
surfaces are precisely those which admit completions by so-called
\emph{zigzag}s, that is, chains of proper nonsingular rational curves.
The automorphism groups of such surfaces have been studied later on
by V.I. Danilov and M.H. Gizatullin \cite{Gi-Da1,Gi-Da2}. Motivated again
by the example of the affine plane due to J.-P. Serre \cite{Ser}, they
established in particular that these automorphism groups can be realized
as fundamental groups of graph of groups constructed from suitable
famillies of projective completions. In principle, this description
would allow to derive a more explicit presentation of these automorphism
groups. This was done by V.I. Danilov and M.H. Gizatullin in the case
of surfaces admitting a completion by an irreducible zigzag \cite{Gi-Da2}.
But in general, the corresponding graphs of groups are infinite and
it becomes very difficult even to extract any explicit description
of potentially interesting subgroups.

A noteworthy geometric feature of affine surfaces $S$ completable
by a zigzag is that they are rational, and admit $\mathbb{A}^{1}$-fibrations
$\pi:S\rightarrow\mathbb{A}^{1}$, that is, surjective morphism with
general fibers isomorphic to the affine line. Actually, except for
the case of $\mathbb{A}^{1}\setminus\left\{ 0\right\} \times\mathbb{A}^{1}$,
it turns out that every such surface admits at least two fibrations
of this type with distinct general fibers (see e.g. \cite{Dub1}). This
motivates an alternative approach consisting of understanding the automorphisms
of these surfaces in terms of their interactions with $\mathbb{A}^{1}$-fibrations.
In particular, the following questions seem natural in this context: 

1) Does the automorphism group ${\rm Aut}\left(S\right)$ of $S$
act transitively on the set of $\mathbb{A}^{1}$-fibrations on $S$
? 

2) Can ${\rm Aut}\left(S\right)$ be generated by automorphisms that
each preserves an $\mathbb{A}^{1}$-fibration ?

\noindent For instance, both questions
are answered affirmatively for the affine plane $\mathbb{A}^{2}$, as consequences of the Abhyankar-Moh Theorem \cite{AbM}
and of the Jung-van der Kulk Theorem \cite{Jun} giving the description of ${\rm Aut}\left(\mathbb{A}^{2}\right)$.
\\

In this article, we develop a general method to address these questions,
based on the study of birational relations between suitably chosen
projective models. Namely, starting with an $\mathbb{A}^{1}$-fibered
surface $\pi:S\rightarrow\mathbb{A}^{1}$, we consider projective
completions $\left(X,B,\bar{\pi}\right)$ of $S$ that we call $1$\emph{-standard}
(see \ref{StandZigDef} below), following the notation of \cite{Gi-Da1}. Here $X$ is a projective surface, $B=X\setminus S$
is a boundary zigzag, and $\pi$ extends to a rational fibration $\bar{\pi}:X\rightarrow\mathbb{P}^{1}$.
We introduce two classes of birational transformations $\phi:\left(X,B,\bar{\pi}\right)\dashrightarrow\left(X',B',\bar{\pi}'\right)$
between such completions that restrict to isomorphisms $X\setminus B\stackrel{\sim}{\rightarrow}X'\setminus B'$.
The first ones, called \emph{fibered modifications} have the property
that they are compatible with the given rational fibrations on $X$
and $X'$ respectively. The second ones, called \emph{reversions} as in \cite{FKZ},
can be thought as the simplest possible birational transformations
between such completions that are not compatible with the rational
fibrations $\bar{\pi}$ and $\bar{\pi}'$. One of the main result
of the article is the fact that these basic  birational transformations
are the building blocks for general birational maps between
$1$-standard completions preserving the complement of the boundaries. More precisely, we establish the following
result (Theorem~\ref{Thm:Factorization}, proved in Section~\ref{Sec:ProofOfTHM}). 
\begin{thm*}
Let $\phi:\left(X,B\right)\dashrightarrow\left(X',B'\right)$ be a birational map between $1$-standard pairs restricting to an isomorphism $X\setminus B\stackrel{\sim}{\rightarrow}X'\setminus B'$. If $\phi$ is not an isomorphism then it can be decomposed into a finite sequence \[ \phi=\phi_{n}\circ\cdots\circ\phi_{1}:\left(X,B\right)=\left(X_{0},B_{0}\right)\stackrel{\phi_{1}}{\dashrightarrow}\left(X_{1},B_{1}\right)\stackrel{\phi_{2}}{\dashrightarrow}\cdots\stackrel{\phi_{n}}{\dashrightarrow}\left(X_{n},B_{n}\right)=\left(X',B'\right)\] of fibered modifications and reversions between $1$-standard pairs $\left(X_{i},B_{i}\right)$, $i=1,\ldots,n$. 

Furthermore, such a factorization of minimal length is unique up to composition by isomorphisms between intermediate pairs. 
\end{thm*}
This leads in particular to a canonical procedure to factor an 
automorphism of an affine surface $S$ completable by a zigzag
considered as a birational transformation of a fixed $1$-standard completion of $S$. Using
this description, we associate to every such surface $S$ a connected
graph $\mathcal{F}_{S}$ with equivalence classes of $1$-standard completions
$\left(X,B\right)$ of $S$ as vertices and with edges being given
by reversions. This graph, which is
in general smaller than the one constructed by V.I. Danilov and M.
H. Gizatullin \cite{Gi-Da1}, encodes all the necessary information to understand
the interactions between automorphisms of $S$ and $\mathbb{A}^{1}$-fibrations
on it. For instance, we establish that under mild assumptions on $S$, 
 the automorphism group ${\rm Aut}\left(S\right)$ is generated
by automorphisms of $\mathbb{A}^{1}$-fibrations if and only if the
associated graph $\mathcal{F}_{S}$ is a tree. In general, we show that
it is also possible to equip $\mathcal{F}_{S}$ with an additional
structure of a graph of groups having ${\rm Aut}\left(S\right)$ as
its fundamental group. \\

The article is organized as follows. Section 1 introduces the basic
definition concerning $1$-standard pairs $\left(X,B\right)$ and
the existing rational fibrations on these. Section 2 contains a detailed
geometric study of fibered modifications and reversions. Section 3
is devoted to the proof of the main theorem above and section 4 presents
the construction and the interpretation of the graph $\mathcal{F}_{S}$
together with its additional structure of graph of groups. 

Finally, in section 5, we apply our general machinery to the study of classical examples
of affine surfaces completable by a zigzag. After explaining how to
recover Jung's Theorem from our description,  we consider normal
surfaces defined by an equation of the form $uv=P\left(w\right)$
in $\mathbb{A}^{3}$. In this case we not only recover generators of their
automorphism groups as obtained by M. Makar-Limanov \cite{ML90}
but we  also show that they can be equipped an additional amalgamated product structure. As a byproduct, we also recover D. Daigle
transitivity Theorem \cite{Dai} asserting that such surfaces admit a
unique equivalence class of $\mathbb{A}^{1}$-fibrations. We also prove that if the degree of $P$ is at least $3$, the group of automorphisms of the corresponding surface is not generated by automorphisms of $\mathbb{A}^1$-fibrations. In the last subsection, we give examples of an affine surfaces with the total inverse properties: in general, they admit infinitely many equivalence classes of $\mathbb{A}^1$-fibrations  but the group is generated by automorphisms of $\mathbb{A}^1$-fibrations.

\section{Preliminaries : Standard zigzags and associated rational fibrations}

 In what follows we fix a field $\k$. All varieties occuring in the sequel are implicitly assumed to be geometrically integral and defined over $\k$, and all morphisms between these are assumed to be defined over $\k$. 

\begin{defn}
A \emph{zigzag} on a normal projective surface $X$ is
a connected SNC-divisor, supported in the smooth locus of $X$, with irreducible components isomorphic to the projective line over $\k$
and whose dual graph is a chain. 

If ${\rm Supp}\left(B\right)=\bigcup_{i=0}^{r}B_{i}$ then the irreducible
components $B_{i}$, $i=0,\ldots,r$, of $B$ can be ordered in such
a way that \[
B_{i}\cdot B_{j}=\begin{cases}
1 & \textrm{if }\left|i-j\right|=1,\\
0 & \textrm{if }\left|i-j\right|>1.\end{cases}\]

A zigzag with such an ordering on the set of its components is called
\emph{oriented} and the sequence $\left((B_0)^2,\ldots,(B_r)^2\right)$
is called the \emph{type} of $B$. The components $B_{0}$ and $B_{r}$
are called the \emph{boundaries} of $B$. For an oriented zigzag $B$,
the same zigzag with the reverse ordering is denoted by $\T{B}$. 

An \emph{oriented sub-zigzag} of an oriented zigzag is an SNC divisor
$B'$ with ${\rm Supp}\left(B'\right)\subset{\rm Supp}\left(B\right)$
which is a zigzag for the induced ordering. 

We say that an oriented zigzag $B$ is composed of sub-zigzags $Z_{1},\ldots,Z_{s}$,
and we write $B=Z_{1}\tr\cdots\tr Z_{s}$, if the $Z_{i}$'s are
oriented sub-zigzags of $B$ whose union is $B$ and the components of $Z_{i}$ precede
those of $Z_{j}$ for $i<j$. 
\end{defn}
\noindent

\begin{defn}\label{StandZigDef}
A zigzag $B$ on a normal projective surface $X$ is called $m$-\emph{standard}
if it can be written as $B=F\tr C\tr E$ where $F$ and $C$ are
smooth irreducible rational curves with self-intersections $F^2=0$ and
$C^2=-m$, $m\in\mathbb{Z}$, and where $E=E_{1}\tr\cdots\tr E_{r}$
is a (possibly empty) chain of irreducible rational curves with self-intersections
$(E_i)^2\leq-2$ for every $i=1,\ldots,r$. 

An {\it $m$-standard pair} is a pair $\left(X,B\right)$ consisting
of a normal rational projective surface $X$ and an $m$-standard zigzag $B$.
A \emph{birational map $\phi: \left(X,B\right)\dasharrow\left(X',B'\right)$  between $m$-standard pairs} is a birational map $\phi:X\dasharrow X'$ which restricts to an isomorphism $X\setminus B\stackrel{{\sim}}{\rightarrow} X'\setminus B'$. 
\end{defn}

\begin{enavant} \label{underlyingSurface} Since it is rational, the underlying projective surface of an $m$-standard
pair $\left(X,B=F\tr C\tr E\right)$ comes equipped
with a rational fibration $\bar{\pi}=\bar{\pi}_{\left|F\right|}:X\rightarrow\mathbb{P}^{1}$
defined by the complete linear system $\left|F\right|$ (see e.g.
 \cite{Dub1}). In the sequel, we will implicitly consider $m$-standard
pairs as equipped with this fibration $\bar{\pi}$. Recall that the generic fiber of a rational fibration $\bar{\pi}$ is isomorphic
to the projective line over the function field of $\mathbb{P}^{1}$,
and that the total transform of the singular fibers of $\bar{\pi}$
in a minimal resolution $\mu:Y\rightarrow X$ of the singularities
of $X$ consist of trees of nonsingular rational curves (see e.g., 
Lemma 1.4.1 p. 195 in \cite{MiyBook} which remains valid over an arbitrary base field). 

The rational fibration $\bar{\pi}$ restricts on the quasi-projective
surface $S=X\setminus B$ to an a faithfully flat morphism $\pi:S\rightarrow \mathbb{A}^1$ with generic fiber isomorphic
to the affine line over the function field of $\mathbb{A}^{1}$. The general fibers of $\pi$ are isomorphic to affine lines
and $\pi$ has finitely many \emph{degenerate} fibers whose
total transforms in a minimal resolution of singularities of $S$
consists of nonempty disjoint unions of trees of rational curves,
with irreducible components isomorphic to either affine or projective lines, possibly defined over finite algebraic extensions of $\k$. 
In contrast, the restriction of $\pi$ to the complement of its degenerate
fibers has the structure of a trivial $\mathbb{A}^{1}$-bundle. 
In what follows, such morphisms will be simply refered to as $\mathbb{A}^1$-\emph{fibrations} or $\mathbb{A}^1$-\emph{fibered surfaces}. 
\end{enavant}

\begin{defn}\label{Def:A1fibredsurfaces}
We say that two $\mathbb{A}^1$-fibered surfaces $(S,\pi)$ and $(S',\pi')$ are \emph{isomorphic} if there exist an isomorphism $\Psi:S\rightarrow S'$ and an automorphism $\psi$ of $\mathbb{A}^1$ such that $\pi'\circ \Psi=\psi\circ \pi$.

On a surface $S$, two $\mathbb{A}^1$-fibrations $\pi,\pi':S\rightarrow \mathbb{A}^1$ are said to be \emph{equivalent} if $(S,\pi)$ and $(S,\pi')$ are isomorphic.
\end{defn}

\begin{enavant}
If $B$ is moreover the support an ample divisor, then $S$ is affine
and $\pi:S\rightarrow\mathbb{A}^{1}$ has a unique degenerate fiber
$\pi^{-1}\left(\bar{\pi}\left(E\right)\right)$ which consists of
a nonempty disjoint union of affine lines, again possibily defined over finite algebraic extensions of $\k$,  when equipped with its reduced
scheme structure. Furthermore, if any, the singularities of $S$ are all supported on the degenerate
fiber of $\pi$  and admit a minimal resolution whose exceptional set consists of a chain of rational curves possibly defined over a finite algebraic extension of $\k$ (this follows from the same argument as in the proof of Lemma 1.4.4 in \cite{MiyBook}). In particular, if $\k$ is algebraically closed of characteristic $0$, then $S$ has at worst Hirzebruch-Jung cyclic quotient singularities. 
\end{enavant}

\begin{enavant} Hereafter, we will mostly consider $1$-standard
pairs $\left(X,B\right)$. The simplest example $\left(\mathbb{F}_{1},F\tr C_0\right)$
consists of the Hirzebruch surface $\rho:\mathbb{F}_{1}=\mathbb{P}\left(\mathcal{O}_{\mathbb{P}^{1}}\oplus\mathcal{O}_{\mathbb{P}^{1}}\left(-1\right)\right)\rightarrow\mathbb{P}^{1}$
and the union of 
a fiber $F$ of $\rho$ and the negative section $C_{0}$ of $\rho$. More generally, we have the following description. 
\end{enavant}

\begin{lem}\label{Lem:GoingDownToF1}
Let $\left(X,B=F\tr C\tr E,\bar{\pi}\right)$ be a $1$-standard
pair and let $\mu:Y\rightarrow X$ be
the minimal resolution of the singularities of $X$. Then there
exists a birational morphism $\eta:Y\rightarrow\mathbb{F}_{1}$, unique
up to automorphisms of $\mathbb{F}_{1}$, that restricts to an isomorphism
outside the degenerate fibers of $\bar{\pi}\circ\mu$, and a commutative
diagram \[\xymatrix@R=1mm@C=1cm{&& Y \ar[dll]_{\mu} \ar[dd]^{\mu\circ \bar{\pi}} \ar[drr]^{\eta} \\ X \ar[drr]_{\bar{\pi}}& & & &\mathbb{F}_1 \ar[dll]^{\rho} \\ && \mathbb{P}^1. }\]

Furthermore, if $\left(X',B'=F'\tr C'\tr E',\bar{\pi}'\right)$
is another $1$-standard pair with associated morphism $\eta:Y'\rightarrow\mathbb{F}_{1}$
then $\left(X,B,\bar{\pi}\right)$ and $\left(X',B',\bar{\pi}'\right)$
are isomorphic if and only if there exists an automorphism of $\mathbb{F}_{1}$
mapping isomorphically $\eta\left(\mu_{*}^{-1}F\right)$  onto $\eta'\left(\left(\mu'\right)_{*}^{-1}F'\right)$
and sending isomorphically the base-points of $\eta^{-1}$ {\upshape (}including infinitely
near ones{\upshape )}  onto those of $\left(\eta'\right)^{-1}$. 
\end{lem}
\begin{proof}
Since $B$ is supported in $X_{{\rm reg}}$, its proper transform in $Y$ coincide with its total transform and is again a $1$-standard zigzag. We may therefore assume  that
$X$ is smooth.

Let us prove the first assertion. By contracting successively all the $\left(-1\right)$-curves
in the degenerate fibers $F_{1},\ldots,F_{s}$ of $\bar{\pi}$, one
obtains a birational morphism $\eta:X\rightarrow\mathbb{F}_{m}$ onto
a certain Hirzebruch surface $\rho_{m}:\mathbb{F}_{m}\rightarrow\mathbb{P}^{1}$,
which maps $C$, $F$ and the $F_{i}$'s onto a section and $s+1$
distinct fibers of $\rho_{m}$ respectively. Let
$r\left(\eta\right)=\left(\eta_{*}C\right)^2\geq-1$. If $r\left(\eta\right)=-1$
then $m=1$ and we are done. Otherwise, since $C^2=-1$
in $X$, it follows that $\eta$ contracts at least one of the irreducible
components of a degenerate fiber, say $F_{1}$, onto the point $p=\eta(F_{1})\in\eta(C)$.
Therefore, $\eta$ factors through the blow-up $\sigma:\tilde{\mathbb{F}}_{m}\rightarrow\mathbb{F}_{m}$
of $p$. Letting $\tau:\tilde{\mathbb{F}}_{m}\rightarrow\mathbb{F}_{m\pm1}$
be the contraction of the strict transform of the fiber $\rho_{m}^{-1}\left(\rho_{m}\left(p\right)\right)$,
we obtain a new birational morphism $\eta'=\tau\circ\sigma^{-1}\circ\eta:X\rightarrow\mathbb{F}_{m\pm1}$
satisfying $r\left(\eta'\right)=r\left(\eta\right)-1$. So the existence
of $\eta:X\rightarrow\mathbb{F}_{1}$ follows by induction. Suppose
that $\eta':X\rightarrow\mathbb{F}_{1}$ is another such morphism.
Then $\eta'\circ\eta^{-1}:\mathbb{F}_{1}\dashrightarrow\mathbb{F}_{1}$ is a birational map which does not blow-up any point of $\eta(C)$ and does not contract any curve intersecting $\eta(C)$. Since $\eta(C)$ is a section and $\eta'\circ\eta^{-1}$ may be decomposed into  elementary links between Hirzebrurch surfaces, $\eta'\circ\eta^{-1}$ is an isomorphism.

 The second assertion follows from the fact that
an isomorphism between $1$-standard pairs $\left(X,B,\bar{\pi}\right)$
and $\left(X',B',\bar{\pi}'\right)$ induces an isomorphism between
$B$ and $B'$ which preserves the orientation, whence descends to
an automorphism of $\mathbb{F}_{1}$. 
\end{proof}

\section{Two basic birational maps between $1$-standard pairs}
\subsection{Base-points and curves contracted}
\indent\newline We will study isomorphisms between the complements of the boundary as birational maps between $1$-standard pairs; we can distinguish two different kind of such maps, according to the following result.

\begin{lem}\label{Lem:BasePoints}
Let $\phi:\left(X,B=F\tr C\tr E\right)\dasharrow(X',B')$ be a birational map between two $1$-standard pairs, which is not an isomorphism, and let $X\stackrel{\sigma}{\leftarrow}Z\stackrel{\sigma'}{\rightarrow}X'$ be a minimal resolution of $\phi$.

Then every curve contracted by $\phi$ and every base-point of $\phi$ is defined over $\k$. Moreover, $\phi$ has a unique proper base-point $q\in B$, and one and exactly one of the following occur:
\begin{itemize}
\item[$a)$]
the strict transform of $C$ in $Z$ is the unique $(-1)$-curve contracted by $\sigma'$, and $q\in F\setminus C$;
\item[$b)$]
the strict transform of $F$ in $Z$ is the unique $(-1)$-curve contracted by $\sigma'$, and $q=F\cap C$.
\end{itemize}
$($Note that in both cases, it is possible that $F$ \emph{and} $C$ are contracted by $\phi$.$)$
\end{lem}
\begin{proof}
To any base-point of respectively $\phi$ and $\phi^{-1}$ is associated a curve contracted by respectively $\phi^{-1}$ and $\phi$. Since any curve contracted by $\phi$ and $\phi^{-1}$ is contained in the boundary, it is defined over $\k$. This implies that all base-points also are defined over $\k$.

 Each $(-1)$-curve in $Z$ which is contracted by $\sigma'$ is the proper transform of either $C$ or $F$. Since $C^2=-1$ and $C\cdot F=1$ in $X$, the two possibilities cannot occur simultaneously, so $\phi^{-1}$ (and thus $\phi$) has at most one proper base-point. If $C$ is the $(-1)$-curve contracted by $\sigma'$, to avoid a positive self-intersection for the curve $F$, there is one base-point on $F\setminus C$ (case $a$). If $F$ is the $(-1)$-curve contracted by $\sigma'$ there is one base-point on $F$; either the base-point is $F\cap C$ (case $b$), or $C$ becomes a non-negative curve, hence the $(0)$-curve of $B'$, but this implies that only one curve is contracted by $\sigma'$, a contradiction. If no $(-1)$-curve is contracted by $\sigma'$, then $\sigma'$ is an isomorphism and the discussion made above shows that so is~$\sigma$.
\end{proof}
\begin{rem}
Because of this result, when dealing with birational maps between $1$-standard pairs the fact that $\k$ is not algebraically closed, and even its characteristic is not relevant. There will only be some distinction in the last section, where the construction of the examples uses the birational morphism that blows-up points of $\mathbb{F}_1$ not necessarily defined over $\k$.
\end{rem}

\begin{defn}
If $p\in X$ is the unique proper base-point of a birational map $\phi:(X,B)\dasharrow (X',B')$ (which induces an isomorphism $X\setminus B\cong X'\setminus B'$), we say that $\phi$ is \emph{centered at $p$} and that $p$ is the \emph{center of $\phi$}.
\end{defn}
\indent In subsections \ref{SubSec:Fibered} and \ref{Sub:ZigZag}, we review two basic classes of
birational transformations between $1$-standard pairs that will play
a central role in the sequel, and are the simplest examples of maps satisfying respectively conditions $a)$ and $b)$ of Lemma~\ref{Lem:BasePoints}.

\subsection{Fibered Modifications}\label{SubSec:Fibered}

\begin{defn}
A a birational map $\phi:\left(X,B,\bar{\pi}\right)\dashrightarrow\left(X',B',\bar{\pi}'\right)$
between $\left(1\right)$-standard pairs is \emph{fibered} if it restricts to an isomorphism
of $\mathbb{A}^{1}$-fibered quasi-projective surfaces \[\xymatrix@R=0.4cm@C=1cm{S=X\setminus B \ar[d]_{\bar{\pi}\mid_{S}} \ar[r]^-{\sim}_-{\phi} & S'=X'\setminus B' \ar[d]^{\bar{\pi}'\mid_{S'}} \\ \mathbb{A}^1 \ar[r]^{\sim} & \mathbb{A}^1.}\] 
We say that  $\phi$ is a \emph{fibered modification} if it is not an isomorphism.
\end{defn}
\begin{example}
Let $\mathbb{F}_1=\{\left((x:y:z),(s:t)\right) \subset \Pn\times\mathbb{P}^1\ |\ yt=zs\}$ be the Hirzebruch surface of index $1$; the projection on the first factor yields a birational morphism $\tau:\mathbb{F}_1\rightarrow \Pn$ which is the blow-up of $(1:0:0)\in\Pn$ and the projection on the second factor yields a $\mathbb{P}^1$-bundle $\rho: \mathbb{F}_1\rightarrow\mathbb{P}^1$.  Denote by $C\subset \mathbb{F}_1$ the exceptional curve $\tau^{-1}((1:0:0))=(1:0:0)\times\mathbb{P}^1$, and by $F\subset \mathbb{F}_1$ the fiber $\rho^{-1}((0:1))$.
The map $(x,y)\mapsto ((x:y:1),(y:1))$ yields an isomorphism $\A^2\rightarrow \mathbb{F}_1\setminus (C\cup F)$.

Then every triangular automorphism
$\Psi$ of $\mathbb{A}^2$ of the form $\left(x,y\right)\mapsto\left(ax+b,cy+P\left(x\right)\right)$,
where $P\in\k\left[x\right]$, preserves the $\mathbb{A}^{1}$-fibration
$\mathrm{pr}_x=\rho|_{\mathbb{A}^2}:\mathbb{A}^2\rightarrow\mathbb{A}^{1}$ and extends to a fibered birational
map $\phi:\left(\mathbb{F}_{1},F\tr C,\rho\right)\dashrightarrow\left(\mathbb{F}_{1},F\tr C,\rho\right)$
of $1$-standard pairs. The latter is a biregular automorphism if
$\Psi$ is affine and a fibered modification otherwise. 
\end{example}
More generally, we have the following description which says in essence
that every fibered birational map between $1$-standard pairs arises
as the lift of a triangular automorphism of $\mathbb{A}^2$ as above. 

\begin{lem}\label{Lem:LiftTriangularC2}
Let $\phi:\left(X,B=F\tr C\tr E,\bar{\pi}\right)\dashrightarrow\left(X',B'=F'\tr C'\tr E',\bar{\pi}'\right)$
be a birational map between $1$-standard pairs and let $X\stackrel{\mu}{\leftarrow}Y\stackrel{\eta}{\rightarrow}\mathbb{F}_{1}$
and $X\stackrel{\mu'}{\leftarrow}Y'\stackrel{\eta'}{\rightarrow}\mathbb{F}_{1}$
be the morphisms constructed in Lemma~$\ref{Lem:GoingDownToF1}$. Then the following are equivalent
: 

$a)$ $\phi$ restricts to an isomorphism of $\mathbb{A}^{1}$-fibered
surfaces $\left(X\setminus B,\bar{\pi}\right)\stackrel{\sim}{\rightarrow}\left(X'\setminus B',\bar{\pi}'\right)$;

$b)$ $\left(\mu'\right)^{-1}\circ\phi\circ\mu:Y\dashrightarrow Y'$
is the lift via $\eta$ and $\eta'$ of an isomorphism of $\mathbb{A}^{1}$-fibered
affine surfaces \[\xymatrix@R=0.4cm@C=1cm{ \mathbb{A}^2=\mathbb{F}_1\setminus (\eta(F)\cup\eta(C)) \ar[d]_{\rho\mid_{\mathbb{A}^2}} \ar[r]_-{\Psi}^-{\sim} & \mathbb{A}^2=\mathbb{F}_1\setminus (\eta'(F')\cup\eta'(C')) \ar[d]^{\rho\mid_{\mathbb{A}^2}} \\ \mathbb{A}^1 \ar[r]_{\psi}^{\sim} & \mathbb{A}^1 }\] 
which maps isomorphically the base-points of $\eta^{-1}$  onto those of $\eta'^{-1}$. 

Furthermore, $\phi:\left(X,B\right)\dashrightarrow$$\left(X',B'\right)$
is an isomorphism if and only if $\Psi$ is affine.
\end{lem}
\begin{proof}
One checks that $\phi:\left(X,B,\bar{\pi}\right)\dashrightarrow\left(X',B',\bar{\pi}'\right)$
restricts to an isomorphism between the $\mathbb{A}^{1}$-fibered
surfaces $(X\setminus B,\bar{\pi})$ and $(X\setminus B,\bar{\pi})$
if and only if its lift $\left(\mu'\right)^{-1}\circ\phi\circ\mu$
restricts to an isomorphism between the $\mathbb{A}^{1}$-fibered
surfaces $(Y\setminus\mu_{*}^{-1}B,\bar{\pi}\circ\mu)$ and $(Y'\setminus\left(\mu'\right)_{*}^{-1}B',\bar{\pi}'\circ\mu')$.
We may thus assume that $X$ and $X'$ are smooth.

Suppose that $\phi$ is not an isomorphism and let $X'\stackrel{\sigma'}{\leftarrow}Z\stackrel{\sigma}{\rightarrow}X$
be a minimal resolution of $\phi$ where $\sigma'$ and $\sigma$
are sequences of blow-ups with centers outside $S'$ and $S$ respectively.

Assume that $\phi$ satisfies (a), which implies that the rational fibrations $\bar{\pi}'$ and $\bar{\pi}$
lift to a same rational fibration $\tilde{\pi}:Z\rightarrow\mathbb{P}^{1}$,
and that the proper transforms of $C'$ and $C$ in $Z$ coincide
with the unique section $\tilde{C}$ of $\tilde{\pi}$ contained in
the boundary $D=\left(\sigma'\right)^{-1}\left(B'\right)=\sigma^{-1}\left(B\right)$
in $Z$. Thus $\phi$ restricts to a birational map $C\dasharrow C'$. The only $(-1)$-curve of $Z$ which is contracted by $\sigma$ (respectively $\sigma'$) is therefore the proper transform of $F$ (respectively of $F'$) in $X$. Consequently, $\phi$ restricts to an isomorphism $E\rightarrow E'$; $\phi$ actually restricts to an isomorphism
of $\mathbb{A}^{1}$-fibered surfaces 
$\left(X\setminus\left(F\cup C\right),\bar{\pi}\right)\simeq\left(X'\setminus\left(F'\cup C'\right),\bar{\pi}'\right)$. Conversely, such isomorphisms extend to birational maps satisfying (a).

Now the equivalence follows from the one-to-one correspondence between
such isomorphisms and those of the form $\Psi:\left(\mathbb{F}_{1}\setminus\eta_{*}\left(F\cup C\right),\rho\right)\stackrel{\sim}{\rightarrow}\left(\mathbb{F}_{1}\setminus\eta'_{*}\left(F'\cup C' \right),\rho\right)$
which map isomorphically the base-points of $\eta^{-1}$ onto those of $\left(\eta'\right)^{-1}$. The last assertion follows from
the fact that $\phi$ extends to an isomorphism $X\stackrel{\sim}{\rightarrow}X'$
if and only if the corresponding automorphism $\Psi$ of $\mathbb{A}^2$
extends to an automorphism of $\mathbb{F}_{1}$ (both conditions are equivalent to say that the proper transform of $F$ is not contracted).
\end{proof}
\begin{enavant}\label{DiscussionFiberedModification} It follows from the above description (Lemmas \ref{Lem:BasePoints} and \ref{Lem:LiftTriangularC2}) that a fibered
modification \[
\phi:\left(X,B=F\tr C\tr E,\bar{\pi}\right)\dashrightarrow\left(X',B'=F'\tr C'\tr E',\bar{\pi}'\right)\]
has a unique proper base-point $q=F\cap C$. Letting $\Psi:\mathbb{A}^2\stackrel{\sim}{\rightarrow}\mathbb{A}^2,\left(x,y\right)\mapsto\left(ax+b,cy+P\left(x\right)\right)$
be the triangular automorphism associated with $\phi$, one checks
that the total transform $\tilde{B}$ of $B$ in a minimal resolution
$\left(X,B\right)\stackrel{\sigma}{\leftarrow}(Z,\tilde{B})\stackrel{\sigma'}{\rightarrow}\left(X',B'\right)$
of $\phi$ is a tree of rational curves with the following dual graph 

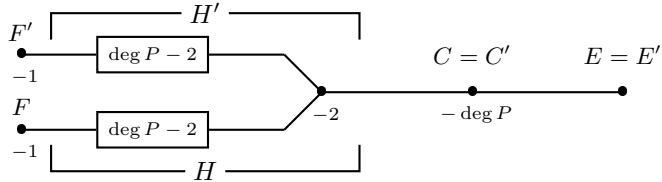
\begin{pspicture}(-2,-0.1)(8,2.3)
\psline(5,1)(9,1)
\rput(9,1){\textbullet}\rput(9,1.5){{\small $E=E'$}}
\rput(7,1){\textbullet}\rput(7,1.5){{\small $C=C'$}}\rput(7.05,0.7){{\scriptsize $-\deg P$}}
\rput(5,1){\textbullet}\rput(5.05,0.7){{\scriptsize $-2$}}
\psline(4.5,0.5)(5,1)\psline(3.5,0.5)(4.5,0.5)
\psframe(2,0.75)(3.5,0.25)\rput(2.75,0.5){{\scriptsize $\deg P -2$}}
\psline(1,0.5)(2,0.5)
\rput(1,0.5){\textbullet}\rput(1,0.8){{\small $F$}}\rput(1.05,0.2){{\scriptsize $-1$}}
\psline(4.5,1.5)(5,1)\psline(3.5,1.5)(4.5,1.5)
\psframe(2,1.25)(3.5,1.75)\rput(2.75,1.5){{\scriptsize$\deg P -2$}}
\psline(1,1.5)(2,1.5)
\rput(1,1.5){\textbullet}\rput(1,1.8){{\small $F'$}}\rput(1.05,1.2){{\scriptsize $-1$}}
\pnode(5.5,0.3){A}\pnode(1.4,0.3){B}\ncbar[angle=-90]{A}{B}
\ncput*{$H$}
\pnode(5.5,1.7){C}\pnode(1.4,1.7){D}\ncbar[angle=90]{C}{D}
\ncput*{$H'$}
\end{pspicture}

\noindent where the two boxes represent chains of $\deg P-2$ $\left(-2\right)$-curves.
Furthermore, the morphisms $\sigma:Z\rightarrow X$ and $\sigma':Z\rightarrow X$
are given by the smooth contractions of the sub-trees $H\cup H'\cup F'$
and $H\cup H'\cup F$ onto the proper base-points $q=F\cap C$ and
$q'=F'\cap C'$ of $\phi$ and $\phi^{-1}$ respectively. 

\end{enavant}

\subsection{Zigzag Reversions}\label{Sub:ZigZag}

\begin{defn}\label{Def:Reversions}
A strictly birational map $\phi:\left(X,B=F\tr C\tr E\right)\dashrightarrow\left(X',B'=F'\tr C'\tr E'\right)$
between $1$-standard pairs is called a\emph{ reversion} if it admits
a resolution of the form
\[\xymatrix@R=0.1cm@C=1cm{ & (Z,\tilde{B}=\T{(C\tr E)} \tr H\tr (C'\tr E')) \ar[dl]_{\sigma} \ar[dr]^{\sigma'} \\ \left(X,\T{B}=\T{(C\tr E)} \tr F \right) \ar@{-->}[rr]^{\phi} & & \left(X', B'=F'\tr (C' \tr F') \right)}\]
where $H$ is a zigzag with boundaries $F$ (left) and $F'$ (right) and where $\sigma:Z\rightarrow X$
and $\sigma':Z\rightarrow X'$ are smooth contractions of the sub-zigzags
$H\tr\left(C'\tr E'\right)$ and $\T{\left(C\tr E\right)}\tr H$
of $\tilde{B}$ onto their left and right boundaries $F$ and $F'$
respectively. 
\end{defn}
\begin{example}\label{Exa:32resol}
Let $n_1,n_2\geq 3$ be two integers, let $Z$ be a normal rational projective surface and let $\tilde{B}\subset Z^{\mathrm{reg}}$ be a zigzag having the following dual graph
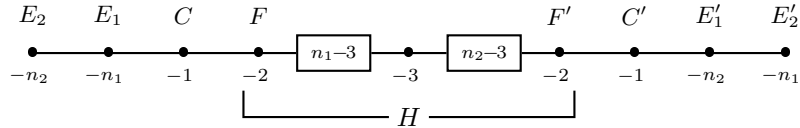
\begin{center}
\begin{pspicture}(0,0)(12,2)
\psline(1,1)(4.5,1)
\rput(1,1){\textbullet}\rput(1,1.5){{\small $E_2$}}\rput(0.95,0.7){{\scriptsize $-n_2$}}
\rput(2,1){\textbullet}\rput(2,1.5){{\small $E_1$}}\rput(1.95,0.7){{\scriptsize $-n_1$}}
\rput(3,1){\textbullet}\rput(3,1.5){{\small $C$}}\rput(2.95,0.7){{\scriptsize $-1$}}
\rput(4,1){\textbullet}\rput(4,1.5){{\small $F$}}\rput(3.95,0.7){{\scriptsize $-2$}}
\psframe(4.5,0.75)(5.5,1.25)\rput(5,1){\scriptsize $n_1\!\!-\!\!3$}
\psline(5.5,1)(6.5,1)
\rput(6,1){\textbullet}\rput(5.95,0.7){{\scriptsize $-3$}}
\psframe(6.5,0.75)(7.5,1.25)\rput(7,1){\scriptsize $n_2\!\!-\!\!3$}
\psline(7.5,1)(11,1)
\rput(8,1){\textbullet}\rput(8,1.5){{\small $F'$}}\rput(7.95,0.7){{\scriptsize $-2$}}
\rput(9,1){\textbullet}\rput(9,1.5){{\small $C'$}}\rput(8.95,0.7){{\scriptsize $-1$}}
\rput(10,1){\textbullet}\rput(10,1.5){{\small $E_1'$}}\rput(9.95,0.7){{\scriptsize $-n_2$}}
\rput(11,1){\textbullet}\rput(11,1.5){{\small $E_2'$}}\rput(10.95,0.7){{\scriptsize $-n_1$}}
\pnode(3.8,0.5){A}\pnode(8.2,0.5){B}\ncbar[angle=-90]{A}{B}
\ncput*{$H$}
\end{pspicture}\end{center}where the boxes represent chains of $n_1-3$ and $n_2-3$ $\left(-2\right)$-curves respectively.
One checks that there exist two birational morphisms
$\sigma:(Z,\tilde{B})\rightarrow\left(X,\T{B}=E_2\tr E_1\tr C\tr F\right)$ and 
$\sigma':(Z,\tilde{B})\rightarrow\left(X',B'=F'\tr C'\tr E'_{1}\tr E'_{2}\right)$
consisting of a sequence of smooth blow-downs of irreducible components
of $\tilde{B}$ starting with those of $C'$ and $C$ respectively. By construction, $\sigma'\circ\sigma^{-1}:\left(X,B\right)\dashrightarrow\left(X',B'\right)$
is a reversion between $1$-standard pairs of type $(0,-1,-n_1,-n_2)$ and $(0,-1,-n_2,-n_1)$ respectively. 
\end{example}
The following lemma summarizes some of the main properties of reversions. 

\begin{lem}[Properties of reversions]\label{Lem:RevProperties}
Let $\phi:\left(X,B=F\tr C\tr E\right)\dashrightarrow\left(X',B'=F'\tr C'\tr E'\right)$
be a reversion and let \[
\left(X,\T{B}\right)\stackrel{\sigma}{\leftarrow}\left(Z,\tilde{B}=\T{\left(C\tr E\right)}\tr H\tr\left(C'\tr E'\right)\right)\stackrel{\sigma'}{\rightarrow}\left(X',B'\right)\]
 be a \emph{minimal} resolution of $\phi$ as in Definition~$\ref{Def:Reversions}$
above. Then, the unique proper base-point of $\phi$ $($respectively $\phi^{-1})$ belongs to $F\setminus C$ $($respectively $F'\setminus C')$; moreover the following hold :

\indent$a)$ The sub-zigzag $H$ of $\tilde{B}$ consists of a
unique smooth rational curve $(F=H=F')$ if and only if $E$ consists
of a possibly empty chain of $\left(-2\right)$-curves. 

\indent $b)$ Otherwise, $($if $E$ contains at least an irreducible
component with self-intersection $\leq-3)$ then $H$ can be written
as $H=F\tr\tilde{H}\tr F'$ and the morphisms $\sigma:Z\rightarrow X$
and $\sigma':Z\rightarrow X'$ are the contractions of the sub-chains
$\tilde{H}\tr F'\tr\left(C'\tr E'\right)$ and $\T{\left(C\tr E\right)}\tr F\tr\tilde{H}$ to a point of $X$ and $X'$ respectively.

\indent $c)$ If $E$ is not empty then the birational morphisms
$\sigma:Z\rightarrow X$ and $\sigma':Z\rightarrow X'$ factor into
unique sequences of smooth blow-downs starting with the contractions
of the $\left(-1\right)$-curve $C'$ and $C$ respectively and ending
with the contractions of the right boundaries of $E'$ and $E$.
\end{lem}
\begin{proof}
Since $C$ is not affected by $\sigma$, it has self-intersection $-1$ in $Z$; since $C$ is contracted by $\sigma'$, the unique proper base-point of $\phi$  belongs to $F\setminus C$ (Lemma \ref{Lem:BasePoints}). Exchanging $C$ and $C'$ yields the analogue result for $\phi^{-1}$.

  Since the exceptional locus a smooth
contraction cannot contain two $\left(-1\right)$-curves which intersect, it follows that $F$ or $F'$ has self-intersection
$-1$ in $Z$ if and only if $H=F=F'$ and $E=E'=\emptyset$; this is a degenerate case of $a)$.

We may now assume that $F^2, (F')^2\leq -2$ in $Z$.
Let us prove that $\sigma'$ does not factor through the contraction $\eta':Z\rightarrow \tilde{X'}$ of $H'\tr C'\tr E'$, for some strict subzigzag $H'\subset (\overline{H\setminus F})$. 
Suppose the contrary. Since $\sigma$ contracts a connected curve -- which is $(\overline{H\setminus F})\tr C'\tr E'$ -- the same holds for $\eta$. Then $\eta(\overline{H\setminus F})$ is a contractible connected curve, containing a unique $(-1)$-curve which is its right boundary. This implies (since $F$ is the right boundary of $H$ and  $F^2\leq -2$ in $Z$) that $F$ has self-intersection $\leq -1$ in $X'$, a contradiction.

This observation proves the following two results: (i) if $C'\tr E'$ is contractible -- which is equivalent to say that each component of $E'$ has self-intersection $-2$ -- then $H=F$; the converse being obvious we obtain assertion a) and b).
(ii) if $E'$ is not empty the last curve contracted by $\sigma'$ is the right boundary of $E'$. The same argument for $\sigma$ achieves to prove $c)$.\end{proof}

\begin{enavant}[Description of reversions between $0$-standard pairs by means of elementary links]\label{DescriptionReversion0}\noindent\\
A \emph{reversion} between $0$-standard pairs was introduced in \cite{FKZ}. Given some pair  with a zigzag of type $(...,n_1,0,n_2,...)$, the blow-up of the point on the $(0)$-curve which also belongs to the next component, followed by the contraction of the proper transform of the $(0)$-curve yields to a pair  with a zigzag of type $(...,n_1+1,0,n_2-1,...)$. Starting from a $0$-standard pair $(X,B)$ of type $(-n_r,...,-n_1,0,0)$, one can then construct a birational map $\varphi_1:(X,B)\dasharrow (X_1,B_1)$ to a pair with a zigzag of type $(-n_r,...,-n_2,0,0,-n_1)$. Repeating this process yields birational maps $\varphi_1,...,\varphi_r$, and a reversion $\phi=\varphi_r\circ...\circ\varphi_1:(X,B)\dasharrow (X_r,B_r)$, where $B_r$ has type $(0,0,-n_r,...,-n_1)$.
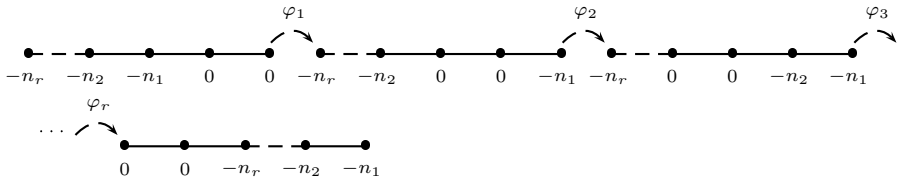
\begin{center}
\begin{tabular}{l}
\begin{pspicture}(0,0.4)(3.76,1.8)
\psline(0.8,1)(3.2,1)
\psline[linestyle=dashed](0,1)(0.8,1)
\rput(0,1){\textbullet}\rput(-0.05,0.7){{\scriptsize $-n_r$}}
\rput(0.8,1){\textbullet}\rput(0.75,0.7){{\scriptsize $-n_2$}}
\rput(1.6,1){\textbullet}\rput(1.55,0.7){{\scriptsize $-n_1$}}
\rput(2.4,1){\textbullet}\rput(2.4,0.7){{\scriptsize $0$}}
\rput(3.2,1){\textbullet}\rput(3.2,0.7){{\scriptsize $0$}}
\rput(3.53,1.55){{\scriptsize $\varphi_1$}}
\parabola[linestyle=dashed]{->}(3.22,1.15)(3.5,1.3)
\end{pspicture}
\begin{pspicture}(0,0.4)(3.76,1.5)
\psline(0.8,1)(3.2,1)
\psline[linestyle=dashed](0,1)(0.8,1)
\rput(0,1){\textbullet}\rput(-0.05,0.7){{\scriptsize $-n_r$}}
\rput(0.8,1){\textbullet}\rput(0.75,0.7){{\scriptsize $-n_2$}}
\rput(1.6,1){\textbullet}\rput(1.6,0.7){{\scriptsize $0$}}
\rput(2.4,1){\textbullet}\rput(2.4,0.7){{\scriptsize $0$}}
\rput(3.2,1){\textbullet}\rput(3.15,0.7){{\scriptsize $-n_1$}}
\rput(3.53,1.55){{\scriptsize $\varphi_2$}}
\parabola[linestyle=dashed]{->}(3.22,1.15)(3.5,1.3)
\end{pspicture}
\begin{pspicture}(0,0.4)(3.76,1.5)
\psline(0.8,1)(3.2,1)
\psline[linestyle=dashed](0,1)(0.8,1)
\rput(0,1){\textbullet}\rput(-0.05,0.7){{\scriptsize $-n_r$}}
\rput(0.8,1){\textbullet}\rput(0.8,0.7){{\scriptsize $0$}}
\rput(1.6,1){\textbullet}\rput(1.6,0.7){{\scriptsize $0$}}
\rput(2.4,1){\textbullet}\rput(2.35,0.7){{\scriptsize $-n_2$}}
\rput(3.2,1){\textbullet}\rput(3.15,0.7){{\scriptsize $-n_1$}}
\parabola[linestyle=dashed]{->}(3.22,1.15)(3.5,1.3)
\rput(3.53,1.55){{\scriptsize $\varphi_3$}}\end{pspicture}\\
\begin{pspicture}(2.6,0.4)(3.76,1.5)
\psline[linestyle=dotted](2.75,1.2)(3.05,1.2)
\parabola[linestyle=dashed]{->}(3.22,1.15)(3.5,1.3)
\rput(3.53,1.55){{\scriptsize $\varphi_r$}}
\end{pspicture}
\begin{pspicture}(0,0.4)(3.2,1.5)
\psline(0,1)(1.6,1)
\psline(2.4,1)(3.2,1)
\psline[linestyle=dashed](1.6,1)(2.4,1)
\rput(0,1){\textbullet}\rput(0,0.7){{\scriptsize $0$}}
\rput(0.8,1){\textbullet}\rput(0.8,0.7){{\scriptsize $0$}}
\rput(1.6,1){\textbullet}\rput(1.55,0.7){{\scriptsize $-n_r$}}
\rput(2.4,1){\textbullet}\rput(2.35,0.7){{\scriptsize $-n_2$}}
\rput(3.2,1){\textbullet}\rput(3.15,0.7){{\scriptsize $-n_1$}}
\end{pspicture}\end{tabular}\end{center}
The construction also decomposes the reversion into birational maps $\varphi_i$, where each $\varphi_i$ preserves the $\mathbb{A}^1_{*}$-fibration on the open part that is given by the $(0)$-curve involved.
However, the disadvantage of the decomposition is that $(\varphi_i)^{-1}$ and $\varphi_{i+1}$ have the same proper base-point, which is the intersection of the two $(0)$-curves of $X_{i+1}$.\footnote[2]{Note also that the same problem holds when dealing with reversion and fibered modification on $0$-standard pairs, which have the same proper base-point. There is thus no analogue of Lemma~\ref{Lem:BasePoints} for $0$-standard pairs.}
\end{enavant}
\begin{enavant}[Description of reversions between $1$-standard pairs by means of elementary links]\label{DescriptionReversion1}\noindent\\
On $1$-standard pairs, the analogue of construction~\ref{DescriptionReversion0} is possible. We start with a pair $(X,B)$ of type $(-n_r,...,-n_2,-n_1,-1,0)$. We choose a point $p\in X$ that belongs to the $(0)$-curve of $B$ but not to its $(-1)$-curve. The contraction of the $(-1)$-curve of $B$ followed by the blow-up of $p$ yields a birational map $\theta_0:(X,B)\dasharrow (X_0,B_0)$ to a pair with a zigzag of type $(-n_r,...,-n_2,-n_1+1,0,-1)$. As before, we can produce a birational map $\varphi_1:(X_0,B_0)\dasharrow (X_1',B_1')$, where $B_1'$ is of type $(-n_r,...,-n_2,-1,0,-n_1+1)$. The blow-down of the $(-1)$-curve followed by the blow-up of the point of intersection of the $(0)$-curve with the curve immediately after it yields a birational map $\theta_1:(X_1',B_1')\dasharrow (X_1,B_1)$ where $B_1$ is a zigzag of type $(-n_r,...,-n_2+1,0,-1,-n_1)$. Repeating this process yields birational maps $\theta_0,\varphi_1,\theta_1,...,\varphi_r,\theta_r$ described by the following figure.
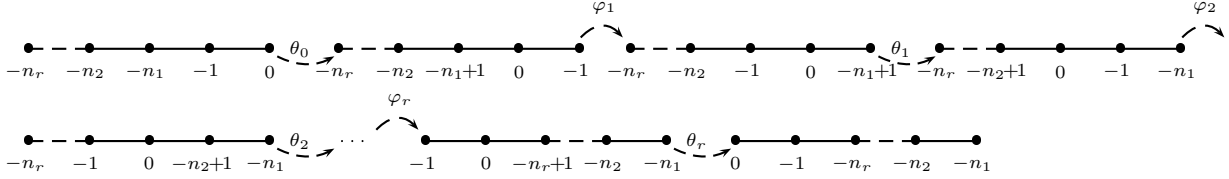
\begin{center}\begin{tabular}{l}
\begin{pspicture}(0,0.4)(4,1.8)
\psline(0.8,1)(3.2,1)
\psline[linestyle=dashed](0,1)(0.8,1)
\rput(0,1){\textbullet}\rput(-0.05,0.7){{\scriptsize $-n_r$}}
\rput(0.8,1){\textbullet}\rput(0.75,0.7){{\scriptsize $-n_2$}}
\rput(1.6,1){\textbullet}\rput(1.55,0.7){{\scriptsize $-n_1$}}
\rput(2.4,1){\textbullet}\rput(2.35,0.7){{\scriptsize $-1$}}
\rput(3.2,1){\textbullet}\rput(3.2,0.7){{\scriptsize $0$}}
\rput(3.6,1){{\scriptsize $\theta_0$}}
\parabola[linestyle=dashed]{->}(3.3,0.9)(3.68,0.8)
\end{pspicture}
\begin{pspicture}(0,0.4)(3.76,1.5)
\psline(0.8,1)(3.2,1)
\psline[linestyle=dashed](0,1)(0.8,1)
\rput(0,1){\textbullet}\rput(-0.05,0.7){{\scriptsize $-n_r$}}
\rput(0.8,1){\textbullet}\rput(0.75,0.7){{\scriptsize $-n_2$}}
\rput(1.6,1){\textbullet}\rput(1.55,0.7){{\scriptsize $-n_1\!\!+\!\!1$}}
\rput(2.4,1){\textbullet}\rput(2.4,0.7){{\scriptsize $0$}}
\rput(3.2,1){\textbullet}\rput(3.15,0.7){{\scriptsize $-1$}}
\rput(3.53,1.55){{\scriptsize $\varphi_1$}}
\parabola[linestyle=dashed]{->}(3.22,1.15)(3.5,1.3)
\end{pspicture}
\begin{pspicture}(0,0.4)(4,1.5)
\psline(0.8,1)(3.2,1)
\psline[linestyle=dashed](0,1)(0.8,1)
\rput(0,1){\textbullet}\rput(-0.05,0.7){{\scriptsize $-n_r$}}
\rput(0.8,1){\textbullet}\rput(0.75,0.7){{\scriptsize $-n_2$}}
\rput(1.6,1){\textbullet}\rput(1.55,0.7){{\scriptsize $-1$}}
\rput(2.4,1){\textbullet}\rput(2.4,0.7){{\scriptsize $0$}}
\rput(3.2,1){\textbullet}\rput(3.15,0.7){{\scriptsize $-n_1\!\!+\!\! 1$}}
\rput(3.6,1){{\scriptsize $\theta_1$}}\parabola[linestyle=dashed]{->}(3.3,0.9)(3.68,0.8)\end{pspicture}
\begin{pspicture}(0,0.4)(3.76,1.5)
\psline(0.8,1)(3.2,1)
\psline[linestyle=dashed](0,1)(0.8,1)
\rput(0,1){\textbullet}\rput(-0.05,0.7){{\scriptsize $-n_r$}}
\rput(0.8,1){\textbullet}\rput(0.75,0.7){{\scriptsize $-n_2\!\!+\!\! 1$}}
\rput(1.6,1){\textbullet}\rput(1.6,0.7){{\scriptsize $0$}}
\rput(2.4,1){\textbullet}\rput(2.35,0.7){{\scriptsize $-1$}}
\rput(3.2,1){\textbullet}\rput(3.15,0.7){{\scriptsize $-n_1$}}
\rput(3.53,1.55){{\scriptsize $\varphi_2$}}
\parabola[linestyle=dashed]{->}(3.22,1.15)(3.5,1.3)\end{pspicture}
\\
\begin{pspicture}(0,0.4)(4,1.5)
\psline(0.8,1)(3.2,1)
\psline[linestyle=dashed](0,1)(0.8,1)
\rput(0,1){\textbullet}\rput(-0.05,0.7){{\scriptsize $-n_r$}}
\rput(0.8,1){\textbullet}\rput(0.75,0.7){{\scriptsize $-1$}}
\rput(1.6,1){\textbullet}\rput(1.6,0.7){{\scriptsize $0$}}
\rput(2.4,1){\textbullet}\rput(2.3,0.7){{\scriptsize $-n_2\!\!+\!\! 1$}}
\rput(3.2,1){\textbullet}\rput(3.15,0.7){{\scriptsize $-n_1$}}
\rput(3.6,1){{\scriptsize $\theta_2$}}\parabola[linestyle=dashed]{->}(3.3,0.9)(3.68,0.8)\end{pspicture}\begin{pspicture}(2.6,0.4)(3.76,1.5)
\psline[linestyle=dotted](2.75,1)(3.05,1)
\parabola[linestyle=dashed]{->}(3.22,1.15)(3.5,1.3)
\rput(3.53,1.55){{\scriptsize $\varphi_r$}}
\end{pspicture}
\begin{pspicture}(0,0.4)(4,1.5)
\psline(0,1)(1.6,1)
\psline(2.4,1)(3.2,1)
\psline[linestyle=dashed](1.6,1)(2.4,1)
\rput(0,1){\textbullet}\rput(-0.05,0.7){{\scriptsize $-1$}}
\rput(0.8,1){\textbullet}\rput(0.8,0.7){{\scriptsize $0$}}
\rput(1.6,1){\textbullet}\rput(1.55,0.7){{\scriptsize $-n_r\!\!+\!\! 1$}}
\rput(2.4,1){\textbullet}\rput(2.35,0.7){{\scriptsize $-n_2$}}
\rput(3.2,1){\textbullet}\rput(3.15,0.7){{\scriptsize $-n_1$}}
\rput(3.6,1){{\scriptsize $\theta_r$}}\parabola[linestyle=dashed]{->}(3.3,0.9)(3.68,0.8)\end{pspicture}
\begin{pspicture}(0,0.4)(3.5,1.5)
\psline(0,1)(1.6,1)
\psline(2.4,1)(3.2,1)
\psline[linestyle=dashed](1.6,1)(2.4,1)
\rput(0,1){\textbullet}\rput(0,0.7){{\scriptsize $0$}}
\rput(0.8,1){\textbullet}\rput(0.75,0.7){{\scriptsize $-1$}}
\rput(1.6,1){\textbullet}\rput(1.55,0.7){{\scriptsize $-n_r$}}
\rput(2.4,1){\textbullet}\rput(2.35,0.7){{\scriptsize $-n_2$}}
\rput(3.2,1){\textbullet}\rput(3.15,0.7){{\scriptsize $-n_1$}}
\end{pspicture}
\end{tabular}\end{center}
Then, the composition $\phi=\theta_r\varphi_r\cdots\theta_1\varphi_0\theta_0$ is a birational  map $\phi:(X,B)\dasharrow(X_r,B_r)$ between two $1$-standard pairs. 
\end{enavant}

\begin{lem}\label{Lem:DescrRev1isaRev}
The map $\phi$ defined in $\S\ref{DescriptionReversion1}$ is a reversion between the two $1$-standard pairs $(X,B)$ and $(X',B')$.
\end{lem}
\begin{proof}
Since $\phi$ is a birational map of pairs, it has one proper base-point only (Lemma~\ref{Lem:BasePoints}), which is $p$: the unique proper base-point of $\theta_0$. Denote by $\sigma:Z\rightarrow X$ the blow-up of the base-points of $\phi$, so that $\sigma'=\phi\sigma$ is a morphism. If $q$ is a base-point of $\phi$, distinct from $p$, then $q$ is infinitely near to $p$ and corresponds to a base-point of some $\varphi_i$ or some $\theta_i$; so $q$ belongs to exactly two components of the total transform of $B$. Consequently, the total transform $\tilde{B}$ of $B$ in $Z$ is a zigzag, equal to $B\tr H$, for some zigzag $H$ (here $B\subset Z$ is the strict transform of $B\subset X$). Doing the same for $\phi^{-1}$ shows that the resolution given by $\sigma$ and $\sigma'$ satisfies the properties of Definition~\ref{Def:Reversions}. 
\end{proof}
\begin{prop}[Unicity of reversions]\label{Prop:UnicityReversions}
For every $1$-standard pair $\left(X,B=F\tr C\tr E\right)$ and
every point $q\in F\setminus C$, there exists a $1$-standard pair
$\left(X',B'\right)$ and a reversion $\phi:\left(X,B\right)\dashrightarrow\left(X',B'\right)$, unique up to an isomorphism at the target,
having $q$ as a unique proper base-point. Furthermore, if $B$ is
of type $\left(0,-1,-n_{1},\ldots,-n_{r}\right)$ then $B'$ is of
type $\left(0,-1,-n_{r},\ldots,-n_{1}\right)$. 
\end{prop}
\begin{proof}
The existence follows from $\S\ref{DescriptionReversion1}$ and Lemma~$\ref{Lem:DescrRev1isaRev}$ (it was also described in the proof of \cite[Proposition 2.10]{Dub1}).

It remains to prove unicity. For $i=1,2$, let $\phi_i:\left(X,B=F\tr C\tr E\right)\dashrightarrow\left(X_i,B_i=F_i\tr C_i\tr E_i\right)$ be a reversion centered at $q\in F\setminus C$,  admitting a minimal resolution $(X,B)\stackrel{\sigma_i}{\leftarrow}(Z_i,\tilde{B}_i)\stackrel{\sigma_i'}{\rightarrow}(X_i,B_i)$ such that $\tilde{B}_i=(\T{E_i}\tr C_i) \tr H_i\tr (C\tr E)$. Denoting by $\eta:(\hat{X},\hat{B})\rightarrow (X,B)$ the blow-up of the common base-points of $\sigma_1$ and $\sigma_2$, we have a commutative diagram
\[\xymatrix@!0{&&&& (W,D)\ar[dll]_{\nu_1} \ar[drr]^{\nu_2} \\ 
&&(Z_1,\tilde{B}_1)\ar[lldd]_{\sigma_1'}\ar[rrdd]_{\sigma_1} \ar[drr]^{\tau_1} & && &(Z_2,\tilde{B}_2)\ar[rrdd]^{\sigma_2'}\ar[lldd]^{\sigma_2} \ar[dll]_{\tau_2} &&\\ 
& & && (\hat{X},\hat{B})\ar[d]^{\eta}& &&& \\ 
(X_1,B_1) & && &(X,B)\ar@{-->}[rrrr]_{\phi_2}\ar@{-->}[llll]^{\phi_1} & & &&(X_2,B_2),}\]
where $\tau_i$ are birational morphisms, where $(Z_1,\tilde{B}_1)\stackrel{\nu_1}{\leftarrow}(W,D)\stackrel{\nu_2}{\rightarrow}(Z_2,\tilde{B}_2)$ is the minimal resolution of $(\tau_2)^{-1}\circ\tau_1$, and where each map is an isomorphism on the open part. We prove now that $\phi_2\circ (\phi_1)^{-1}$ is an isomorphism.

We first prove that either $\tau_1$ or $\tau_2$ is an isomorphism. Suppose the contrary; then, for $i=1,2$ the map $(\tau_i)^{-1}$ has a unique proper base-point $p_i\in \hat{X}$ (because so is $\sigma_i^{-1}$, Lemma~\ref{Lem:RevProperties}). Recall that $\hat{B}=A\tr F\tr C\tr E$ for some non-empty subzigzag $A\subset \hat{B}$ corresponding to the exceptional divisor of $\eta$. Furthermore $p_1,p_2$ are two distinct singular points of $\hat{B}$, belonging to the same component $D\subset \hat{B}$, which is the unique $(-1)$-curve of $A$. 
Assume that $p_1$ belongs to the component of $A$ which precedes $D$, which implies that the point $\tau_2^{-1}(p_1)\in \tilde{B}_2=(\T{E_2}\tr C_2) \tr H_2\tr (C\tr E)$ belongs to a component which precedes $C_2\subset \tau_2^{-1}(p_2)$ and thus belongs to $E_2$. Since $\nu_2$ blows-up $\tau_2^{-1}(p_1)$ and sends $C_1$ onto this point, the map $\phi_2\circ (\phi_1)^{-1}$ sends the curve $C_1$ on a point of $E_2\subset B_2$ which is a proper base-point of $\phi_1\circ (\phi_2)^{-1}$, contrary to Lemma~\ref{Lem:BasePoints}.

We may now suppose that $\tau_2$ is an isomorphism, and consider that it is the identity. We prove that so is $\tau_1$. Suppose on the contrary that $(\tau_1)^{-1}$ has a unique proper base-point $p_1 \in \hat{X}=Z_2$. Since $\sigma_1$ contracts a chain which contains only one $(-1)$-curve, the point $p_1$ belongs to the unique $(-1)$-curve of $Z_2$ contracted by $\eta=\sigma_2$, i.e. $p_1\in C_2$. Since $\phi_2(\phi_1)^{-1}$ contracts $C_1$ on $\sigma_2'(p_1)$, the point $\sigma_2'(p_1)$ is a proper base-point of $\phi_1(\phi_2)^{-1}$ and consequently belongs to $F_2$, so $p_1=F_2\cap C_2\in Z_2$. Hence the strict transform of $\T{E_2}\tr C_2$ precedes $C_1$ in $\tilde{B}_1=(\T{E_1}\tr C_1) \tr H_1\tr (C\tr E)$. The fact that $\phi_1$ is a reversion implies that no component at the left of $C_1$ in $\tilde{B}_1$ is contracted by $\sigma_1'$, hence $\phi_1(\phi_2)^{-1}$ does not contract any curve of $B_2$ except perhaps $F_2$; this implies that $\phi_1(\phi_2)^{-1}$ is a fibered modification. But $\phi_2(\phi_1)^{-1}$ cannot be a fibered modification, as it contracts $C_1$. 

The contradiction shows that $\tau_1$ is an isomorphism, hence $\tau_1(\tau_2)^{-1}(\tilde{B}_2)=\tilde{B}_1$. This means that neither $\phi_1(\phi_2)^{-1}$ nor $\phi_2(\phi_1)^{-1}$ contracts any curve, hence both maps are isomorphisms.
 \end{proof}
 If the type of the subzigzag $E$ of $B=F\tr C\tr E$ is not a palindrome, then the composition of two reversions cannot be a reversion. However, the following shows that this may occur.
 \begin{lem}[Composition of two reversions]\label{Lem:TwoReversionsGiveAreversion}
For $i=1,2$, let $\phi_i:(X,B)\dasharrow (X_i,B_i)$ be a reversion, and assume that  every irreducible curve of $B$ has self-intersection $\geq -2$. If the proper base-points of $\phi_1$ and $\phi_2$ are distinct $($respectively equal$)$ the map $\phi_{2}\circ(\phi_{1})^{-1}$ is a reversion $($respectively an isomorphism$)$.
 \end{lem}
 \begin{proof}
 Denote by $r\geq 0$ the number of components of $E$ (each one is a  (-2)-curve).
 For $i=1,2$, let $(X,B=F\tr C\tr E)\stackrel{\sigma_i}{\leftarrow}(Z_i,\tilde{B}_i)\stackrel{\sigma_i'}{\rightarrow}(X_i,B_i=F_i\tr C_i\tr E_i)$ be a minimal resolution of $\phi_i$, such that $\tilde{B}_i=(\T{E_i}\tr C_i) \tr H_i\tr (C\tr E)$. Observe that $H_i$ is the proper transform of $F$ and $F_i$ by respectively $(\sigma_i)^{-1}$ and $(\sigma_i')^{-1}$, and that $E_i$ is a chain of $r$ (-2)-curves. 
 We therefore have a commutative diagram
 \[\xymatrix@!0{&  &(Z_i,\tilde{B}_i)\ar[rrd]^{\sigma_i'}\ar[lld]_{\sigma_i}  \\ 
(X,B)\ar[rrd]_{\nu_i}\ar@{-->}[rrrr]_{\phi_i} & & &&(X_i,B_i)\ar[lld]^{\nu_i'}\\
 && (W_i,B_i'), }\] 
where $\nu_i$ and $\nu_i'$ contract the curves $E\tr C$ and $E_i\tr C_i$ respectively. Since $\nu_1$ and $\nu_2$ contract the same curves, we may assume that $\nu_1=\nu_2=\nu$ and $(W_1,B_1')=(W_2,B_2')=(W,B')$. This yields the following commutative diagram:
 \[\xymatrix@!0{&&(Z_1,\tilde{B}_1)\ar[lld]_{\sigma_1'}\ar[rrd]^{\sigma_1}  & && &(Z_2,\tilde{B}_2)\ar[rrd]^{\sigma_2'}\ar[lld]_{\sigma_2}  &&\\ 
(X_1,B_1)\ar[rrrrd]_{\nu_1'} & && &(X,B)\ar[d]^{\nu}\ar@{-->}[rrrr]_{\phi_2}\ar@{-->}[llll]^{\phi_1} & & &&(X_2,B_2),\ar[lllld]^{\nu_2'}\\
&& & &(W,B')& &  &&}\] 
where the proper base-point of $(\nu_i')^{-1}$ is equal to the image by $\nu$ of the proper base-point of $\phi_i$ (and $(\sigma_i)^{-1}$). Consequently, if these two base-points are equal then $\phi_{2}\circ(\phi_{1})^{-1}=(\nu_2')^{-1}\circ \nu_1'$ is an isomorphism, and otherwise it is a reversion. \end{proof}
\begin{rem} By definition, a reversion $\phi:\left(X,B,\bar{\pi}\right)\dashrightarrow\left(X',B',\bar{\pi}'\right)$
restricts to an isomorphism $\phi:S=X\setminus B\stackrel{\sim}{\rightarrow}S'=X'\setminus B'$
of quasi-projective surfaces, which, in contrast with the case of
fibered modifications, is never an isomorphism of $\mathbb{A}^{1}$-fibered
surfaces between $\left(S,\bar{\pi}\mid_{S}\right)$ and $\left(S',\bar{\pi}'\mid_{S}\right)$.
Indeed, it is easily seen that the rational fibrations $\bar{\pi}:X\rightarrow\mathbb{P}^{1}$
and $\bar{\pi}':X'\rightarrow\mathbb{P}^{1}$ lift to rational fibrations
with distinct general fibers on the minimal resolution $(Z,\tilde{B})$
of $\phi$. This implies that the induced $\mathbb{A}^{1}$-fibrations
$\bar{\pi}\mid_{S}$ and $\bar{\pi}'\circ\phi\mid_{S}$ on $S$ have
distinct general fibers. 

\end{rem}

\subsection{Summary on the base-points and curves contracted}\indent\\
Recall that the center of a birational map $(X,B)\dasharrow (X',B')$ is its unique proper base-point.
\begin{lem}\label{Lem:BasePointsRevFibred}
Let $\phi:(X,B=F\tr C\tr E)\dasharrow (X',B'=F'\tr C'\tr E')$ be a birational map.

$a)$  If $\phi$ is a fibered modification, it is centered at $p=F\cap C$, and $F'$ is the only irreducible component of $B'$ contracted by $\phi^{-1}$.

$b)$ If $\phi$ is a reversion, it is centered at a point $p\in F\setminus C$, and $\phi^{-1}$ contracts the curves $C'$ and $E'$ on $p$, and also contracts $F'$ on $p$ if and only if some irreducible component of $E'$ has self-intersection $\leq -3$.
\end{lem}
\begin{proof}
Follows respectively from \ref{DiscussionFiberedModification} and Lemma~\ref{Lem:RevProperties}.\end{proof}

\section{Factorization of birational maps between $1$-standard pairs }\label{Sec:ProofOfTHM}

This section is devoted to the proof of the following result.

\begin{thm}\label{Thm:Factorization}
Let $\phi:\left(X,B\right)\dashrightarrow\left(X',B'\right)$ be a
birational map between $1$-standard pairs restricting to an isomorphism
$X\setminus B\stackrel{\sim}{\rightarrow}X'\setminus B'$. If
$\phi$ is not an isomorphism then it can be decomposed into a finite
sequence \[
\phi=\phi_{n}\circ\cdots\circ\phi_{1}:\left(X,B\right)=\left(X_{0},B_{0}\right)\stackrel{\phi_{1}}{\dashrightarrow}\left(X_{1},B_{1}\right)\stackrel{\phi_{2}}{\dashrightarrow}\cdots\stackrel{\phi_{n}}{\dashrightarrow}\left(X_{n},B_{n}\right)=\left(X',B'\right)\]
of fibered modifications and reversions between $1$-standard pairs
$\left(X_{i},B_{i}\right)$, $i=1,\ldots,n$. 

Furthermore, such a factorization of minimal length is unique, which means that if  $$\phi=\phi_{n}'\circ\cdots\circ\phi_{1}':\left(X,B\right)=\left(X'_{0},B'_{0}\right)\stackrel{\phi_{1}'}{\dashrightarrow}\cdots\stackrel{\phi_{n}'}{\dashrightarrow}\left(X'_{n},B'_{n}\right)=\left(X',B'\right)$$ is another factorization,  then there exist isomorphisms of pairs $\alpha_i:(X_i,B_i)\rightarrow (X'_i,B'_i)$  for $i=1,...,n$ such that $\alpha_{i}\circ \phi_i=\phi_i'\circ \alpha_{i-1}$  for $i=2,...,n$.
\end{thm}
\begin{enavant} Let us compare Theorem~\ref{Thm:Factorization} with the existing results in the literature. Since $\phi$ restricts
to an isomorphism between $X\setminus B$ and $X'\setminus B'$, we
know that it can be factored into a sequence of smooth blow-ups and
contraction with centers on the successive boundaries. A refined description
of such factorizations, based on a careful study of base-points of
the birational maps under consideration, was obtained by V. Danilov and M. Gizatullin
\cite{Gi-Da1}. Namely, they established that one can always find a factorization
as above with the additional property that the boundaries of all intermediate
pairs consist of a certain type of zigzags called standard in
\emph{loc. cit.} Moreover, such a factorization of minimal length
is unique up to composition by automorphisms of the intermediate projective
surfaces preserving the boundaries. In general, the intermediate pairs which arise in a Danilov-Gizatullin
factorization \[
\phi:\left(X,B\right)=\left(X_{0},B_{0}\right)\stackrel{\psi_{1}}{\dashrightarrow}\left(X_{1},B_{1}\right)\stackrel{\psi_{2}}{\dashrightarrow}\ \cdots\stackrel{\psi_{s}}{\dashrightarrow}\left(X_{s},B_{s}\right)=\left(X',B'\right)\]
of $\phi$ of minimal length are not all $1$-standard. However, there
is an obvious way to concatenate these maps into a sequence of birational
maps $\phi_{j+1}:\left(X_{\alpha_{j}},B_{\alpha_{j}}\right)\dashrightarrow\left(X_{\alpha_{j+1}},B_{\alpha_{j+1}}\right)$
between all successive $1$-standard pairs $\left(X_{\alpha_{j}},B_{\alpha_{j}}\right)$
among the pairs $\left(X_{i},B_{i}\right)$ occurring in the factorization.
Theorem~\ref{Thm:Factorization} would follow provided that we show that the birational maps
obtained by this procedure are either reversions or fibered modifications.
This is the case, and the uniqueness properties actually imply that
a minimal factorization as in Theorem~\ref{Thm:Factorization} coincides with a one obtained
from a minimal Danilov-Gizatullin factorization by the above procedure. But
a proof of this fact would require to redo a careful analysis of the
base-points of the birational maps $\phi:\left(X,B\right)\dashrightarrow\left(X',B'\right)$
under consideration. So we find it simpler and more enlightening to
give a complete and self-contained proof.\end{enavant}
We proceed in two steps.
First we show in \S\ref{SubSec:ElmLink} below that every birational map $\phi:\left(X,B\right)\dashrightarrow\left(X',B'\right)$
between $1$-standard pairs restricting to an isomorphism $\phi:X\setminus B\stackrel{\sim}{\rightarrow}X'\setminus B'$
can be decomposed in an essentially unique sequence of elementary
birational maps between a certain class of pairs which strictly
contains the $1$-standard ones. Then we check in \S\ref{SubSec:Concatenation} that these elementary
birational maps can be concatenated into sequences of reversions and
fibered modifications between the $1$-standard pairs occurring in
the factorization.

\subsection{Elementary birational links between almost standard pairs}\label{SubSec:ElmLink}

\indent\newline\noindent Here we construct an enlargement of the class
of $1$-standard pairs consisting of pairs $\left(X,B\right)$
with a boundary zigzag $B$ of a more general type. We show that within
this class every birational map $\phi:\left(X,B\right)\dashrightarrow\left(X',B'\right)$
between $1$-standard pairs restricting to an isomorphism $\phi:X\setminus B\stackrel{\sim}{\rightarrow}X'\setminus B'$
can be decomposed into an essentially unique sequence of suitable
elementary birational links consisting of either smooth blow-ups or
contractions. Some of the results of this subsection are closely related to
those of in \cite{FKZ}.

\begin{defn}
A pair $\left(X,B\right)$ consisting of a normal rational projective surface
$X$ and a nonempty zigzag $B=B_{1}\tr\cdots\tr B_{s}$ (each $B_i$ being irreducible) supported
in the regular part of $X$ is called \emph{almost standard} if the
following hold :

a) There exists a unique irreducible component $B_{m}$ of $B$ with
non-negative self-intersection, called \emph{the positive curve of $B$};

b) There exists at most one irreducible component $B_{l}$ of $B$
with self-intersection $(B_l)^2=-1$. Furthermore, if it exists it is called \emph{the $(-1)$-curve of $B$} and
$B_{m}\cdot B_{l}=1$, i.e., $l=m\pm 1$.
\end{defn}
\noindent

\begin{defn}Let $(X,B)$ be an almost standard pair and let $B_m$ be the positive curve of $B$.

A birational map $\phi:\left(X,B\right)\dashrightarrow\left(X',B'\right)$
between almost standard pairs is called an \emph{elementary link}
if it consists of one of the following four operations :

I) The contraction of the $\left(-1\right)$-curve of $B$ if it exists, 

II) If $B$ contains a $\left(-1\right)$-curve $B_{m\pm1}$, the blow-up of the intersection point $p$ of $B_m$ with  $B_{m\pm1}$, immediately followed by the contraction of the
strict transform of $B_{m}$ when $(B_m)^2=0$ in $X$.

\begin{center}
\begin{pspicture}(-0.3,0.4)(2.3,1.8)
\psline(-0.25,1)(-0.15,1)
\psline(0,1)(0.8,1)
\psline(0.95,1)(1.05,1)
\rput(0,1){\textbullet}\rput(-0.05,0.7){{\scriptsize $a$}}
\rput(-0.05,1.3){\small $B_m$}
\rput(0.8,1){\textbullet}\rput(0.75,0.7){{\scriptsize $-1$}}
\rput(0.75,1.3){\small $B_{m\pm 1}$}
\psline{<-}(1.3,1)(2.1,1)
\end{pspicture}
\begin{pspicture}(-0.3,0.4)(3.5,1.8)
\psline(-0.25,1)(-0.15,1)
\psline(0,1)(1.6,1)
\psline(1.75,1)(1.85,1)
\rput(0,1){\textbullet}\rput(-0.05,0.7){{\scriptsize $a\!-\!1$}}
\rput(-0.05,1.3){\small $B_m$}
\rput(0.8,1){\textbullet}\rput(0.75,0.7){{\scriptsize $-1$}}
\rput(1.6,1){\textbullet}\rput(1.55,0.7){{\scriptsize $-2$}}
\rput(1.55,1.3){\small $B_{m\pm 1}$}
\rput(2.6,1){$\left(\begin{array}{c} \\ \\ \end{array}\right.$}\psline{->}(2.6,1)(3.4,1)
\rput(3,1.25){\scriptsize $a=0$}
\end{pspicture}
\begin{pspicture}(-0.3,0.4)(1.3,1.8)
\psline(-0.25,1)(-0.15,1)
\psline(0,1)(0.8,1)
\psline(0.95,1)(1.05,1)
\rput(0,1){\textbullet}\rput(-0.05,0.7){{\scriptsize $0$}}
\rput(0.8,1){\textbullet}\rput(0.75,0.7){{\scriptsize $-2$}}
\rput(0.75,1.3){\small $B_{m\pm 1}$}
\rput(1.2,1){$\left.\begin{array}{c} \\ \\ \end{array}\right)$}
\end{pspicture}\end{center}
III) If $B$ contains no $\left(-1\right)$-curve and if $B_{m}$ is not a boundary of $B$, the blow-up of one
of the two points $p=B_{m-1}\cap B_{m}$ or $p=B_{m}\cap B_{m+1}$,
immediately followed by the contraction of the strict transform of
$B_{m}$ when $(B_m)^2=0$ in $X$. 

\begin{center}
\begin{pspicture}(-1.1,0.4)(2.3,1.8)
\psline(-1.05,1)(-0.95,1)
\psline(-0.8,1)(0.8,1)
\psline(0.95,1)(1.05,1)
\rput(-0.8,1){\textbullet}\rput(-0.85,0.7){{\scriptsize $-b$}}
\rput(-0.85,1.3){\small $B_{m\pm 1}$}
\rput(0,1){\textbullet}\rput(-0.05,0.7){{\scriptsize $a$}}
\rput(-0.05,1.3){\small $B_m$}
\rput(0.8,1){\textbullet}\rput(0.75,0.7){{\scriptsize $-c$}}
\rput(0.75,1.3){\small $B_{m\mp 1}$}
\psline{<-}(1.3,1)(2.1,1)
\end{pspicture}
\begin{pspicture}(-1.1,0.4)(3.5,1.8)
\psline(-1.05,1)(-0.95,1)
\psline(-0.8,1)(1.6,1)
\psline(1.75,1)(1.85,1)
\rput(-0.8,1){\textbullet}\rput(-0.85,0.7){{\scriptsize $-b$}}
\rput(-0.85,1.3){\small $B_{m\pm 1}$}
\rput(0,1){\textbullet}\rput(-0.05,0.7){{\scriptsize $a\!-\!1$}}
\rput(-0.05,1.3){\small $B_m$}
\rput(0.8,1){\textbullet}\rput(0.75,0.7){{\scriptsize $-1$}}
\rput(1.6,1){\textbullet}\rput(1.55,0.7){{\scriptsize $-c\!-\!1$}}
\rput(1.55,1.3){\small $B_{m\mp 1}$}
\rput(2.6,1){$\left(\begin{array}{c} \\ \\ \end{array}\right.$}\psline{->}(2.6,1)(3.4,1)
\rput(3,1.25){\scriptsize $a=0$}
\end{pspicture}
\begin{pspicture}(-1.1,0.4)(1.3,1.8)
\psline(-1.05,1)(-0.95,1)
\psline(-0.8,1)(0.8,1)
\psline(0.95,1)(1.05,1)
\rput(-0.8,1){\textbullet}\rput(-0.85,0.7){{\scriptsize $-b\!+\!1$}}
\rput(-0.85,1.3){\small $B_{m\pm 1}$}
\rput(0,1){\textbullet}\rput(-0.05,0.7){{\scriptsize $0$}}
\rput(0.8,1){\textbullet}\rput(0.75,0.7){{\scriptsize $-c\!-\!1$}}
\rput(0.75,1.3){\small $B_{m\pm 1}$}
\rput(1.2,1){$\left.\begin{array}{c} \\ \\ \end{array}\right)$}
\end{pspicture}\end{center}

IV) If $B$ contains no $\left(-1\right)$-curve and if $B_{m}$ is a boundary of $B$, the blow-up of an arbitrary
point $p\in B_{m}$, immediately followed by the contraction of the
strict transform of $B_{m}$ when $(B_m)^2=0$ in $X$. 

As before, the elementary links of type II), III) and IV) are said to be \emph{centered at $p$.}
\end{defn}
\begin{prop}\label{Prp:DecompElementary}
Let $\phi:\left(X,B\right)\dashrightarrow\left(X',B'\right)$ be a
birational map between almost standard pairs restricting to an isomorphism
$X\setminus B\stackrel{\sim}{\rightarrow}X'\setminus B'$. Then
$\phi$ is either an isomorphism or it can be factored into a finite
sequence \[
\phi=\varphi_{r}\circ\cdots\circ\varphi_{1}:\left(X,B\right)=\left(X_{1},B_{1}\right)\stackrel{\varphi_{1}}{\dashrightarrow}\left(X_{2},B_{2}\right)\stackrel{\varphi_{2}}{\dashrightarrow}\cdots\stackrel{\varphi_{r}}{\dashrightarrow}\left(X_{r},B_{r}\right)=\left(X',B'\right)\]
of elementary links between almost standard pairs. 

\end{prop}
\begin{proof}
We proceed by induction on the total number of base-points $s\left(\phi,\phi^{-1}\right)$
of $\phi$ and $\phi^{-1}$. If $s\left(\phi,\phi^{-1}\right)=0$ then $\phi$
is an isomorphism. We assume thus that $s\left(\phi,\phi^{-1}\right)>0$,
and let $\left(X,B\right)\stackrel{\sigma}{\leftarrow}(Z,\tilde{B})\stackrel{\sigma'}{\rightarrow}\left(X',B'\right)$
be the minimal resolution of $\phi$, where the birational morphisms
$\sigma$ and $\sigma'$ consist of blow-ups of the successive base
points of $\phi$ and $\phi^{-1}$ respectively. The map $\sigma'$ contracts at most two $(-1)$-curves of $B$, namely the proper transforms
of the positive curve $B_{m}$ and of the unique possible $\left(-1\right)$-curve
$B_{l}$ of $B$ if it exists. If the proper transforms of $B_{m}$
and $B_{l}$ in $Z$ are both $\left(-1\right)$-curves then $B_{m}\cdot B_{l}=1$
in $Z$ necessarily and so, they cannot be both exceptional for $\sigma$.
Therefore $\phi$ has at most one proper base-point. In turn, this
implies that $\phi$ and all its successive lifts to the intermediate
pairs occurring in the decomposition of $\sigma$ into a sequence
of smooth blow-ups have a unique proper base-point. A similar description
holds for $\phi^{-1}$. 

If the proper base-point of $\phi$ (respectively $\phi^{-1}$) corresponds to the proper transform in $Z$ of the unique possible $\left(-1\right)$-curve of $B'$ (respectively of $B$) we factor $\phi$ by the contraction of this $(-1)$-curve; this decreases the total number of base-points of $\phi$ (respectively $\phi^{-1}$). We may thus assume that the unique proper base-points $p$ and $q$ of $\phi$ and $\phi^{-1}$ respectively correspond to the positive curve of $B'$ and $B$ respectively (if $\phi$ contracts the positive curve, it has to have a base-point, so either both have base-points or $\phi$ is an isomorphism, a case eliminated before).

We decrease $s\left(\phi,\phi^{-1}\right)$ by performing an elementary link with center at  $p\in B_m$. The existence of such a link
is clear if $B_{m}$ is a boundary of $B$ and if there is no $\left(-1\right)$-curve in $B$. Otherwise,
we distinguish two cases : 

a) If $B$ contains a $(-1)$-curve $B_l$, then $B_{m}$ intersects it and $p=B_{m}\cap B_{l}$ necessarily. Indeed, otherwise,
after the contraction of the proper transform of $B_{m}$ by $\sigma'$,
$B_{l}$ would be a curve with non-negative self-intersection and
not contracted by $\phi^{-1}$, which contradicts our assumptions.

b) If $B_{m}$ intersects two irreducible components $B_{m-1}$ and
$B_{m+1}$ of $B$ with self-intersection $\leq-2$, then $p=B_{m-1}\cap B_{m}$
or $p=B_{m}\cap B_{m+1}$. Indeed, otherwise, after the contraction
of the proper transform of $B_{m}$ by $\sigma'$, the total transform
of $B$ would not $B$ an SNC divisor, which is absurd. 

We conclude that in any of these two cases, $\phi$ can be factored
through an elementary link with center at $p$, of type II) in case
a) and of type III) in case b). This completes the proof.
\end{proof}
\begin{example}
(Factorization of fibered modifications). Let $\phi:\left(X,B=F\tr C\tr E\right)\dashrightarrow\left(X',B'\right)$
be a fibered modification between $1$-standard pairs lifting a triangular
automorphism $\phi:\left(x,y\right)\mapsto\left(ax+b,cy+P\left(x\right)\right)$
of $\mathbb{A}^2$, where $d=\deg P\geq2$ as in Lemma~\ref{Lem:LiftTriangularC2}. It follows
from the description of the resolution of such a birational map
given in \ref{DiscussionFiberedModification} that $\phi$ factors into a sequence of $d-1$ elementary
links of type IV) with centers at the intersection of the positive
curve with the proper transform of $C$, followed by a sequence of
$d-1$ elementary link of type IV) with centers at points outside the
proper transform of $C$. One easily checks that this factorization
is obtained by the lift (described in Lemma \ref{Lem:LiftTriangularC2}) of the factorization of the corresponding birational map $\phi:\left(\mathbb{F}_{1},F_{1}\tr C_{0}\right)\dashrightarrow\left(\mathbb{F}_{1},F_{1}\tr C_{0}\right)$,
which consists of a sequence of $d-1$ elementary transformations
$\sigma_{i}:\left(\mathbb{F}_{i},F_{i}\tr C_{0}\right)\dashrightarrow\left(\mathbb{F}_{i+1},F_{i+1}\tr C_{0}\right)$
with center at $F_{i}\cap C_{0}$, $i=1,\ldots,d-1$, followed by
a sequence of $d-1$ elementary transformations $\sigma_{i}':\left(\mathbb{F}_{i+1},F'_{i+1}\tr C_{0}\right)\dashrightarrow\left(\mathbb{F}_{i},F'_{i}\tr C_{0}\right)$
with center at a point of $F'_{i+1}\setminus C_{0}$, $i=d-1,\ldots,1$. 
\end{example}
\noindent

\begin{example}\label{Exa:ReversionsElLinks}
(Factorization of reversions). Let $\phi:\left(X,B=F\tr C\tr E\right)\dashrightarrow\left(X',B'\right)$
be a reversion between $1$-standard pairs, where $B$ is of type $(0,-1,-n_1,...,-n_r)$. According to Lemma~\ref{Lem:DescrRev1isaRev} and Proposition~\ref{Prop:UnicityReversions}, $\phi$ may be decomposed  as $\phi=\theta_r\varphi_r...\varphi_0\theta_0$, where $\varphi_i:(X_{i-1},B_{i-1})\dasharrow (X_i',B_i')$ and $\theta_i:(X_i',B_i')\dasharrow (X_i,B_i)$ are defined in $\S\ref{DescriptionReversion0}$ (note that $(X,B)=(X_0',B_0')$).
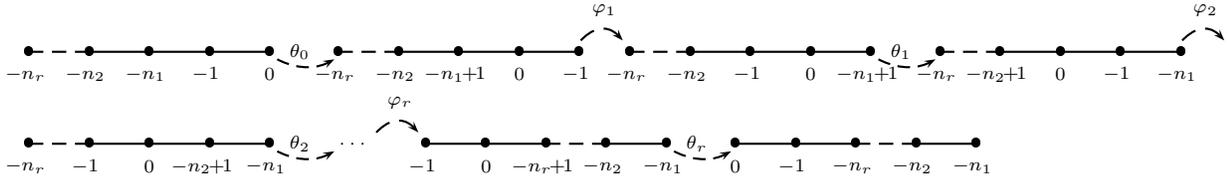
\begin{center}\begin{tabular}{l}
\begin{pspicture}(0,0.4)(4,1.5)
\psline(0.8,1)(3.2,1)
\psline[linestyle=dashed](0,1)(0.8,1)
\rput(0,1){\textbullet}\rput(-0.05,0.7){{\scriptsize $-n_r$}}
\rput(0.8,1){\textbullet}\rput(0.75,0.7){{\scriptsize $-n_2$}}
\rput(1.6,1){\textbullet}\rput(1.55,0.7){{\scriptsize $-n_1$}}
\rput(2.4,1){\textbullet}\rput(2.35,0.7){{\scriptsize $-1$}}
\rput(3.2,1){\textbullet}\rput(3.2,0.7){{\scriptsize $0$}}
\rput(3.6,1){{\scriptsize $\theta_0$}}
\parabola[linestyle=dashed]{->}(3.3,0.9)(3.68,0.8)
\end{pspicture}
\begin{pspicture}(0,0.4)(3.76,1.5)
\psline(0.8,1)(3.2,1)
\psline[linestyle=dashed](0,1)(0.8,1)
\rput(0,1){\textbullet}\rput(-0.05,0.7){{\scriptsize $-n_r$}}
\rput(0.8,1){\textbullet}\rput(0.75,0.7){{\scriptsize $-n_2$}}
\rput(1.6,1){\textbullet}\rput(1.55,0.7){{\scriptsize $-n_1\!\!+\!\!1$}}
\rput(2.4,1){\textbullet}\rput(2.4,0.7){{\scriptsize $0$}}
\rput(3.2,1){\textbullet}\rput(3.15,0.7){{\scriptsize $-1$}}
\rput(3.53,1.55){{\scriptsize $\varphi_1$}}
\parabola[linestyle=dashed]{->}(3.22,1.15)(3.5,1.3)
\end{pspicture}
\begin{pspicture}(0,0.4)(4,1.5)
\psline(0.8,1)(3.2,1)
\psline[linestyle=dashed](0,1)(0.8,1)
\rput(0,1){\textbullet}\rput(-0.05,0.7){{\scriptsize $-n_r$}}
\rput(0.8,1){\textbullet}\rput(0.75,0.7){{\scriptsize $-n_2$}}
\rput(1.6,1){\textbullet}\rput(1.55,0.7){{\scriptsize $-1$}}
\rput(2.4,1){\textbullet}\rput(2.4,0.7){{\scriptsize $0$}}
\rput(3.2,1){\textbullet}\rput(3.15,0.7){{\scriptsize $-n_1\!\!+\!\! 1$}}
\rput(3.6,1){{\scriptsize $\theta_1$}}\parabola[linestyle=dashed]{->}(3.3,0.9)(3.68,0.8)\end{pspicture}
\begin{pspicture}(0,0.4)(3.76,1.5)
\psline(0.8,1)(3.2,1)
\psline[linestyle=dashed](0,1)(0.8,1)
\rput(0,1){\textbullet}\rput(-0.05,0.7){{\scriptsize $-n_r$}}
\rput(0.8,1){\textbullet}\rput(0.75,0.7){{\scriptsize $-n_2\!\!+\!\! 1$}}
\rput(1.6,1){\textbullet}\rput(1.6,0.7){{\scriptsize $0$}}
\rput(2.4,1){\textbullet}\rput(2.35,0.7){{\scriptsize $-1$}}
\rput(3.2,1){\textbullet}\rput(3.15,0.7){{\scriptsize $-n_1$}}
\rput(3.53,1.55){{\scriptsize $\varphi_2$}}
\parabola[linestyle=dashed]{->}(3.22,1.15)(3.5,1.3)\end{pspicture}
\\
\begin{pspicture}(0,0.4)(4,1.5)
\psline(0.8,1)(3.2,1)
\psline[linestyle=dashed](0,1)(0.8,1)
\rput(0,1){\textbullet}\rput(-0.05,0.7){{\scriptsize $-n_r$}}
\rput(0.8,1){\textbullet}\rput(0.75,0.7){{\scriptsize $-1$}}
\rput(1.6,1){\textbullet}\rput(1.6,0.7){{\scriptsize $0$}}
\rput(2.4,1){\textbullet}\rput(2.3,0.7){{\scriptsize $-n_2\!\!+\!\! 1$}}
\rput(3.2,1){\textbullet}\rput(3.15,0.7){{\scriptsize $-n_1$}}
\rput(3.6,1){{\scriptsize $\theta_2$}}\parabola[linestyle=dashed]{->}(3.3,0.9)(3.68,0.8)\end{pspicture}\begin{pspicture}(2.6,0.4)(3.76,1.5)
\psline[linestyle=dotted](2.75,1)(3.05,1)
\parabola[linestyle=dashed]{->}(3.22,1.15)(3.5,1.3)
\rput(3.53,1.55){{\scriptsize $\varphi_r$}}
\end{pspicture}
\begin{pspicture}(0,0.4)(4,1.5)
\psline(0,1)(1.6,1)
\psline(2.4,1)(3.2,1)
\psline[linestyle=dashed](1.6,1)(2.4,1)
\rput(0,1){\textbullet}\rput(-0.05,0.7){{\scriptsize $-1$}}
\rput(0.8,1){\textbullet}\rput(0.8,0.7){{\scriptsize $0$}}
\rput(1.6,1){\textbullet}\rput(1.55,0.7){{\scriptsize $-n_r\!\!+\!\! 1$}}
\rput(2.4,1){\textbullet}\rput(2.35,0.7){{\scriptsize $-n_2$}}
\rput(3.2,1){\textbullet}\rput(3.15,0.7){{\scriptsize $-n_1$}}
\rput(3.6,1){{\scriptsize $\theta_r$}}\parabola[linestyle=dashed]{->}(3.3,0.9)(3.68,0.8)\end{pspicture}
\begin{pspicture}(0,0.4)(3.5,1.5)
\psline(0,1)(1.6,1)
\psline(2.4,1)(3.2,1)
\psline[linestyle=dashed](1.6,1)(2.4,1)
\rput(0,1){\textbullet}\rput(0,0.7){{\scriptsize $0$}}
\rput(0.8,1){\textbullet}\rput(0.75,0.7){{\scriptsize $-1$}}
\rput(1.6,1){\textbullet}\rput(1.55,0.7){{\scriptsize $-n_r$}}
\rput(2.4,1){\textbullet}\rput(2.35,0.7){{\scriptsize $-n_2$}}
\rput(3.2,1){\textbullet}\rput(3.15,0.7){{\scriptsize $-n_1$}}
\end{pspicture}
\end{tabular}\end{center}

If $n_i\geq 3$, then  $(X_{i-1},B_{i-1})$ and $(X_i',B_i')$  are almost-standard pairs and $\varphi_i:(X_{i-1},B_{i-1})\dasharrow(X_i',B_i')$ is the composition of a link of type II and $n_i-3$ links of type III.

If $n_i=2$, then $\varphi_i$ is an isomorphism between two pairs which are not almost-standard (there are two $(-1)$-curves in the boundary). Let  $n_i,n_i+1,...,n_i+m$ be a sequence of multiplicities equal to $2$ (with $m\geq 0$), such that either $n_{i-1}$ -- respectively  $n_{i+m+1}$ -- does not exist  ($i=1$ or $i+m=r$) or is strictly bigger than $2$. Then, the map $\theta_{i+m}\varphi_{i+m}\cdots \theta_i\varphi_i\theta_{i-1}$  is the composition of $m+2$ links of type I and $m+2$ links of type III or IV.

The remaining maps to decompose are the $\theta_i$ which are between two almost standard pairs. Then $\theta_i$ is the composition of a link of type I and a link of type III (respectively IV), if $i>0$ (respectively if $i=0$). 
\end{example}

\subsection{Concatenating elementary links into birational maps between $1$-standard
pairs}\label{SubSec:Concatenation}

\begin{enavant} Given a birational map $\phi:$$\left(X,B\right)\dashrightarrow\left(X',B'\right)$
between $1$-standard pairs, restricting to an isomorphism $X\setminus B\stackrel{\sim}{\rightarrow}X'\setminus B'$,
it follows from Proposition~\ref{Prp:DecompElementary} that there exists a factorization
\[
\phi=\varphi_{r}\circ\cdots\circ\varphi_{1}:\left(X,B\right)=\left(X_{1},B_{1}\right)\stackrel{\varphi_{1}}{\dashrightarrow}\left(X_{2},B_{2}\right)\stackrel{\varphi_{2}}{\dashrightarrow}\cdots\stackrel{\varphi_{r}}{\dashrightarrow}\left(X_{r},B_{r}\right)=\left(X',B'\right)\]
of $\phi$ into a finite sequence of elementary links between almost
standard pairs. By concatenating these elementary links into birational
maps $\phi_{j}=\varphi_{\alpha_{j+1}-1}\circ\cdots\varphi_{\alpha_{j}+1}\circ\varphi_{\alpha_{j}}:\left(X_{\alpha_{j}},B_{\alpha_{j}}\right)\dashrightarrow\left(X_{\alpha_{j+1}},B_{\alpha_{j+1}}\right)$
between all successive $1$-standard pairs $\left(X_{\alpha_{j}},B_{\alpha_{j}}\right)$
among the pairs $\left(X_{i},B_{i}\right)$ occurring in the factorization
$\phi=\varphi_{r}\circ\cdots\circ\varphi_{1}$, we obtain a new factorization
of $\phi=\phi_{n}\circ\cdots\circ\phi_{1}$ into a finite sequence
of birational maps between $1$-standard pairs. The following lemma
gives the first part of the proof of Theorem~\ref{Thm:Factorization}.

\end{enavant}

\begin{lem}\label{Lem:ExistenceFactorization}
The birational maps $\phi_{j}:\left(X_{\alpha_{j}},B_{\alpha_{j}}\right)\dashrightarrow\left(X_{\alpha_{j+1}},B_{\alpha_{j+1}}\right)$
defined above are either reversions or fibered modifications. 
\end{lem}
\begin{proof}
There is only two possible elementary links starting with a $1$-standard
pair $\left(X,B=F\tr C\tr E\right)$, namely the contraction of
the $\left(-1\right)$-curve $C$, or the blow-up of the point $F\cap C$
followed by the contraction of the proper transform of $F$. It is
enough to show that each possibility gives rise to a birational map
which is reversion in the first case and a fibered modification on
the second one.

a) If $\varphi_{1}$ is the contraction of the $\left(-1\right)$-curve
$C$ then one checks easily that the only possible sub-sequence of
elementary links occurring in the decomposition of $\phi$ before
we reach the first $1$-standard pair $\left(X_{\alpha_{1}},B_{\alpha_{1}}\right)$
coincides with the one described in Example~\ref{Exa:ReversionsElLinks} above; indeed at each step there are only two possible links, one being the inverse of the last link produced. This shows that
if $\varphi_{1}$ is the contraction of the $\left(-1\right)$-curve
$C$, then $\phi_{1}:\left(X,B\right)\dashrightarrow\left(X_{\alpha_{1}},B_{\alpha_{1}}\right)$
is a reversion. 

b) If $\varphi_{1}:\left(X,B\right)\dashrightarrow\left(X_{1},B_{1}\right)$
is the blow-up of the point $F\cap C$ followed by the contraction
of the proper transform of $F$, then the proper transform of $C$
has self-intersection $-2$, and intersects the $(0)$-curve $F_i$ produced, which is the boundary of $B_i$, for $i=1$.
 Until the self-intersection of (the proper transform of) $C$ becomes $-1$ again, the elementary links $\varphi_{i+1}:\left(X_{i},B_{i}\right)\dashrightarrow\left(X_{i+1},B_{i+1}\right)$,
$i=1,\ldots,d-1$ consist necessarily of a sequence of the blow-up of a point of the $(0)$-curve $F_i$ of $B_i$ -- having self-intersection $0$ -- followed by the contraction of the proper transform of this curve. Consequently, the map $\phi_1$ does not contract the curve $C$, which is a section of the fibration on $(X,B)$, and thus $\phi_{1}:\left(X_{\alpha_{1}},B_{\alpha_{1}}\right)\dashrightarrow\left(X_{\alpha_{2}},B_{\alpha_{2}}\right)$ is a fibered modification. 
\end{proof}
As a consequence of the descriptions, we recover  \cite[Corollary 2]{Gi-Da1} :
\begin{cor}\label{Cor:TypeRev} 
If $(X,B)$ and $(X',B')$ are two $1$-standard pairs of type $(0,-1,-n_1,...,-n_r)$ and $(0,-1,-n_1',...,-n_s')$ such that $X\setminus B\cong X'\setminus B'$, then $r=s$ and either $n_i=n_i'$ for each $i$ or $n_i=n_{r+1-i}'$ for each $i$.
\end{cor}
\begin{proof}
Denote by $\phi:X\dasharrow X'$ the birational map obtained by extension of the isomorphism.
Lemma~\ref{Lem:ExistenceFactorization} yields a decomposition of $\phi$ into fibered modifications and reversions; the fibered modifications do not change the type of the zigzag and the reversions reverse the order of the $n_i$.
\end{proof}
Now that the existence of the factorization of Theorem~\ref{Thm:Factorization} is proved, it remains to deduce the unicity.  It is a consequence of the following lemma, which completes  the proof of  the theorem.
\begin{lem}\label{Lem:UnicityFacto}
Let $\phi:\left(X,B=F\tr C\tr E\right)\dashrightarrow\left(X',B'\right)$
be a strictly birational map between $1$-standard pairs restricting
to an isomorphism $X\setminus B\stackrel{\sim}{\rightarrow}X'\setminus B'$
and let \[
\phi=\phi_{n}\circ\cdots\circ\phi_{1}:\left(X,B\right)=\left(X_{0},B_{0}\right)\stackrel{\phi_{1}}{\dashrightarrow}\left(X_{1},B_{1}\right)\stackrel{\phi_{2}}{\dashrightarrow}\cdots\stackrel{\phi_{n}}{\dashrightarrow}\left(X_{n},B_{n}\right)=\left(X',B'\right)\]
be a decomposition of $\phi$, for $n\geq 1$, satisfying that $\phi_i$ is either a reversion or a fibered modification. Then, the following are equivalent:
\begin{itemize}
\item[({\upshape 1)}]
the decomposition above is minimal $($i.e.\ there does not exist another such decomposition with less than $n$ factors$)$;
\item[({\upshape 2)}]
for any $i<n$, the centers of $(\phi_i)^{-1}$ and $\phi_{i+1}$ are distinct, and if $\phi_i$ and $\phi_{i+1}$ are reversions then $E$ contains at least one curve of self-intersection $\leq -3$.
\end{itemize}

Furthermore, if the conditions are satisfied, the following hold:
\begin{itemize}
\item[({\upshape a)}]
the map $\phi$ is not an isomorphism, and the centers of $\phi$ and $\phi_1$ $($respectively of $\phi^{-1}$ and $(\phi_n)^{-1})$ are equal;
\item[({\upshape b)}]
if  $\phi=\phi_{n}'\circ\cdots\circ\phi_{1}':\left(X,B\right)=\left(X'_{0},B'_{0}\right)\stackrel{\phi_{1}'}{\dashrightarrow}\cdots\stackrel{\phi_{n}'}{\dashrightarrow}\left(X'_{n},B'_{n}\right)=\left(X',B'\right)$ is another factorization,  there exist isomorphisms $\alpha_i:(X_i,B_i)\rightarrow (X'_i,B'_i)$  for $i=1,...,n$ such that $\alpha_{i}\circ \phi_i=\phi_i'\circ \alpha_{i-1}$  for $i=2,...,n$.
\end{itemize}
\end{lem}
\begin{proof}
For any $i$, we write $(X_i,B_i=F_i\tr C_i\tr E_i)$, and recall that the type of $E_i$ is equal to the type of $E$ or of $\T{E}$ (Corollary~\ref{Cor:TypeRev}).

We now prove the implication $(1)\Rightarrow (2)$, or in fact its contraposition. First assume that $(\phi_i)^{-1}$ and $\phi_{i+1}$ have the same proper base-point $p\in B_i$. According to Lemma~\ref{Lem:BasePointsRevFibred}, either $p=F_i\cap C_i$ and both $(\phi_i)^{-1}$ and $\phi_{i+1}$ are fibered modifications, or $p\in F_i\setminus C_i$ and $(\phi_i)^{-1}$ and $\phi_{i+1}$ are reversions. In the first case, $\phi_{i+1}\circ \phi_{i}$ is a fibered modification and in the second this map is an isomorphism; the decomposition is thus not minimal. Assume now that $\phi_i$ and $\phi_{i+1}$ are reversions and $E$ (and thus $E_i$) is a chain of $(-2)$-curves. Lemma~\ref{Lem:TwoReversionsGiveAreversion} shows that $\phi_{i+1}\circ \phi_{i}$ is either a reversion or an isomorphism, and once again the decomposition is not minimal. 

We now prove $(2)\Rightarrow ({\mathrm a})$. Since $(2)$ is symmetric, it suffices to assume $(2)$ and to prove by induction on $n$ that $\phi$ is not an isomorphism and that the center of $\phi^{-1}$ and $(\phi_n)^{-1}$ are equal. If $n=1$, this is obvious. If $n>1$, the map $\psi=\phi_{n-1}\circ...\circ \phi_1$ contracts some curve on the center $p\in B_{n-1}$ of $(\phi_{n-1})^{-1}$, by induction hypothesis. Then $p$ is not a base-point of $\phi_n$, and $\phi_n$ contracts a curve $\Gamma \subset B_{n-1}$ that contains $p$. Consequently,  the map $\phi=\phi_n\circ \psi$ contracts a curve on $\phi_n(\Gamma)=\phi_n(p)$. This point is furthermore the center of $(\phi_n)^{-1}$ (Lemma~\ref{Lem:BasePointsRevFibred}).
  
Assume now that two decompositions $\phi=\phi_n\circ...\circ \phi_1=\phi_m'\circ ...\circ \phi_1'$ satisfying $(2)$ exist. Then, the identity map factors as $\phi'_m \circ ...\circ \phi'_1 \circ (\phi_1)^{-1} \circ... \circ (\phi_n)^{-1}$. Since it is an isomorphism, condition $(2)$ is not satisfied for this decomposition. Three possibilities occur; in each one we prove that $\phi'_1 \circ (\phi_1)^{-1}$ is an isomorphism.

{\it $i)$ both $\phi'_1$ and $(\phi_1)^{-1}$ are fibered modifications.} In this case, $\phi'_1 \circ (\phi_1)^{-1}$ is either a fibered modification or an isomorphism; the first case is not possible as it yields a decomposition of the identity satisfying $(2)$.

{\it $ii)$ both $\phi'_1$ and $(\phi_1)^{-1}$ are reversions with the same center.} Proposition~\ref{Prop:UnicityReversions} shows that $\phi'_1 \circ (\phi_1)^{-1}$ is an isomorphism.

{\it $iii)$ both $\phi'_1$ and $(\phi_1)^{-1}$ are reversions with distinct centers, and each irreducible curve of $E$ has self-intersection $\geq -2$.} This case is not possible, as it implies that  $\phi'_1 \circ (\phi_1)^{-1}$ is a reversion (Lemma~\ref{Lem:TwoReversionsGiveAreversion}) and yields a decomposition of the identity satisfying $(2)$.

Denote by $\alpha_1$ the isomorphism $\phi'_1\circ (\phi_1)^{-1}$ and replace it in the decomposition of the identity written above. Writing $\psi'_2=\phi'_2\circ \alpha_1$, which is again a reversion or a fibered reversion, we find as before that $\psi'_2\circ (\phi_2)^{-1}$ is an isomorphism, that we denote by $\alpha_2$. By induction, we define $\psi'_r=\phi'_r\circ \alpha_{r-1}$ and obtain an isomorphism $\alpha_{r}=\psi'_r\circ (\phi_r)^{-1}$ for $r=2,...,m+1$. The last relation obtained is $\alpha_{m}=\mathrm{id}$, which shows that $m=n$. Choosing $\alpha_0=\mathrm{id}$ we find that $\alpha_{i}\circ \phi_1=\phi'_i\circ \alpha_{i-1}$ for $i=1,...,n$.

This proves the two remaining implications needed, which are $(2)\Rightarrow (1)$ and $(2)\Rightarrow (b)$.
\end{proof}

\section{Graphs associated to pairs and fibrations}
In this section, we associate a graph to every normal quasi-projective surface $S$ admitting a completion by a $1$-standard pair. The graph reflects the $\mathbb{A}^1$-fibrations on $S$ and the links between these.

\begin{defn}
To every normal quasi-projective surface $S$ we associate the oriented graph $\mathcal{F}_S$, defined in the following way:

a) A vertex of $\mathcal{F}_S$ is an equivalence class of $1$-standard pairs $(X,B)$ such that $X\setminus B\cong S$, where two $1$-standard pairs $(X_1,B_1,\overline{\pi_1})$, $(X_2,B_2,\overline{\pi_2})$ define the same vertex if and only if the $\mathbb{A}^{1}$-fibered surfaces $(X_1\setminus B_1,\pi_1)$ and $(X_2\setminus B_2,\pi_2)$ are isomorphic.

b) Any arrow of $\mathcal{F}_S$ is an equivalence class of reversions. If $\phi:(X,B)\dasharrow (X',B')$ is a reversion, then the class $[\phi]$ of $\phi$ is an arrow starting from the class $[(X,B)]$ of $(X,B)$ and ending at the class $[(X',B')]$ of $(X',B')$. To follow the notation of \cite{Ser}, we write $o([\phi])=[(X,B)]$ and $t([\phi])=[(X',B')]$ for respectively the \emph{origin} and \emph{target}. Two reversions $\phi_1:(X_1,B_1)\dasharrow (X_1',B_1')$ and $\phi_2:(X_2,B_2)\dasharrow (X_2',B_2')$ define the same arrow if and only if there exist two isomorphisms  $\theta:(X_1,B_1)\rightarrow (X_2,B_2)$ and  $\theta':(X_1',B_1')\rightarrow (X_2',B_2')$, such that $\phi_2\circ\theta=\theta'\circ\phi_1$.
\end{defn}
\begin{rem}
Note that, as in \cite[2.1]{Ser}, this graph is oriented, and that any arrow $a$ admits an inverse arrow  $\bar{a}$, which is the class of $\theta^{-1}$ for any $\theta$ such that $a=[\theta]$. However, contrary to the definition of \cite{Ser}, here it is possible that $a=\bar{a}$.  The factorization theorem yields  the following basic properties for the graph $\mathcal{F}_S$.\end{rem}
\begin{prop}\label{Prp:FS}Let $S$ be a normal quasi-projective surface with a non-empty graph $\mathcal{F}_S$. Then, the following hold.

$a)$ The graph $\mathcal{F}_S$ is connected.

$b)$ There is a natural bijection between the set of vertices of $\mathcal{F}_S$ and the set of equivalence classes of $\mathbb{A}^1$-fibrations on $S$ $($see Definition~$\ref{Def:A1fibredsurfaces})$.

$c)$ Assume that $(X,B)$ is a $1$-standard pair with $X\setminus B\cong S$ and that $B$ contains at least one curve of self intersection $\leq -3$. Then, the graph $\mathcal{F}_S$ is a tree if and only if $\Aut(S)$ is generated by automorphisms of $\mathbb{A}^1$-fibrations on $S$. Moreover, we have a natural exact sequence

\[1\rightarrow H \rightarrow \Aut(S) \rightarrow \Pi_1(\mathcal{F}_S)\rightarrow 1,\]
where $H$ is the $($normal$)$ subgroup of $\Aut(S)$ generated by all automorphisms of $\mathbb{A}^1$-fibrations and $\Pi_1(\mathcal{F}_S)$ is
the fundamental group of the graph $\mathcal{F}_S$.
\end{prop}
\begin{proof}The connectedness is a direct consequence of Theorem~\ref{Thm:Factorization}.
In the sequel, we fix a vertex $v$ of $\mathcal{F}_S$ and $a$ $1$-standard pair $(X_v,B_v)$ with $X_v\setminus B_v=S$ and $v=[(X,B)]$.

If $\alpha:S\rightarrow \mathbb{A}^1$ is a $\mathbb{A}^1$-fibration, then there exists a $1$-standard pair $(X,B,\overline{\pi})$ and an isomorphism $(S,\alpha)\rightarrow (X\setminus B,\pi)$ of $\mathbb{A}^1$-fibered surfaces. The isomorphism class of $(S,\alpha)$ gives the one of $(X\setminus B,\pi)$, which is equal to the vertex $[(X,B)]$. This yields $b)$.

Given any birational map $\phi:(X_v,B_v)\dasharrow (X_v,B_v)$, we use Theorem \ref{Thm:Factorization} to write $\phi=\theta_{n+1} r_n\cdots \theta_2r_1\theta_1$, where $n\geq 0$, each $r_i$ is a reversion and each $\theta_i$ is a fibered birational map between $1$-standard pairs (which may be the identity). We associate to $\phi$ the element $[r_n][r_{n-1}]\cdots [r_2][r_1]$ of the fundamental group  $\Pi_1(\mathcal{F}_S,v)$. Observe that because $B$ contains at least a curve of self-intersection $\leq -3$, the element of $\Pi_1(\mathcal{F}_S,v)$ does not depend of the choice of the decomposition (Lemma~\ref{Lem:UnicityFacto}) and the map defined is a surjective homomorphism $\nu:\Aut(S)\rightarrow \Pi_1(\mathcal{F}_S,v)$. 

Given an $\mathbb{A}^1$-fibration $\beta:S\rightarrow \mathbb{A}^1$, let $\psi:(X_v,B_v) \dasharrow (X,B,\overline{\pi})$ be a birational map of pairs such that $\psi$ restricts to an isomorphism $(S,\beta)\rightarrow (X\setminus B,\pi)$. Then, the group $\Aut(S,\beta)$ is equal to $\psi^{-1} \Aut(X\setminus B,\pi)\psi$. By construction, this group is contained in the kernel of $\nu$.

Take an element $\phi=\theta_{n+1} r_n\cdots \theta_2r_1\theta_1$ as before, and assume that $\nu(\phi)=1$. We prove by induction on the number of reversions in the decomposition (here $n$) that $\phi$ belongs to $H$. If $n=0$, $\phi$ is a fibered birational map of $(X_v,B_v)$. Otherwise,  $[r_{i+1}][r_{i}]$ vanishes in $\Pi_1(\mathcal{F}_S)$, for some $i$, which means that $r_{i+1}=\gamma\circ (r_i)^{-1}\circ \delta$ for certain isomorphisms of $1$-standard pairs $\gamma$ and $\delta$. Writing $\phi=\varphi'r_{i+1}\theta_{i+1}r_i\theta_i\varphi$, we have
 $\phi=\varphi' \gamma (r_i)^{-1}\delta \theta_{i+1}r_i\theta_i\varphi=\varphi'\gamma \theta_i \varphi (r_i\theta_i \varphi)^{-1} \delta \theta_{i+1}(r_i\theta_i \varphi)$. Since $(r_i\theta_i \varphi)^{-1} \delta \theta_{i+1}(r_i\theta_i \varphi)\in H$, we may conclude by applying induction hypothesis to $\varphi'\gamma \theta_i \varphi$.
\end{proof}

Then, we give to the graph $\mathcal{F}_S$ a natural structure of graph of groups. Before doing it in Definition~\ref{Def:GraphGroupFS}, we recall the notion of graph of groups, following \cite[4.4]{Ser}.
\begin{defn}\label{Def:GraphGroup}
Let $\mathcal{G}$ be a graph. 

$\bullet$ A \emph{graph of groups} structure on $\mathcal{G}$ is given by the choice of

\indent\indent $a)$ a group $G_v$, for any vertex $v$ of $\mathcal{G}$;

\indent\indent $b)$  a group $G_a$ and an injective morphism $\rho_a:G_a\rightarrow G_{t(a)}$, for any arrow $a$ of $\mathcal{G}$;

\indent\indent $c)$  an anti-isomorphism $\bar{}:G_a\rightarrow G_{\bar{a}}$, for any arrow $a$, such that $\bar{\bar x}=x$ for any $x\in G_a$.
 
$\bullet$ A \emph{path} in the graph of groups is a sequence $g_na_{n-1}g_{n-1}\cdots a_2g_2a_1g_1$, where $a_i$ is an arrow from $v_i$ to $v_{i+1}$ and $g_i\in G_{v_i}$. The path \emph{starts} at $v_1$ and \emph{ends} at $v_n$, and is \emph{closed} if and only if $v_1=v_{n}$.

$\bullet$ The \emph{fundamental group} of the graphs of groups \emph{at the vertex $v$} consists of closed paths starting and ending at $v$, modulo the relations $\rho_a(h)\cdot a=a\cdot \rho_{\bar{a}}(\bar{h})$ and $a\bar{a}=1$ for any arrow $a$ and any $h\in G_a$.
\end{defn}
Note that $\rho_a(g)$ is written $g^a$ in \cite{Ser}; furthermore, the two groups $G_a$ and $G_{\bar a}$ are said to be equal, which yields the same structure as our definition, but is less convenient for the following definition.
\begin{defn}\label{Def:GraphGroupFS}
Let $S$ be a normal quasi-projective surface and let $\mathcal{F}_S$ its associated graph. Then, a \emph{graph of groups structure} on $\mathcal{F}_S$ is given by the choice of

$a)$ for any vertex $v$ of $\mathcal{F}_S$, a fixed $1$-standard pair $(X_v,B_v,\overline{\pi_v})$ in the class $v$. The group $G_v$ is then equal to $\Aut(X_v\setminus B_v,\pi_v)$;

$b)$ for any arrow $a$ of $\mathcal{F}_S$,  a reversion $r_a$ in the class of $a$, which is $(X_a,B_a,\overline{\pi_a})\stackrel{r_a}{\dasharrow} (X'_a,B'_a,\overline{\pi'_a})$, and also an isomorphism  $\mu_a:(X'_a\setminus B'_a,\pi'_a)\rightarrow (X_{t(a)}\setminus B_{t(a)},\pi_{t(a)})$. The group $G_a$ is then equal to \[\{(\phi,\phi')\in \Aut(X_a,B_a)\times \Aut(X'_a,B'_a)\ |\ r_a\circ \phi=\phi'\circ r_a\},\] and $\rho_a:G_a\rightarrow G_{t(a)}$ is given by $\rho_a((\phi,\phi'))=\mu_a\circ \phi'\circ (\mu_a)^{-1}$. 

We ask further that $r_{\bar{a}}=(r_a)^{-1}$, and $\bar{}:G_a\rightarrow G_{\bar a}$ is the map $(\phi,\phi')\mapsto (\phi',\phi)$.
 
\end{defn}
\begin{enavant}
In most cases (in particular when the subzigzag $E$ of $B=F\tr C\tr E$ is not a palindrome) $a\not=\bar{a}$ for any arrow $a$ of $\mathcal{F}_S$, and it is clear that a graph of groups structure exists on $\mathcal{F}_S$. If $a=\bar{a}$ for a certain $a$, then we may choose that $(X_a,B_a,\overline{\pi_a})=(X'_a,B'_a,\overline{\pi'_a})$ and we have $(r_a)^{-1}=\lambda \circ r_a \circ \mu$ for some elements $\lambda,\mu\in \Aut(X_a,B_a)$.  Replacing $r_a$ by $\mu r_a$ we may choose that $\mu=1$. Consequently, $(r_a)^2=r_{\bar a}r_a\in \Aut(X_a,B_a)$. But, it is not clear that this one can always be chosen to be the identity. However, we will see that this property is satisfied for all the cases that we deal with in the sequel.\end{enavant}

\begin{thm}
Let $(X,B)$ be a $1$-standard pair such that at least one component of $B$ has self-intersection $\leq -3$, and let $S=X\setminus B$. 

 If $\mathcal{F}_S$ admits a structure of graph of groups, then the fundamental group of the graph of groups obtained is naturally isomorphic to $\Aut(S)$.
\end{thm}
\begin{proof}
Let us fix a graph of groups structure for $\mathcal{F}_S$, as in Definition~\ref{Def:GraphGroupFS}.


 We will work with \emph{$g$-sequences} $s=g_{n} a_n g_{n-1}\cdots a_1 g_0$, where $n\geq 1$, $g_i\in G_{v_i}$ for $i=0,...,n$ and $a_i$ is an arrow  for $i=1,...,n$, satisfying $o(a_i)=v_{i-1}, t(a_i)=v_{i}$. We write $t(s)=v_n$ and $o(s)=v_1$. There is a natural way of concatening $g$-sequences $s_1,s_2$ to $s_2s_1$ satisyfying $t(s_1)=o(s_2)$, by multiplying the last term of $s_1$ with the first of $s_2$. Then, to any $g$-sequence $s$, we can associate a birational map $\psi_s:(X_{o(s)},B_{o(s)})\dasharrow (X_{t(s)},B_{t(s)})$, by saying that $\psi_a=\mu_{a}\circ r_a\circ (\mu_{\bar a})^{-1}:(X_{o(a)},B_{o(a)})\dasharrow (X_{t(a)},B_{t(a)})$ for any arrow $a$ of $\mathcal{F}_S$, that  $\psi_g=g:(X_v,B_v)\dasharrow (X_v,B_v)$ for any $g\in G_v$, and that $\psi_{ss'}=\psi_s\circ \psi_{s'}$, for any $g$-sequences $s,s'$ with $t(s)=o(s')$.

Let us prove now that for any two vertices $v,v'$ and any birational map $\phi:(X_v,B_v)\dasharrow (X_{v'},B_{v'})$  there exists a $g$-sequence $s$ such that $\phi=\psi_s$. We decompose $\phi$ into a minimal sequence of fibered modifications and reversions, using Theorem~\ref{Thm:Factorization}, and proceed by induction on the number of reversions that occur in the decomposition. If there is no reversion, $\phi$ is a fibered birational map, thus $v=v'$ and $\phi\in G_{v}$. Otherwise, $\phi=\phi'\circ \theta_2\varphi\theta_1$, where $\theta_1$ and $\theta_2$ are fibered birational maps (which may be isomorphisms), $\varphi$ is a reversion and the decomposition of $\phi'$ involves less reversions than the one of $\phi$. Up to isomorphisms, which change the maps $\theta_1$ and $\theta_2$, we may assume that $\varphi=r_a$ for some arrow $a$ starting from $v$. Since $r_a$ is a birational map starting from $(X_a,B_a)$ and both $\theta_1$ and $(\mu_{\bar a})^{-1}$ are fibered modification or isomorphisms $(X_v,B_v)\dasharrow (X_a,B_a)$, the map $\mu_{\bar a}\theta_1$ belongs to $G_v$. We write $\phi=\phi' \theta_2(\mu_a)^{-1} \psi_a (\mu_{\bar a}\theta_1)$ and use induction hypothesis on the map $\phi' \theta_2(\mu_a)^{-1}:(X_{t(a)},B_{t(a)})\dasharrow (X_v,B_v)$ to conclude.

Let us fix a vertex $w$ of $\mathcal{F}_S$, write $S=X_w\setminus B_w$, and denote by $\Lambda$ the group of $g$-sequences $s$ such that $t(s)=o(s)=w$. The map $s\mapsto \psi_s$ yields a surjective homomorphism $\Psi:\Lambda\rightarrow \Aut(X_w\setminus B_w)=\Aut(S)$.
The fundamental group of the graph of groups at $w$ is the quotient of the group $\Lambda$ by the relations $\rho_a(h)\cdot a=a\cdot \rho_{\bar{a}}(\bar{h})$, and $a\bar{a}=1$ for any arrow $a$ and any $h\in G_a$. To prove the theorem we prove that these relations generates the kernel of $\Psi$.

Let $a$ be an arrow. Then, $\psi_a=\mu_{a}\circ r_a\circ (\mu_{\bar a})^{-1}$ and $\psi_{\bar a}=\mu_{\bar a}\circ r_{\bar a}\circ (\mu_{ a})^{-1}$. Since $r_{\bar a}=(r_a)^{-1}$, we have $\psi_{a\bar a}=\mathrm{id}$. Recall that $G_a=\{(\phi,\phi')\in \Aut(X_a,B_a)\times \Aut(X'_a,B'_a)\ |\ r_a\circ \phi=\phi'\circ r_a\}$. If $h=(\phi,\phi')\in G_a$, then $\rho_a(h)=\mu_a\circ \phi'\circ (\mu_a)^{-1}\in G_{t(a)}$ and $\rho_{\bar a}(\bar h)=\mu_{\bar a}\circ \phi\circ (\mu_{\bar a})^{-1}\in G_{o(a)}$ by definition. Consequently, the birational maps
 $\psi_{\rho_a(h)a}=\mu_a\circ \phi' \circ r_a \circ (\mu_{\bar a})^{-1}$ and $\psi_{a\rho_{\bar{a}}(\bar{h})}=\mu_a \circ r_a \circ \phi\circ (\mu_{\bar a})^{-1}$ are equal. This shows that each relation of the fundamental group is satisfied in $\Aut(S)$.

Let $s=g_{n} a_n g_{n-1}\cdots a_1 g_0\in \Lambda$ as above, and suppose that $\psi_s=\mathrm{id}$. We prove by induction on $n$ that $s$ is trivial in the fundamental group. If $n=0$, then $s=g_0\in G_v$, and $\psi_{g_0}=\mathrm{id}$ means that $g_0=1$. Assume now that $n>0$. 
We fix $\varphi_0=(\mu_{\overline{a_1}})^{-1}\circ g_0$, $\varphi_i=(\mu_{\overline{a_{i+1}}})^{-1}\circ g_i\circ \mu_{a_i}$ for $1\leq 1\leq n-1$ and $\varphi_n=g_n \circ \mu_{a_n}$. Since $\psi_{a_i}=\mu_{a_i}\circ r_{a_i}\circ (\mu_{\overline{a_i}})^{-1}$ for each $i$, $\psi_s$ decomposes as $\psi_s=\varphi_n r_{a_n}\cdots \varphi_1 r_{a_1} \varphi_0$, where each $r_{a_i}$ is a reversion and each $\varphi_i$ is a fibered birational map. Because $\psi_s$ is the identity, there are simplifications in this decomposition, which means (by Theorem~\ref{Thm:Factorization}, and because the boundary contains at least a curve of self-intersection $\leq -3$) that for some $j\in \{1,...,n-1\}$ the map $\varphi_j$ is an isomorphism of $1$-standard pairs which sends the proper base-point of $(r_{a_j})^{-1}$ on the one of $r_{a_{j+1}}$. Consequently, $(r_{a_j})^{-1}$ and $r_{a_{j+1}}\varphi_j$ are two reversions centred at the same point, so $r_{a_{j+1}}\varphi_j=\theta (r_{a_j})^{-1}$ for some isomorphism of pairs $\theta$. This means that $\overline{a_j}=[(r_{a_j})^{-1}]=[r_{a_{j+1}}]=a_{j+1}$, whence $(r_{a_j})^{-1}=r_{\overline{a_j}}=r_{a_{j+1}}$. Moreover, $(\varphi_j,\theta)\in G_{a_{j+1}}=G_{\overline{a_j}}$, so the element $h=(\theta,\varphi_j)$ belongs to $G_{a_{j}}=G_{\overline{a_{j+1}}}$. Thus we have 
$\varphi_j=(\mu_{a_{j}})^{-1}\circ g_j\circ \mu_{a_j}$, which means that 
$\rho_{a_j}(h)=\mu_{a_j} \circ \varphi_j \circ (\mu_{a_j})^{-1}=g_j$. 
Since $\overline{a_j}\rho_{a_j}(h)a_j=\rho_{\overline{a_j}}(\bar{h})$
in the fundamental group, we may replace $a_{j+1} g_j a_j$ by $\rho_{\bar{a_j}}(\bar{h}) \in G_{t(a_{j+1})}$ in the decomposition of $s$, and reduce its length. By induction, we find that $s$ is trivial in the fundamental group.
 \end{proof}

\section{Explicit examples of affine surfaces}\label{Sec:Explicit}
In this section, we apply the tools used before (especially Lemma~\ref{Lem:GoingDownToF1} and Theorem~\ref{Thm:Factorization}) to describe examples of affine surfaces. 

\subsection{Explicit form - notation}\label{Explicitform}
According to Lemma~\ref{Lem:GoingDownToF1} the resolution of singularities of any $1$-standard pair may be obtained by some blow-up of points on a fiber of $\mathbb{F}_1$. We embedd $\mathbb{F}_1$ into $\Pn\times\mathbb{P}^1$ as \[\mathbb{F}_1=\{\left((x:y:z),(s:t)\right) \subset \Pn\times\mathbb{P}^1\ |\ yt=zs\};\] the projection on the first factor yields the birational morphism $\tau:\mathbb{F}_1\rightarrow \Pn$ which is the blow-up of $(1:0:0)\in\Pn$ and the projection on the second factor yields a $\mathbb{P}^1$-bundle $\rho: \mathbb{F}_1\rightarrow\mathbb{P}^1$. We denote by $F,L\subset \Pn$ the lines with equations $z=0$ and $y=0$ respectively. We also call $F,L\subset \mathbb{F}_1$ their proper transforms on $\mathbb{F}_1$, and denote by $C\subset \mathbb{F}_1$ the exceptional curve $\tau^{-1}((1:0:0))=(1:0:0)\times\mathbb{P}^1$. The affine line $L\setminus C\subset \mathbb{F}_1$ and its image $L\setminus (1:0:0)\subset \Pn$ will be called $L_0$.

In the sequel, we associate to any $1$-standard pair $(X,B,\overline{\pi})$, its minimal resolution of singularities $\mu:\left(Y,B,\overline{\pi}\circ\mu\right)\rightarrow\left(X,B,\overline{\pi}\right)$ 
  and a birational morphism $\eta:Y\rightarrow \mathbb{F}_1$ which is the blow-up of a finite number of points. Each of these points belongs -- as proper or infinitely near point -- to the affine line $L_0=L\setminus C\subset \mathbb{P}_1$, is defined over $\kk$ but not necessarily over $\k$; however, the set of all points blown-up by $\eta$ is defined over $\k$.
  We have $D=F\tr C\tr\hat{E}$, for some (possibly reducible) curve $\hat{E}\subset Y$ contained in $\eta^{-1}(L_0)$, and where $\overline{\pi}\circ \mu=\rho\circ\eta$, as in the diagram of Lemma~\ref{Lem:GoingDownToF1}.

We fix embeddings of $\mathbb{A}^2={\rm Spec}\left(\k\left[x,y\right]\right)$ into $\mathbb{F}_1$ and $\Pn$ as $(x,y)\mapsto \left((x:y:1),(y:1)\right)$ and $(x,y)\mapsto (x:y:1)$, which give natural isomorphisms $\mathbb{A}^2\stackrel{\sim}{\rightarrow} \mathbb{F}_1\setminus (F\cup C)$ and $\mathbb{A}^2\stackrel{\sim}{\rightarrow} \Pn\setminus F$. The restriction to the line $y=0$ yields a canonical isomorphism $\mathbb{A}^1\rightarrow L_0$ which sends $\alpha\in \k$ on $(\alpha:0:1)\in L_0\subset \Pn$. The group of affine automorphisms of $\mathbb{A}^2$ -- which is the group of automorphisms that extend to automorphisms of $\Pn$ -- is denoted by $\Aff$ and the group of triangular or de Jonqui\`eres automorphisms -- automorphisms of the fibered surface $(\mathbb{A}^2,\rho|_{\mathbb{A}^2})$ -- is denoted by $\dJo$. Explicitly, we have \[\begin{array}{rcl}
\Aff&=&\left\{(x,y)\mapsto (a_1x+a_2y+a_3,b_1x+b_2y+b_3) \ |\ a_i,b_i\in \k, a_1b_2\not=a_2b_1)\right\},\\
\dJo&=&\left\{(x,y)\mapsto (ax+P(y),by+c) \ |\ a,b\in\k^{*}, c\in \k, P\in\k[y])\right\}.\end{array}\]

Two $1$-standard pairs are isomorphic (respectively induce isomorphic affine fibered surfaces) if and only their corresponding set of points blown-up are equivalent after the action of some element of $\Aff\cap\dJo$ (respectively of $\dJo$); this follows from Lemma~\ref{Lem:LiftTriangularC2} and is explained more precisely in Lemma~\ref{Lem:IsoFibTriangular} below.

\subsection{Links between $1$-standard pairs -- isomorphisms of fibrations}
Here we describe the links between $1$-standard pairs obtained from isomorphisms of affine fibered-surfaces. In general it is possible that for two non-isomorphic $1$-standard pairs $(X,B,\overline{\pi})$ and $(X',B',\overline{\pi}')$, the affine $\mathbb{A}^1$-fibered surfaces $(X\setminus B,\pi)$ and $(X'\setminus B',\pi')$ are isomorphic; the following simple result describes the situation.
\begin{lem}\label{Lem:IsoFibTriangular}
For $i=1,2$, let $(X_i,B_i,\overline{\pi_i})$ be a $1$-standard pair, with a minimal resolution of singularities $\mu_i:\left(Y_i,B_i,\overline{\pi_i}\circ\mu_i\right)\rightarrow\left(X_i,B_i,\overline{\pi_i}\right)$ and let $\eta_i:Y_i\rightarrow \mathbb{F}_1$ be a birational morphism as in $\S \ref{Explicitform}$ above. 
 Then, the following relations are equivalent:
\begin{enumerate}
\item
the $\mathbb{A}^1$-fibered surface $(X_1\setminus B_1,\pi_1)$ and $(X_2\setminus B_2,\pi_2)$ $($respectively the pairs $(X_1,B_1,\overline{\pi_1})$ and $(X_2,B_2,\overline{\pi_2}))$ are isomorphic;
\item
there exists an element of $\dJo$ $($respectively of $\dJo\cap \Aff)$ which sends the points blown-up by $\eta_1$ onto those blown-up by $\eta_2$ and sends the curves contracted by $\mu_1$ onto those contracted by $\mu_2$.
\end{enumerate}
\end{lem}
\begin{proof}
Follows directly from Lemma~\ref{Lem:LiftTriangularC2}.
\end{proof}

Recall that each point of the exceptional curve obtained by blowing-up a point $p$ on a surface is \emph{in the first neighbourhood of $p$}, and that if $q$ is in the $m$-th neighbourhood of $p$, then any point in the first neighbourhood of $q$ is in \emph{the $(m+1)$-th neighbourhood of $p$}; by convention, a point $p$ is in its $0$-th neighbourhoud.
\begin{cor}\label{Cor:AutoFibInfinite}
Let $(X,B,\overline{\pi})$ be a $1$-standard pair with $X\setminus B$ affine. The group of automorphisms of the affine fibered-space $(X\setminus B,\pi)$ contains a subgroup isomorphic to the polynomial group $(\k[x],+)$, acting algebraically on $X$ and trivially on the base of the fibration.
\end{cor}
\begin{proof}
Consider as before a minimal resolution of singularities $\mu:\left(Y,B,\overline{\pi}\circ\mu\right)\rightarrow\left(X,B,\overline{\pi}\right)$ and let $\eta:Y\rightarrow \mathbb{F}_1$ be a birational morphism as in \S \ref{Explicitform}. Denote by $m$ the maximal height of the points blown-up by $\eta$ (so that each point blown-up belongs to the $i$-th neighbourhood of a point of $\mathbb{F}_1$, for $i\leq m$). Then, the algebraic subgroup $H$ (of infinite dimension) of $\dJo$ equal to $\left\{(x,y)\mapsto (x+y^{m+1}\cdot P(y),y) \ |\ P\in\k[y])\right\}$ acts trivially on the line $y=0$ and on the set of points that belong to the $i$-th neighbourhood of this line, for $i\leq m$. Consequently, $H$ fixes any point blown-up by $\eta$ and then preserves any irreducible curve contracted by $\mu$; according to Lemma~\ref{Lem:LiftTriangularC2}, $\mu\eta^{-1}$ conjugates $H$ to a group of automorphims of the affine fibered-surface $(X\setminus B,\pi)$.
\end{proof}
\subsection{The simplest case: $1$-standard pairs of type $(0,-1)$ / the affine plane}
The simplest $1$-standard zigzag is of type $(0,-1)$. Assuming that the complement of the boundary is affine implies that it is isomorphic to $\A^2$, as the following simple result shows:
\begin{lem}\label{Lem:k10C2}
Let $(X,B,\overline{\pi})$ be a $1$-standard pair of type $(0,-1)$, and assume that $S=X\setminus B$ is affine, then $(X,B)\cong (\mathbb{F}_1,F\tr C)$, and thus $X\setminus B\cong \mathbb{A}^2$.
\end{lem}
\begin{proof}
Let $(X,B,\overline{\pi})\stackrel{\mu}{\leftarrow}(Y,B,\overline{\pi}\mu)$ be a minimal resolution of the singularities of $X$ and let $Y\stackrel{\eta}{\rightarrow}\mathbb{F}_1$ be the morphism defined in \S \ref{Explicitform}. Suppose that some point $p\in L_0\subset \mathbb{F}_1$ is blown-up by $\eta$; then $\eta^{-1}(p)$ is a connected tree of smooth rational curves of negative intersection, and at least one irreducible curve in $\eta^{-1}(p)$ is a $(-1)$-curve, which is not contracted by $\mu$ by minimality. The image of this $(-1)$-curve by $\mu$ is thus a projective curve, which does not intersect $B$ because no singular point of $X$ belongs to $B$. This contradicts the fact that $S$ is affine, so $\eta$ -- and consequently $\mu$ -- is an isomorphism.
\end{proof}
As a direct consequence of our approach, we find the following well-known results for this simple case. Recall the notation of \S \ref{Explicitform} for the natural isomorphism $\mathbb{A}^2\stackrel{\sim}{\rightarrow} \mathbb{F}_1\setminus (F\cup C)$, such that $\rho:\mathbb{F}_1\rightarrow \mathbb{P}^1$ restricts on $\mathbb{A}^2$ to the projection on the first factor.
\begin{prop}
Let $(X,B,\overline{\pi})=(\mathbb{F}_1,F\tr C,\rho)$, and let $S=\A^2=X\setminus B$, with the natural isomorphism $(x,y)\mapsto ((x:y:1),(y:1))$ $($as in $\S \ref{Explicitform})$, such that $\overline{\pi}$ extends the map $\pi:S\rightarrow \A^1$ which is $(x,y)\mapsto x$. Then:
\begin{enumerate}
\item
the automorphism $\Psi_1:(x,y)\mapsto (y,x)$ of $S$ extends to a reversion $\Psi_1:(X,B)\dasharrow (X,B)$;
\item
letting $J=\Aut(S,\pi)$, $A=<\Aut(X,B),\Psi_1>$, the group $\Aut(S)$ is the free product of $A$ and $J$ amalgamated over their intersection:
\[\Aut(S)=A\star_{A\cap J} J.\]
\item
the following equalities occur:
\[\begin{array}{rcrcl}
\Aff&=&A&=&\left\{(x,y)\mapsto (a_1x+a_2y+a_3,b_1x+b_2y+b_3) \ |\ a_i,b_i\in \k, a_1b_2\not=a_2b_1)\right\},\\
\dJo&=&J&=&\left\{(x,y)\mapsto (ax+P(y),by+c) \ |\ a,b\in\k^{*}, c\in \k, P\in\k[y])\right\};\end{array}\]
\item
there exist infinitely many $\mathbb{A}^1$-fibrations on $S$, but only one up to automorphisms of $S$;
\item
the group $\Aut(S)$ is generated by automorphisms of $\mathbb{A}^1$-fibrations.
\end{enumerate}
\end{prop}
\begin{rem}
Assertion (2) is the famous \emph{Jung's theorem}, proved from many different manners since the original proof \cite{Jun} of Jung.  We refer to \cite{Lam} and \cite{Bob} for the proofs which are the closest to our approach.
\end{rem}
\begin{proof}
Any reversion or fibered modification that starts from $(X,B)$ gives a $1$-standard pair with a zigzag of type $(-1,0)$ (Corollary~\ref{Cor:TypeRev}), and thus which is isomorphic to $(X,B)$ (Lemma \ref{Lem:k10C2}). This implies assertion~(4). 

Observe that a reversion consists of the contraction of the $(-1)$-curve of $\mathbb{F}_1$, followed by the blow-up of a point of $F$; it is therefore the lift of an automorphism of $\Pn$, which sends the point $(1:0:0)$ on another point of $F$ (the line $z=0$), which yields a new fibration. Since the map $\Psi_1:(x:y:z)\mapsto (y:x:z)$ of $\Pn$ is an example of such map, assertion (1) is clear.

The group $\Aut(X,B)$ is the lift of the group of automorphisms of $\Pn$ that fix $(1:0:0)$ and leave $F$ invariant; since $\Psi_1$ exchanges the two points $(1:0:0)$, $(0:1:0)\in F$, the group $A$ is equal to $\Aut(\Pn,F)=\Aff$. 
The equality $\dJo=\Aut(S,\pi)$ being obvious, $(3)$ follows directly.

Let us prove now that $A$ and $J$ generate the group $\Aut(S)$.
Any element $\phi\in \Aut(S)$ extends to a birational map $\phi:(X,B)\dasharrow (X,B)$, which belongs either to $\Aut(X,B)\subset A$ or factorises as $\phi_n\circ ...\circ \phi_1$, where each $\phi_i$ is a reversion or a fibered modification (Theorem \ref{Thm:Factorization}). Since each pair which occurs in this decomposition is isomorphic to $(X,B)$, we may assume that $\phi_i\in \Aut(S)$. The fibered modifications belong to $J$ and the reversions are equal to $\alpha \Psi_1 \beta$, for some $\alpha,\beta \in \Aut(X,B)$ (follows from Proposition~\ref{Prop:UnicityReversions} and from the transitivity of the action of $\Aut(X,B)$ on $L$), this yields the equality $\Aut(S)=<A,J>$. 

Since $J$ contains $\Aut(X,B)$ we also have $\Aut(S)=<\Psi_1,J>$. 
Note that elements of $A$ are products of reversions, and then are either reversions or elements of $\Aut(X,B)$ (Lemma~\ref{Lem:TwoReversionsGiveAreversion}). To prove the amalgamated product structure, we take an element $g=g_n\circ \cdots\circ g_1\circ g_0\in \Aut(S)$, where $n\geq 1$ and the $g_i$ belong alternatively to $A\setminus J$ 
or to $J\setminus A$, and prove that $g$ is not the identity.
Since both $A$ and $J$ contain $\Aut(X,B)$, elements of $A\setminus J$  are reversions and elements of $J\setminus A$ are fibered modifications. The fact that $g$ is not trivial -- and furthermore is not an automorphism -- follows from Theorem~\ref{Thm:Factorization}, or more precisely of 
Lemma~\ref{Lem:UnicityFacto}. Assertion (2) is now clear.

It remains to prove assertion (5). Since $\Psi_1$ corresponds to $(x,y)\mapsto (y,x)$, it preserves the $\mathbb{A}^1$-fibration $(x,y)\mapsto x+y$. The equality $\Aut(S)=<\Psi_1,J>$ yields the assertion. 
\end{proof}

\subsection{$1$-standard pairs of type $(0,-1,-n)$ / surfaces with equation $uv=P(w)$ in $\mathbb{A}^3$}\label{SubSec:01k}
In (\ref{etaP}) below, we construct  a $1$-standard pair of type $(0,-1,-n)$ asssociated to any polynomial $P\in\k[w]$ of degree $n$. Then, Lemma~\ref{Lem:etaPmuP} shows that any such pair is obtained by this way. Lemma~\ref{Lem:Rev01k} provides an isomorphism of the surface with the hypersurface of $\mathbb{A}^3$ with equation $uv=P(w)$. 
\begin{enavant} Surfaces defined by an equation of the form $uv=P(w)$ have been intensively studied during the last decade, with a particular focus on the classification of additive group actions on them. In particular, L. Makar-Limanov \cite{ML90} determined by careful algebraic analysis of the coordinate ring a set of generators of their automorphism group.  Every surface with equation $uv=P(w)$ admits at least two $\mathbb{A}^1$-fibration over $\mathbb{A}^1$  induced respectively by the restrictions of the projections ${\rm pr}_u$ and ${\rm pr}_v$. The latter obviously differ by the composition of the involution of the surface which exchanges $u$ and $v$. In  \cite{Dai},  D. Daigle used similar algebraic methods as L. Makar-Limanov to show every  $\mathbb{A}^1$-fibration over $\mathbb{A}^1$ on these  is of the form ${\rm pr}_u\circ \phi$, where $\phi$ is an automorphism of the surface. Here we recover these results as corollaries of the  description of birational maps between $1$-standard pairs associated with these surfaces. It follows from a general description due to V. I. Danilov and M.H. Gizatullin \cite{Gi-Da2} (see also S. Lamy \cite{Lam2} for a self-contained proof) that the automorphism group of smooth affine quadric with equation $uv=w^2-1$ admits the structure of an amalgamated product analogous to the one of the automorphism group of the plane.  In Theorem~\ref{Thm:01kAuto}, we show that this holds more generaly for every surface with equation  $uv=P(w)$. 
\end{enavant}

We keep the notation of \S \ref{Explicitform}. Recall that $L_0=L\setminus C\subset \A^2\subset \mathbb{F}_1$ is identified with $\A^1$ via the inclusion $\alpha\mapsto (\alpha,0)$ of $\A^1$ into $\A^2\subset \mathbb{F}_1$.

\begin{enavant}\label{etaP}
To any polynomial $P\in \k[x]$ of degree $n\geq 2$, we associate a birational morphism $\eta_P:Y\rightarrow \mathbb{F}_1$ which is the blow-up of $n$ points. For each root $\alpha\in \kk$ of $P$ of multiplicity $r$, the point $\alpha \in F_0(\kk)\subset \mathbb{F}_1$ is blown-up by $\eta_P$, and for $i=1,...,r-1$, the point in the $i$-th neighbourhood of $\alpha$ that belongs to the proper transform of $F$ is also blown-up. It follows from the definition of $\eta_{P}$ that it is defined over $\k$. 
In this construction,  any irreducible curve of $Y$ contracted by $\eta$ has self-intersection~$-1$ or $-2$; the curves of self-intersection~$-1$ intersect $E$ and the others do not intersect $E$; furthermore $L^2=-n$ in $Y$.

The contraction of every irreducible curve contracted by $\eta_P$ which has self-intersection $-2$ gives rise to a birational morphism $\mu_P:Y\rightarrow X$ to a $1$-standard pair $(X,B=F\tr C\tr L)$ with a zigzag of type $(0,-1,-n)$. The following figure describes the situation. In the sequel, all the figures will represent all curves and their intersections over $\kk$.

\begin{figure}[ht]
\begin{pspicture}(1,0.7)(3,1.5)
\psline(1,1)(3,1)
\rput(1,1){\textbullet}\rput(1,1.3){{\small $F$}}\rput(0.95,0.8){{\scriptsize 0}}
\rput(2,1){\textbullet}\rput(2,1.3){{\small $C$}}\rput(1.95,0.8){{\scriptsize $-1$}}
\rput(3,1){\textbullet}\rput(3,1.3){{\small $L$}}\rput(2.95,0.8){{\scriptsize $0$}}
\end{pspicture}
\begin{pspicture}(0.4,0.7)(1.6,1.5)
\psline{<-}(0.6,1)(1.4,1)
\rput(1,1.2){$\eta_P$}
\end{pspicture}
\begin{pspicture}(1,0.7)(5.5,1.5)
{\darkgray 
\psline[linecolor=darkgray](3,1)(4.5,1)
\rput(4,1){\textbullet}\rput(3.95,0.8){{\scriptsize $-1$}}
\psline[linecolor=darkgray](3,1)(4,1.9)
\psline[linecolor=darkgray](4,1.9)(4.5,1.9)
\psframe[linecolor=darkgray](4.5,0.75)(5.5,1.25)
\rput(5,1){\scriptsize  $r_l-1$}\rput(5,1.55){$\vdots$}
\rput(4,1.9){\textbullet}\rput(3.95,1.6){{\scriptsize $-1$}}
\psframe[linecolor=darkgray](4.5,1.65)(5.5,2.15)\rput(5,1.9){\scriptsize  $r_1-1$}
}
\psline(1,1)(3,1)
\rput(1,1){\textbullet}\rput(1,1.3){{\small $F$}}\rput(0.95,0.8){{\scriptsize 0}}
\rput(2,1){\textbullet}\rput(2,1.3){{\small $C$}}\rput(1.95,0.8){{\scriptsize $-1$}}
\rput(3,1){\textbullet}\rput(3,1.3){{\small $E$}}\rput(2.95,0.8){{\scriptsize $-n$}}
\end{pspicture}
\begin{pspicture}(0.4,0.7)(1.6,2)
\psline{->}(0.6,1)(1.4,1)
\rput(1,1.2){$\mu_P$}
\end{pspicture}
\begin{pspicture}(1,0.7)(4,1.5)
{\darkgray 
\psline[linecolor=darkgray](3,1)(4,1)
\rput(4,1){\textbullet}
\psline[linecolor=darkgray](3,1)(4,1.9)
\rput(4,1.9){\textbullet}\rput(4,1.55){$\vdots$}
}
\psline(1,1)(3,1)
\rput(1,1){\textbullet}\rput(1,1.3){{\small $F$}}\rput(0.95,0.8){{\scriptsize 0}}
\rput(2,1){\textbullet}\rput(2,1.3){{\small $C$}}\rput(1.95,0.8){{\scriptsize $-1$}}
\rput(3,1){\textbullet}\rput(3,1.3){{\small $E$}}\rput(2.95,0.8){{\scriptsize $-n$}}
\end{pspicture}
\caption{\small The description of the morphisms $(\mathbb{F}_1,F\tr C\tr L)\stackrel{\eta_P}{\leftarrow}(Y,B)\stackrel{\mu_P}{\rightarrow}(X,B)$, where $r_1,...,r_l$ are the multiplicities of the roots of $P$, the degree of $P$ is $n=\sum_{i=1}^{l} r_{i}$, and where a block with label $t$ consists of a zigzag of $t$ $(-2)$-curves.}
\end{figure}

Moreover, the birational morphism $\eta_P\circ (\mu_P)^{-1}:X\rightarrow \mathbb{F}_1$ is locally given by the blow-up of the ideal $(P(x),y)$ in $\A^2$.
\end{enavant}
\begin{lem}[Isomorphism classes of surfaces of type $(0,-1,-n)$]\label{Lem:etaPmuP}
Let $(X,B=F\tr C\tr E,\overline{\pi})$ be a $1$-standard pair, such that $B$ is of type $(0,-1,-n)$ $(n\geq 2)$, with a minimal resolution of singularities $\mu:\left(Y,B,\overline{\pi}\circ\mu\right)\rightarrow\left(X,B,\overline{\pi}\right)$ and let $\eta:Y\rightarrow \mathbb{F}_1$ be a birational morphism as in $\S \ref{Explicitform}$ above. Assuming that $X\setminus B$ is affine, the following hold:

\begin{enumerate}
\item
the morphisms $\eta$, $\mu$ are equal to the morphisms $\eta_P$, $\mu_P$ defined in $\ref{etaP}$, for some polynomial $P$ of degree~$n$;
\item
any $1$-standard pair $(X',B',\overline{\pi}')$ such that $(X\setminus B,\pi)\cong(X'\setminus B',\pi')$ is isomorphic to $(X,B,\pi)$.
\item
The isomorphism classes of $1$-standard pairs $(X,B)$ of type $(0,-1,-n)$ with $X\setminus B$ affine correspond to polynomials in $\k[X]$ of degree $n$ up to affine automorphisms of the line.
\end{enumerate}
\end{lem}
\begin{proof}
Since $S=X\setminus B$ is affine, only one fiber of $\overline{\pi}$ is singular and each singularity of $X$ is solved by a chain of rational curves of $Y$ (see \S \ref{underlyingSurface}). Note that $E\subset X,Y$ is the proper transform of $L\subset \mathbb{F}_1$, and has self-intersection $-n$ in $X$ and $Y$. Denote by $f\subset Y$ the unique singular fiber of $\overline{\pi}\circ\mu$; then $f$ contains $E$, which is in the boundary $B=F\tr C\tr E$ of $Y$, and $f\setminus E$ is contained in the affine part $S$. Denote by $\Gamma$ a connected component of $\overline{f\setminus E}$. Then, $\Gamma$ contains one irreducible curve $\Gamma_0$ not contracted by $\mu$ which intersect $E$, and a (possibly empty) set of connected chains of smooth rational curves, each of self-intersection $\leq -2$, contracted by $\mu$. Since $\Gamma$ is contracted by $\eta$, it contains a $(-1)$-curve, which is necessarily $\Gamma_0$, and therefore $\overline{\Gamma\setminus \Gamma_0}$ is a chain of smooth rational curves of self-intersection $-2$. This shows that each point blown-up by $\eta$ belongs -- as a proper or infinitely near point -- to $L$.

 Denote by $a_1,...,a_l\in L_0(\kk)\subset\A^2\subset \mathbb{F}_1$ the proper base-points of $\eta^{-1}$. For $i=1,...,l$, denote by $r_i\in \mathbb{N}$ the number of components of $\eta^{-1}(a_i)\subset Y$. Since $F^2=0$ in $\mathbb{F}_1$ and $F^2=E^2=-n$ in $Y$, we have $\sum_{i=1}^r r_i=n$. Then, $\eta$ and $\mu$ correspond to the morphisms $\eta_P$ and $\mu_P$ defined in \ref{etaP}, for the polynomial $P=\prod_{i=1}^{l}(x-a_i)^{r_i}\in \k[x]$. This gives the first assertion.
 
 Let us prove the remaining assertions. Denote by $\Base(\eta_P^{-1})$ the set of points blown-up by $\eta_P$, which belong to $L_0\subset \mathbb{F}_1$ as proper or infinitely near points.  Let $\alpha \in \dJo$. According to Lemma~\ref{Lem:IsoFibTriangular}, to prove $(2)$ it suffices to show that there exists $\beta \in \Aff\cap\dJo$ such that $\beta^{-1}\alpha$ fixes each point of $\Base(\eta_P^{-1})$. The map $\alpha$ restricts to an automorphism of the affine line $L_0=L\setminus C\subset\mathbb{F}_1$, which extends to an element $\beta\in \Aff\cap\dJo=\Aut(\mathbb{F}_1,F\tr C)$. Then, $\beta^{-1}\alpha$ acts trivially on $L\subset \mathbb{F}_1$ and consequently fixes $a_i$ for $i=1,...,r$; it also fixes each point of $\Base(\eta_P^{-1})$, since these points belong to the proper transform of $L$. This yields $(2)$. Assertion $(3)$ follows directly from Lemma~\ref{Lem:IsoFibTriangular}.\end{proof}

\begin{lem}[Reversions between pairs of type $(0,-1,-n)$]\label{Lem:Rev01k}Let $P,P'\in \k[x]$ be two polynomials of degree $n\geq 2$, and let 

\vspace{-0.5 cm}

\[\begin{array}{rcccl}(X,B=F\tr C\tr E,\overline{\pi})&\stackrel{\mu_P}{\longleftarrow}&(Y,B,\overline{\pi}\mu_P)&\stackrel{\eta_P}{\longrightarrow}&(\mathbb{F}_1,F\tr C\tr L)\\
(X',B'=F'\tr C'\tr E',\overline{\pi}')&\stackrel{\mu_{P'}}{\longleftarrow}&(Y',B',\overline{\pi}'\mu_{P'})&\stackrel{\eta_{P'}}{\longrightarrow}&(\mathbb{F}_1,F\tr C\tr L)\end{array}\]
 be the corresponding construction made in $\ref{etaP}$. Suppose that there exists a reversion $\phi:(X,B)\dasharrow (X',B')$ centred at $p\in F\setminus C$, with $\phi^{-1}$ centred at $p'\in F'\setminus C'$. Then, the following hold:

$(1)$ Let $a_1,...,a_l\in\kk$ and $a_1',...,a_{l'}'\in\kk$ be the roots of $P$ and $P'$ respectively. For $i=1,...,l$, let $r_i\in \mathbb{N}$ be the multiplicity of $a_i$, which is the number of components of $\eta_P^{-1}(a_i)\subset Y$; we denote by $A_i$ the component of self-intersection $-1$ and by $\mathcal{A}_i$ the union of the $r_i-1$ components of self-intersection $-2$. We also denote by $D_i$ the strict transform by $(\tau\eta_P)^{-1}$ of the line of $\Pn$ passing through $p$ and $a_i$. 
Doing the same with primes for $P'$, we get on $Y$ and $Y'$ the following dual graphs of curves:

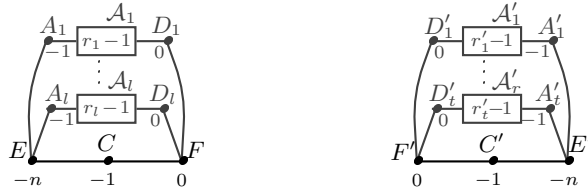
\begin{center}\begin{pspicture}(1,0.7)(4,3.1)
{\darkgray
\psline[linecolor=darkgray](1,1)(1.25,1.7)\rput(1.25,1.7){\textbullet}\rput(1.4,1.55){\scriptsize $-1$}\rput(1.35,1.9){\small $A_l$}
\psline[linecolor=darkgray,linearc=2](1,1)(0.95,1.9)(1.2,2.6)\rput(1.2,2.6){\textbullet}\rput(1.35,2.45){\scriptsize $-1$}\rput(1.3,2.8){\small $A_1$}
\psline[linecolor=darkgray](3,1)(2.75,1.7)\rput(2.75,1.7){\textbullet}\rput(2.65,1.55){\scriptsize $0$}\rput(2.75,1.9){\small $D_l$}
\psline[linecolor=darkgray,linearc=2](3,1)(3.05,1.9)(2.8,2.6)\rput(2.8,2.6){\textbullet}\rput(2.7,2.45){\scriptsize $0$}\rput(2.8,2.8){\small $D_1$}
\psline[linecolor=darkgray](1.2,2.6)(1.6,2.6)
\psline[linecolor=darkgray](2.4,2.6)(2.8,2.6)
\psframe[linecolor=darkgray](1.6,2.4)(2.4,2.8)
\psline[linecolor=darkgray](1.2,1.7)(1.6,1.7)
\psline[linecolor=darkgray](2.4,1.7)(2.8,1.7)
\psframe[linecolor=darkgray](1.6,1.5)(2.4,1.9)
\rput(2,1.7){\scriptsize $r_l\!-\!1$}\rput(2.2,2.07){\small $\mathcal{A}_l$}
\rput(1.9,2.3){$\vdots$}
\rput(2,2.6){\scriptsize $r_1\!-\!1$}\rput(2.2,2.97){\small $\mathcal{A}_1$}
}
\psline(1,1)(3,1)
\rput(0.82,1.15){{\small $E$}}
\rput(2,1.25){{\small $C$}}
\rput(3.18,1.15){{\small $F$}}
\rput(1,1){\textbullet}\rput(0.95,0.75){{\scriptsize $-n$}}
\rput(2,1){\textbullet}\rput(1.95,0.75){{\scriptsize $-1$}}
\rput(3,1){\textbullet}\rput(3,0.75){{\scriptsize $0$}}
\end{pspicture}\hspace{2cm}
\begin{pspicture}(1,0.7)(4,3.1)
{\darkgray
\psline[linecolor=darkgray](1,1)(1.25,1.7)\rput(1.25,1.7){\textbullet}\rput(1.35,1.55){\scriptsize $0$}\rput(1.35,1.9){\small $D_t'$}
\psline[linecolor=darkgray,linearc=2](1,1)(0.95,1.9)(1.2,2.6)\rput(1.2,2.6){\textbullet}\rput(1.3,2.45){\scriptsize $0$}\rput(1.3,2.8){\small $D_1'$}
\psline[linecolor=darkgray](3,1)(2.75,1.7)\rput(2.75,1.7){\textbullet}\rput(2.55,1.55){\scriptsize $-1$}\rput(2.75,1.9){\small $A_t'$}
\psline[linecolor=darkgray,linearc=2](3,1)(3.05,1.9)(2.8,2.6)\rput(2.8,2.6){\textbullet}\rput(2.6,2.45){\scriptsize $-1$}\rput(2.8,2.8){\small $A_1'$}
\psline[linecolor=darkgray](1.2,2.6)(1.6,2.6)
\psline[linecolor=darkgray](2.4,2.6)(2.8,2.6)
\psframe[linecolor=darkgray](1.6,2.4)(2.4,2.8)
\psline[linecolor=darkgray](1.2,1.7)(1.6,1.7)
\psline[linecolor=darkgray](2.4,1.7)(2.8,1.7)
\psframe[linecolor=darkgray](1.6,1.5)(2.4,1.9)
\rput(2,1.7){\scriptsize $r_t'\!\!-\!\!1$}\rput(2.2,2.07){\small $\mathcal{A}_r'$}
\rput(1.9,2.3){$\vdots$}
\rput(2,2.6){\scriptsize $r_1'\!\!-\!\!1$}\rput(2.2,2.97){\small $\mathcal{A}_1'$}
}
\psline(1,1)(3,1)
\rput(0.82,1.15){{\small $F'$}}
\rput(2,1.25){{\small $C'$}}
\rput(3.18,1.25){{\small $E'$}}
\rput(1,1){\textbullet}\rput(1,0.75){{\scriptsize $0$}}
\rput(2,1){\textbullet}\rput(1.95,0.75){{\scriptsize $-1$}}
\rput(3,1){\textbullet}\rput(2.95,0.75){{\scriptsize $-n$}}
\end{pspicture}
\end{center}

$(2)$ The numbers $l$ and $l'$ are equal, and after renumbering $r_i=r_i'$ for each $i$ and there exists an automorphism of $L_0$ which sends $a_i$ on $a_i'$ for each $i$. Moreover $(X,B)$ is isomorphic to $(X',B')$.
 
 $(3)$ 
Let $\psi=(\mu_{P'})^{-1}\phi\mu_P:Y\dasharrow Y'$ be the lift of the reversion. Then, $\psi$ restricts to an isomorphism from respectively $A_i$, $D_i$ and $\mathcal{A}_i$ to $D_i'$, $A_i'$ and $\mathcal{A}_i'$. And $\psi$ decomposes as in the following diagram
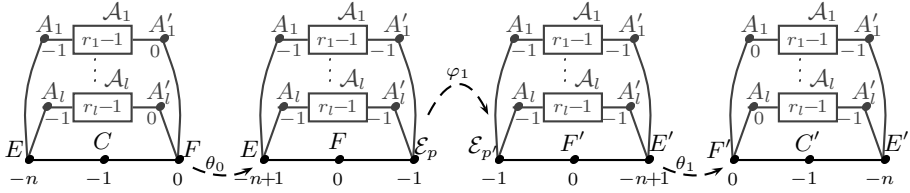
\begin{center}\begin{pspicture}(1,0.7)(4,3.1)
{\darkgray
\psline[linecolor=darkgray](1,1)(1.25,1.7)\rput(1.25,1.7){\textbullet}\rput(1.4,1.55){\scriptsize $-1$}\rput(1.35,1.9){\small $A_l$}
\psline[linecolor=darkgray,linearc=2](1,1)(0.95,1.9)(1.2,2.6)\rput(1.2,2.6){\textbullet}\rput(1.35,2.45){\scriptsize $-1$}\rput(1.3,2.8){\small $A_1$}
\psline[linecolor=darkgray](3,1)(2.75,1.7)\rput(2.75,1.7){\textbullet}\rput(2.65,1.55){\scriptsize $0$}\rput(2.75,1.9){\small $A_l'$}
\psline[linecolor=darkgray,linearc=2](3,1)(3.05,1.9)(2.8,2.6)\rput(2.8,2.6){\textbullet}\rput(2.7,2.45){\scriptsize $0$}\rput(2.8,2.8){\small $A_1'$}
\psline[linecolor=darkgray](1.2,2.6)(1.6,2.6)
\psline[linecolor=darkgray](2.4,2.6)(2.8,2.6)
\psframe[linecolor=darkgray](1.6,2.4)(2.4,2.8)
\psline[linecolor=darkgray](1.2,1.7)(1.6,1.7)
\psline[linecolor=darkgray](2.4,1.7)(2.8,1.7)
\psframe[linecolor=darkgray](1.6,1.5)(2.4,1.9)
\rput(2,1.7){\scriptsize $r_l\!\!-\!\!1$}\rput(2.2,2.07){\small $\mathcal{A}_l$}
\rput(1.9,2.3){$\vdots$}
\rput(2,2.6){\scriptsize $r_1\!\!-\!\!1$}\rput(2.2,2.97){\small $\mathcal{A}_1$}
}
\psline(1,1)(3,1)
\rput(0.82,1.15){{\small $E$}}
\rput(2,1.25){{\small $C$}}
\rput(3.18,1.15){{\small $F$}}
\rput(1,1){\textbullet}\rput(0.95,0.75){{\scriptsize $-n$}}
\rput(2,1){\textbullet}\rput(1.95,0.75){{\scriptsize $-1$}}
\rput(3,1){\textbullet}\rput(3,0.75){{\scriptsize $0$}}
\rput(3.5,0.95){{\scriptsize $\theta_0$}}
\parabola[linestyle=dashed]{->}(3.2,0.9)(3.6,0.8)
\end{pspicture}
\begin{pspicture}(1,0.7)(4,3.1)
{\darkgray
\psline[linecolor=darkgray](1,1)(1.25,1.7)\rput(1.25,1.7){\textbullet}\rput(1.4,1.55){\scriptsize $-1$}\rput(1.35,1.9){\small $A_l$}
\psline[linecolor=darkgray,linearc=2](1,1)(0.95,1.9)(1.2,2.6)\rput(1.2,2.6){\textbullet}\rput(1.35,2.45){\scriptsize $-1$}\rput(1.3,2.8){\small $A_1$}
\psline[linecolor=darkgray](3,1)(2.75,1.7)\rput(2.75,1.7){\textbullet}\rput(2.55,1.55){\scriptsize $-1$}\rput(2.75,1.9){\small $A_l'$}
\psline[linecolor=darkgray,linearc=2](3,1)(3.05,1.9)(2.8,2.6)\rput(2.8,2.6){\textbullet}\rput(2.6,2.45){\scriptsize $-1$}\rput(2.8,2.8){\small $A_1'$}
\psline[linecolor=darkgray](1.2,2.6)(1.6,2.6)
\psline[linecolor=darkgray](2.4,2.6)(2.8,2.6)
\psframe[linecolor=darkgray](1.6,2.4)(2.4,2.8)
\psline[linecolor=darkgray](1.2,1.7)(1.6,1.7)
\psline[linecolor=darkgray](2.4,1.7)(2.8,1.7)
\psframe[linecolor=darkgray](1.6,1.5)(2.4,1.9)
\rput(2,1.7){\scriptsize $r_l\!\!-\!\!1$}\rput(2.2,2.07){\small $\mathcal{A}_l$}
\rput(1.9,2.3){$\vdots$}
\rput(2,2.6){\scriptsize $r_1\!\!-\!\!1$}\rput(2.2,2.97){\small $\mathcal{A}_1$}
}
\psline(1,1)(3,1)
\rput(0.82,1.15){{\small $E$}}
\rput(2,1.25){{\small $F$}}
\rput(3.18,1.15){{\small $\mathcal{E}_p$}}
\rput(1,1){\textbullet}\rput(0.95,0.75){{\scriptsize $-n\!\!+\!\!1$}}
\rput(2,1){\textbullet}\rput(2,0.75){{\scriptsize $0$}}
\rput(3,1){\textbullet}\rput(2.95,0.75){{\scriptsize $-1$}}
\rput(3.6,2.15){{\scriptsize $\varphi_1$}}
\parabola[linestyle=dashed]{->}(3.2,1.6)(3.6,2)
\end{pspicture}
\begin{pspicture}(1,0.7)(4,3.1)
{\darkgray
\psline[linecolor=darkgray](1,1)(1.25,1.7)\rput(1.25,1.7){\textbullet}\rput(1.4,1.55){\scriptsize $-1$}\rput(1.35,1.9){\small $A_l$}
\psline[linecolor=darkgray,linearc=2](1,1)(0.95,1.9)(1.2,2.6)\rput(1.2,2.6){\textbullet}\rput(1.35,2.45){\scriptsize $-1$}\rput(1.3,2.8){\small $A_1$}
\psline[linecolor=darkgray](3,1)(2.75,1.7)\rput(2.75,1.7){\textbullet}\rput(2.55,1.55){\scriptsize $-1$}\rput(2.75,1.9){\small $A_l'$}
\psline[linecolor=darkgray,linearc=2](3,1)(3.05,1.9)(2.8,2.6)\rput(2.8,2.6){\textbullet}\rput(2.6,2.45){\scriptsize $-1$}\rput(2.8,2.8){\small $A_1'$}
\psline[linecolor=darkgray](1.2,2.6)(1.6,2.6)
\psline[linecolor=darkgray](2.4,2.6)(2.8,2.6)
\psframe[linecolor=darkgray](1.6,2.4)(2.4,2.8)
\psline[linecolor=darkgray](1.2,1.7)(1.6,1.7)
\psline[linecolor=darkgray](2.4,1.7)(2.8,1.7)
\psframe[linecolor=darkgray](1.6,1.5)(2.4,1.9)
\rput(2,1.7){\scriptsize $r_l\!\!-\!\!1$}\rput(2.2,2.07){\small $\mathcal{A}_l$}
\rput(1.9,2.3){$\vdots$}
\rput(2,2.6){\scriptsize $r_1\!\!-\!\!1$}\rput(2.2,2.97){\small $\mathcal{A}_1$}
}
\psline(1,1)(3,1)
\rput(0.82,1.15){{\small $\mathcal{E}_{p'}$}}
\rput(2,1.25){{\small $F'$}}
\rput(3.18,1.25){{\small $E'$}}
\rput(1,1){\textbullet}\rput(0.95,0.75){{\scriptsize $-1$}}
\rput(2,1){\textbullet}\rput(2,0.75){{\scriptsize $0$}}
\rput(3,1){\textbullet}\rput(2.95,0.75){{\scriptsize $-n\!\!+\!\!1$}}
\rput(3.5,0.95){{\scriptsize $\theta_1$}}\parabola[linestyle=dashed]{->}(3.2,0.9)(3.6,0.8)
\end{pspicture}
\begin{pspicture}(1,0.7)(4,3.1)
{\darkgray
\psline[linecolor=darkgray](1,1)(1.25,1.7)\rput(1.25,1.7){\textbullet}\rput(1.35,1.55){\scriptsize $0$}\rput(1.35,1.9){\small $A_l$}
\psline[linecolor=darkgray,linearc=2](1,1)(0.95,1.9)(1.2,2.6)\rput(1.2,2.6){\textbullet}\rput(1.3,2.45){\scriptsize $0$}\rput(1.3,2.8){\small $A_1$}
\psline[linecolor=darkgray](3,1)(2.75,1.7)\rput(2.75,1.7){\textbullet}\rput(2.55,1.55){\scriptsize $-1$}\rput(2.75,1.9){\small $A_l'$}
\psline[linecolor=darkgray,linearc=2](3,1)(3.05,1.9)(2.8,2.6)\rput(2.8,2.6){\textbullet}\rput(2.6,2.45){\scriptsize $-1$}\rput(2.8,2.8){\small $A_1'$}
\psline[linecolor=darkgray](1.2,2.6)(1.6,2.6)
\psline[linecolor=darkgray](2.4,2.6)(2.8,2.6)
\psframe[linecolor=darkgray](1.6,2.4)(2.4,2.8)
\psline[linecolor=darkgray](1.2,1.7)(1.6,1.7)
\psline[linecolor=darkgray](2.4,1.7)(2.8,1.7)
\psframe[linecolor=darkgray](1.6,1.5)(2.4,1.9)
\rput(2,1.7){\scriptsize $r_l\!\!-\!\!1$}\rput(2.2,2.07){\small $\mathcal{A}_l$}
\rput(1.9,2.3){$\vdots$}
\rput(2,2.6){\scriptsize $r_1\!\!-\!\!1$}\rput(2.2,2.97){\small $\mathcal{A}_1$}
}
\psline(1,1)(3,1)
\rput(0.82,1.15){{\small $F'$}}
\rput(2,1.25){{\small $C'$}}
\rput(3.18,1.25){{\small $E'$}}
\rput(1,1){\textbullet}\rput(1,0.75){{\scriptsize $0$}}
\rput(2,1){\textbullet}\rput(1.95,0.75){{\scriptsize $-1$}}
\rput(3,1){\textbullet}\rput(2.95,0.75){{\scriptsize $-n$}}
\end{pspicture}
\end{center}
where $\theta_0,\varphi_1$ and $\theta_1$ correspond to the maps described in $\S\ref{DescriptionReversion1}$, and $\mathcal{E}_p$ and $\mathcal{E}_{p'}$ correspond to the exceptional curve contracted on $p$ and $p'$ respectively, which are the proper transforms of respectively $E'$ and $E$.

$(4)$ We have the following commutative diagram of birational maps
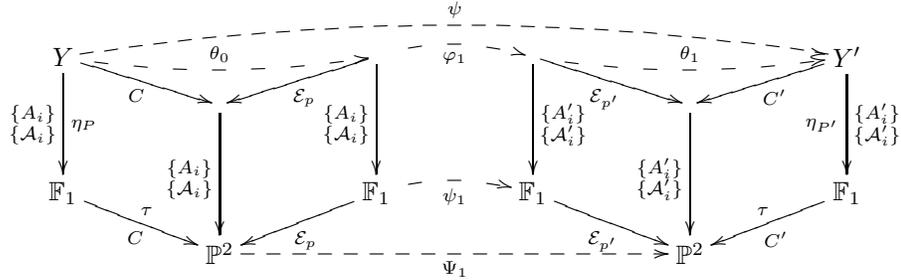
\[\xymatrix@R=3mm@C=1.5cm{Y\ar@/^1pc/@{-->}[rrrrr]^{\psi}\ar@/_0.4pc/@{-->}[rr]^{\theta_0}
\ar[ddd]^{\eta_{P}}_{\scriptsize \begin{array}{c}
\{A_i\}\\ \{\mathcal{A}_i\}\end{array}\!\!\!\!}\ar[rd]_{\scriptsize C}
& &\ar@/^0.4pc/@{-->}[r]_{\varphi_1}
\ar[ddd]_{\scriptsize \begin{array}{c}
\{A_i\}\\ \{\mathcal{A}_i\}\end{array}\!\!\!\!}
\ar[ld]^{\scriptsize \mathcal{E}_p}
& \ar[rd]_{\scriptsize \mathcal{E}_{p'}} \ar@/_0.4pc/@{-->}[rr]^{\theta_1}
\ar[ddd]^{\!\!\!\!\scriptsize \begin{array}{c}
\{A_i'\}\\ \{\mathcal{A}_i'\}\end{array}}
& & Y'\ar[ld]^{\scriptsize C'}\ar[ddd]_{\eta_{P'}}^{\!\!\!\!\scriptsize \begin{array}{c}
\{A_i'\}\\ \{\mathcal{A}_i'\}\end{array}}\\
&\ar[ddd]_{\scriptsize \begin{array}{c}
\{A_i\}\\ \{\mathcal{A}_i\}\end{array}\!\!\!\!}&&&\ar[ddd]_{\scriptsize \begin{array}{c}
\{A_i'\}\\ \{\mathcal{A}_i'\}\end{array}\!\!\!\!}&\\
&&&&&\\
\mathbb{F}_1\ar[rd]_{\scriptsize C}^{\tau}&& \mathbb{F}_1\ar@/^0.4pc/@{-->}[r]_{\psi_1}\ar[ld]^{\scriptsize \mathcal{E}_p}
& \mathbb{F}_1\ar[rd]_{\scriptsize \mathcal{E}_{p'}}&&\mathbb{F}_1\ar[ld]^{\scriptsize C'}_{\tau}&\\
&\Pn\ar@{-->}[rrr]_{\Psi_1}&&&\Pn
}\]
where each straight line is a birational morphism which contracts the curves $($or proper transform of curves$)$ written above the arrow, where $\psi_1$ is a birational map which preserves the ruling of $\mathbb{F}_1$, and where $\Psi_1$ is given by 
\[\Psi_1:(x:y:z)\dasharrow (xyz^{n-2}:G(x,z):yz^{n-1}),\]
up to automorphisms of $(Y,B)$ and $(Y',B')$.

$(5)$
letting $\chi:\Pn\dasharrow (\mathbb{P}^1)^3$ be the rational map \[\chi:(x:y:z)\dasharrow \left( (y:z), (G(x,z):yz^{n-1}), (x:z)\right),\]
the map $\chi \circ \tau\eta_P(\mu_P)^{-1}$ restricts to an embedding of $X\setminus B$ to the hypersurface of \[\mathbb{A}^{3}=\left\{(u:1),(v:1),(w:1)\in (\mathbb{P}^1)^3, u,v,w\in \k\right\}\] given by \[\left\{(u,v,w)\in \mathbb{A}^{3}\ |\ uv=P(w)\right\}.\] 
The restrictions of the three canonical projections $\mathbb{A}^3\rightarrow\mathbb{A}^1$ give respectively the $\mathbb{A}^1$-fibration $\pi$, the $\mathbb{A}^1$-fibration $\pi\circ (\Psi_1|_{X\setminus B})$ obtained by means of the reversion $\psi_0$, and the $\mathbb{A}^1_{*}$-fibration given by the pencil of lines of $\Pn$ passing through the point where the reversion $\Psi_1$ is centered.
\end{lem}
\begin{proof}Each curve $D_i$ has self-intersection $0$ in $Y$, and intersects $B$ transversally and only at $p$. Moreover, it intersects also the curve corresponding to the blow-up of $a_i$; this latter curve is $A_i$ if $r_i=1$ and the component at the right side of $\mathcal{A}_i$ otherwise. Since $D_i$ does not intersect any other curve $A_j$ or in $\mathcal{A}_j$, we obtain $(1)$.

Decompose $\psi$ into $\psi=\theta_1\varphi_1\theta_0$, as in $\S\ref{DescriptionReversion1}$, and compute the diagram of $(3)$ above for the $A_i,D_i,\mathcal{A}_i$. This shows that $D_i$ is sent by $\psi$ on a curve of self-intersection $-1$, intersecting $B'$ only at $E'$, and transversally. In consequence, any $D_i$ is sent on a curve $A_j'$. Moreover, if $r_i>0$ the singular point of $\mu_P(D_i)$ is sent by $\phi$ onto the singular point of $\mu_P(A_j')$. Since $\sum_{t=1}^l r_t=\sum_{t=1}^{l'} r_t'=n$, we obtain the equality $l=l'$ and get that after renumbering $D_i$ is sent on $A_i'$ and $\mathcal{A}_i$ is sent on $\mathcal{A}_i'$, whence $r_i=r_i'$. Assertion $(3)$ is now proved.

The projection by $\eta_P(p)\in \p^2$ on the line $L=\eta_P(E)$ induces an isomorphism $\mathcal{E}_p\to E$ which sends $\mathcal{E}_p\cap D_i$ onto $a_i$, and sends $\mathcal{E}_p\cap F$ onto $E\cap C$. Using the diagram of $(3)$, the map $\theta_1\varphi_1$ restricts to an isomorphism $\mathcal{E}_p\to E'$ which sends $\mathcal{E}_p\cap D_i=$ onto $E'\cap A_i'$ and sends $\mathcal{E}_p\cap F$ onto $E'\cap C'$. Combining the two isomorphisms, and since  $E'$ is the proper transform of $L$ by $\eta_{P'}$, we obtain $(2)$.

In the decomposition of $\psi$ given in $(3)$, $\theta_0$ decomposes as the contraction of $C$, followed by the blow-up of $p$. Moreover, each of these two steps do not change the self-intersection of any of the components of $\{A_i\}_{i=1}^l,\{\mathcal{A}_i\}_{i=1}^l$, contracted by $\eta_P$, which are thus still contractible in the surfaces obtained from $Y$ by contracting $C$ and blowing-up $p$. Doing the same with $(\theta_1)^{-1}$, we obtain the diagram of $(4)$. 
Since any curve contracted by the map $\phi_1\colon \mathbb{F}_1\dasharrow \mathbb{F}_1$ is a fibre of the ruling, $\phi_1$ preserves the ruling. 
The lift of the group of automorphisms of $\p^2$ of the form $(x:y:z)\to (x+\lambda z:y:z)$, $\lambda\in \k$ gives a group of automorphisms of $(Y,B)$ or $(X,B)$ which acts transitively on the $\k$-points of $F$. Thus, we may assume, up to automorphisms of $(Y,B)$ and $(Y',B')$, that $p=p'=(0:1:0)$. 
 It remains to observe that $\Psi_1$ can be given in this case by the map $\Psi_0:(x:y:z)\dasharrow (xyz^{n-2}:G(x,z):yz^{n-1})$. The map $\Psi_0$ clearly preserves the lines passing through $(0:1:0)$, and this point is a base-point of $\Psi_0$ of multiplicity $n-1$. One can moreover check that it has $2n-2$ other base-points defined as follows. a) The base-points which corresponds to the $n-2$ base-points of $\varphi_1$, all infinitely near of $p=(0:1:0)$ and lying on $L$, b) the $n$ points blown-up by $\eta_P$, which are $\{(a_i:0:1)\}_{i=1}^l$ and points infinitely near, all on $F$. Thus, $\Psi_1$ and $\Psi_0$ have the same base-points, and $(4)$ is now proved.
 
Letting $\zeta=\tau\eta_P(\mu_P)^{-1}:X\dasharrow \Pn$, we prove now that $\chi\circ \zeta$ restricts to an embedding of $S=X\setminus B$ into $\mathbb{A}^{3}$. The first coordinate of $(\mathbb{P}^1)^3$ corresponds to the projection of $\Pn$ by $(1:0:0)$ and then restricts exactly to $\pi:S\rightarrow\mathbb{A}^1$; the second coordinate is obtained by means of the reversion, it restricts to the $\mathbb{A}^{1}$-fibration $(\overline{\pi}\Psi_1)|_{X\setminus B}$; the last one corresponds to the projection of $\Pn$ by $(0:1:0)$ and restricts to a $\mathbb{A}^1_{*}$-fibration on $S$. This last map separates the points of the different regular fibers of $\pi$ and separates the components of the reduced fiber. Since each of these components is a section of the $\mathbb{A}^{1}$-fibration $(\overline{\pi}\Psi_1)|_{X\setminus B}$, the map $(\chi\circ \zeta)|_{X\setminus B}$ is an embedding of $S$ into $\mathbb{A}^{3}\subset (\mathbb{P}^1)^3$. Taking coordinates $(u:1),(v:1),(w:1)$ on $\mathbb{A}^{3}\subset (\mathbb{P}^1)^3$, we deduce from the explicit form of $\chi$ that the image is the surface with equation $uv=P(w)$. The reversion $\Psi_1$ computed before corresponds to the automorphism $(u,v,w)\mapsto (v,u,w)$ of the surface, and the fibration $\pi$ is the projection on the first factor.
\end{proof}

\begin{thm}\label{Thm:01kAuto}
Let $(X,B=F\tr C\tr E,\overline{\pi})$ be a $1$-standard pair, such that $B$ is of type $(0,-1,-n)$ $(n\geq 2)$ and such that the surface $S=X\setminus B$ is \emph{affine}. 

Then, there exists an isomorphism of fibered-surfaces from $(S,\pi)$ to the hypersurface of $\mathbb{A}^{3}$ given by \[\left\{(u,v,w)\in \mathbb{A}^{3}\ |\ uv=P(w)\right\},\] for some polynomial $P$ of degree $n$, equipped by the $u$-fibration; and any such surface is obtained in this way. Furthermore,  the following assertions hold:
\begin{enumerate}
\item
the isomorphism class of the surface is given by the polynomial $P$, up to a multiple and up to an automorphism of $\mathbb{A}^{1}={\rm Spec}\left(\k\left[w\right]\right)$;
\item
there exist infinitely many fibrations on $S$, but only one up to an automorphism of $S$;
\item
the graph $\mathcal{F}_S$ is \hspace{0.5cm} $\xymatrix@R=1mm@C=2cm{
\ar@(ul,dl)\bullet
}$;
\item
if $n\geq 3$, the group $\Aut(S)$ is not generated by the automorphisms of $\mathbb{A}^1$-fibration;
\item
if $n=2$, the group $\Aut(S)$ is generated by the automorphisms of $\mathbb{A}^1$-fibration;
\item
the involution $(u,v,w)\mapsto (v,u,w)$ on $S$ corresponds to a reversion $\Psi_1:(X,B)\dasharrow (X,B)$;
\item
the group 
 $\Aut(S)$ is the free product of $A=<\Aut(X,B),\Psi_1>$ and $J=\Aut(S,\pi)$, amalgamated over their intersection $A\cap J=\Aut(X,B)$:
\[\Aut(S)=A\star_{A\cap J} J;\]
\begin{enumerate}
\item
if $n=2$, the contraction of $C\tr E$ gives a birational morphism of pairs $(X,B)\rightarrow (Z,D)$, which conjugates $A$ to the group $\Aut(Z,D)$. Moreover if $P$ has two distinct roots in $\kk$, $Z$ is a smooth quadric in $\mathbb{P}^3$, and $D$ is an hyperplane section; if one adds that the two roots are defined over $\k$, then $Z$ is isomorphic to $\mathbb{P}^1\times\mathbb{P}^1$ and $D$ becomes a diagonal. If $P$ has only one multiple root (necessarily defined over $\k$),  $Z$ is the weighted plane $\mathbb{P}(1,1,2)$ obtained by contracting the $(-2)$-curve of $\mathbb{F}_2$ and $D$ is the image of a section of self-intersection $2$.
\item
 if $n>3$, we denote by $\Aut(X,B,\Psi_1)$ the subgroup of $\Aut(X,B)$ which fixes the unique proper base-point of $\Psi_1$; then $A$ is the free product of $\Aut(X,B)$ and $A_0=<\Aut(X,B,\Psi_1),\Psi_1>$, amalgamated over their intersection $A_0\cap \Aut(X,B)=\Aut(X,B,\Psi_1)$, and we also have another amalgamated product structure for $\Aut(S)$:
\[\Aut(S)=A_0 \star_{A_0\cap J} J.\]
 \end{enumerate}
 \item
denoting by $H,I,T,T_0,Sp$ the following subgroups of automorphisms of $S$:
\[\begin{array}{lcl}
H&=&\{(u,v,w)\mapsto (au,a^{-1}v,w)\ |\ a\in\k^{*}\};\\
I&=&\{(u,v,w)\mapsto (v,u,w), (u,v,w)\mapsto (u,v,w)\};\\
T&=&\{(u,v,w)\mapsto (u,u^{-1}\cdot (P(w+uq(u))-P(w)),w+uq(u))\ | \ q\in \k[u]\};\\
T_0&=&\{(u,v,w)\mapsto (u,u^{-1}\cdot (P(w+au)-P(w)),w+au)\ | \ a\in \k\};\\
Sp&=&\{(u,v,w)\mapsto (u,cv,aw+b)\ | \ a,c\in \k^{*}, b\in \k, P(aw+b)=c P(w)\};
\end{array}\]
then, $H\cong \k^{*}$, $I\cong \mathbb{Z}/2\mathbb{Z}$, $T\cong \k[u]$, $T\supset T_0\cong \k$ and $Sp\cong\left\{\begin{array}{rl}\k^{*}&\mbox{if } p \mbox{ has only one root,}\\ \mathbb{Z}/m\mathbb{Z}&\mbox{otherwise {\upshape (}$m=1$ in general{\upshape )}.}\end{array}\right.$ 

Furthermore, the following occur:
\begin{enumerate}
\item
 $J=<H,Sp,T>\cong \k[u]\rtimes (\k^{*}\times Sp)$ is the group of automorphisms of $(S,\pi)$.
\item
$\Aut(X,B)=<H,Sp,T_0>\cong \k\rtimes (\k^{*}\times Sp)$;
\item

$\Aut(X,B,\Psi_1)=<H,Sp>\cong (\k^{*}\times Sp)$;
\item

$A_0=<H,Sp,I>\cong (\k^{*}\rtimes \mathbb{Z}/{2}\mathbb{Z})\times Sp$;
\item
$A=<H,Sp,I,T_0>\cong \left\{\begin{array}{ll}\k\rtimes (\k^{*}\times Sp)&\mbox{if } n=2,\\ \Aut(X,B)\star_{\Aut(X.B)\cap A_0}A_0&\mbox{otherwise.}\end{array}\right.$
\end{enumerate}
\end{enumerate} 
\end{thm}
\begin{proof}
According to Lemma \ref{Lem:etaPmuP}, there exist two morphisms $(X,B=F\tr C\tr E,\overline{\pi})\stackrel{\mu_P}{\leftarrow}(Y,B,\overline{\pi}\mu_P)\stackrel{\eta_P}{\rightarrow}(\mathbb{F}_1,F\tr C\tr L)$ as in \ref{etaP}, for some polynomial $P$ of degree $n$. 
 The isomorphism between $S$ and the surface $\{(u,v,w)\in \A^3 \ |\ uv=P(w)\}$ follows from Lemma~\ref{Lem:Rev01k}. 
 The isomorphism class of the pair $(X,B)$ is determined by the points blown-up by $\eta_P$, up to an action of $\Aff\cap \dJo$ (Lemma~\ref{Lem:Rev01k}), and consequently by the polynomial $P$ up to multiple and to an automorphism of $L_0=\mathbb{A}^1$ (follows from the description of $\eta_P$ and $\mu_P$, made in \ref{etaP}).
 
Since any reversion or fibered modification that starts from $(X,B)$ yields an isomorphic $1$-standard pair (Lemma~\ref{Lem:Rev01k}), any $1$-standard pair  $(X',B')$ such that $X\setminus B\cong X'\setminus B'$ is isomorphic to $(X,B)$ (Theorem~\ref{Thm:Factorization}). This implies -- with the discussion made above -- the assertions (1) and (2); it also shows that the graph $\mathcal{F}_S$ contains only one vertex; we prove now that it contains only one arrow. The group of automorphism of $\Pn$ that fix each point of $L$, and preserve the line $F$ lift to a subgroup of $\Aut(X,B)$ which acts transitively on $F\setminus C$. Consequently, if $\phi:(X,B)\dasharrow (X',B')$ is a reversion, there exists $\alpha \in \Aut(X,B)$ such that $\phi\alpha$ is a reversion centered at the same point as $\Psi_1$. Proposition~\ref{Prop:UnicityReversions} implies that $\Psi_1=\beta\phi\alpha$, for some isomorphism $\beta:(X',B')\rightarrow (X,B)$. This yields assertions (3) and thus (4) (using Proposition~\ref{Prp:FS}).

Let us prove assertion $(5)$. Assume that $n=2$, and let $\alpha$ be an element of $\Aut(X,B)$ which does not fix the proper base-point of $\Psi_1$. The reversions $\Psi_1^{-1}=\Psi_1$ and $\Psi_1\alpha$ have thus distinct base-points, so $\Psi_1\alpha\Psi_1$ is a reversion (Lemma~\ref{Lem:TwoReversionsGiveAreversion}), equal to $\beta \Psi_1\gamma$, for some $\beta,\gamma\in \Aut(X,B)$. Consequently $\Psi_1=(\alpha^{-1}\Psi_1)(\beta\Psi_1\gamma)=\alpha^{-1}(\Psi_1\beta\Psi_1^{-1})\gamma$; since $\Psi_1\beta\Psi_1^{-1}$ preserves the fibration $\Psi_1\pi$, the reversion $\Psi_1$ is generated by automorphisms of $\mathbb{A}^1$-fibrations. The equality $\Aut(S)=<\Aut(X,B),\Psi_1,J>$ yields assertion $(5)$. 

Assertion~$(6)$ follows from Lemma~\ref{Lem:Rev01k}. It remains to prove the main assertions, i.e.\ $(7)$ and $(8)$.


Let us write $I=<\Psi_1>$ and $J=\Aut(S,\pi)$ (automorphisms of $S$ which preserve the fibration $\pi$). We prove now that $\Aut(X,B),I,J$ generate $\Aut(S)$. Any element $g\in \Aut(S)$ extends to a birational map $g:(X,B)\dasharrow (X,B)$; either $g$ belongs to $\Aut(X,B)$ or it may be written -- using Theorem \ref{Thm:Factorization}) -- as 
\[
 g= g_{n}\circ\cdots\circ g_{1}:\left(X,B\right)=\left(X_{0},B_{0}\right)\stackrel{ g_{1}}{\dashrightarrow}\left(X_{1},B_{1}\right)\stackrel{ g_{2}}{\dashrightarrow}\cdots\stackrel{ g_{n}}{\dashrightarrow}\left(X_{n},B_{n}\right)=\left(X,B\right)\]
 where $g_i$ is a reversion or a fibered modification. We proved previously that each $(X_i,B_i)$ is isomorphic to $(X,B)$, we may thus assume, by changing the $g_i$, that $(X_i,B_i)=(X,B)$. Consequently, $g_i$ may be viewed as an element of $\Aut(S)$. If it is a fibered modification, it belongs to $J$. Otherwise, it is a reversion; since $\Aut(X,B)$ acts transitively on $F\setminus E$, $g_i=\alpha \Psi_1 \beta$, for some $\alpha,\beta\in \Aut(X,B)$. This achieves the proof of the equality $\Aut(S)=<\Aut(X,B),I,J>$.
 
Writing $A=<\Aut(X,B),I>$, the group $\Aut(S)$ is generated by $A$ and $J$. Let us prove that it is an amalgamated free product. Let $g=a_n\circ j_n\circ...\circ a_1\circ j_1$, where each $a_i\in A\setminus J$ and $j_i\in J\setminus A$. Then, $a_i$ is a product of reversions which is not an isomorphism, and $j_i$ is a fibered modification. Theorem \ref{Thm:Factorization} (or more precisely 
Lemma~\ref{Lem:UnicityFacto}) implies that $g$ does not belong to $\Aut(X,B)$ and then is not the identity. This shows that $\Aut(S)=A \star_{A\cap J} J$.

Assume that $n=2$. Then, $C\tr E$ is a zigzag of type $(-1,-2)$; the contraction of this zigzag gives rise to birational morphism of pairs $\nu:(X,B)\rightarrow (Z,D)$ for some projective surface $Z$, and some curve $D$. Furthermore, $\nu$ induces an isomorphism $F\rightarrow D$. Let us describe the pair $(Z,D)$, using the maps $(Z,D)\stackrel{\nu}{\leftarrow}(X,B)\stackrel{\mu_P}{\leftarrow}(Y,B)\stackrel{\eta_P}{\rightarrow}(\mathbb{F}_1,F\tr C\tr L)\stackrel{\tau}{\rightarrow}(\mathbb{P}^2,F\tr L)$. If $P$ has two distinct roots in $\kk$, then $\eta_P$ is the blow-up of two distinct points of $L_0=L\setminus C\subset \mathbb{F}_1$ and $\mu_P$ is an isomorphism. Since both $\nu$ and $\tau$ contract the same curve $C$, the birational map $(\mathbb{P}^2,F\tr L)\dasharrow (Z,D)$ consists of the blow-up of two distinct points of $F\setminus L$, followed by the contraction of $F$. This implies that $Z$ is isomorphic to a smooth quadric in $\mathbb{P}^3$ and that $D$  (which is the image of a line $F\subset  \mathbb{P}^2$) is an hyperplane section of $Z$. Moreover, if the two roots of $P$ are defined over $\k$, $Z$ is isomorphic to $\mathbb{P}^1\times\mathbb{P}^1$ and $D$ is a diagonal (i.e.\ a curve of bidegree $(1,1)$). If $P$ has one root of multiplicity two, then $\eta_P$ is the blow-up of a point $p_1\in L_0=L\setminus C\subset \mathbb{F}_1$, followed by the blow-up of the point $p_2$ in the first neighborhood of $p_1$, which belongs to the proper transform of the line $L$. Furthermore, $\mu_P$ consists of the contraction of the exceptional curve $E_{p_1}$ of $p_1$ (which is a $(-2)$-curve) on the unique singular point of $X$. Once again, both $\nu$ and $\tau$ contract the same curve $C$. The map $(\mathbb{P}^2,F\tr L)\dasharrow (Z,D)$ is therefore the composition of the blow-up of $p_1$, $p_2$ and the contraction of the two curves $E_{p_1}$ and $L$. The blow-up of $p_1$ goes to a surface isomorphic to $\mathbb{F}_1$, where the exceptional section is $E_1$ and where $C$ becomes a section of self-intersection $1$. Then, the blow-up of $p_2$ followed by the contraction of $F$ is an elementary link $\mathbb{F}_1\dasharrow \mathbb{F}_2$; the curves $E_1$ and $F$ become sections of self-intersection $-2$ and $2$ respectively. The contraction of $E_2$ gives the birational morphism $\mathbb{F}_2\rightarrow \mathbb{P}(1,1,2)=Z$. Now that $Z$ is described in each case, let us prove that $\Aut(Z,D)=\nu A\nu^{-1}$. Since each of the three curves $F$, $C$, $E$ is preserved by any automorphism of $(X,B)$, the group $\nu^{-1}\Aut(X,B)\nu$ is contained in $\Aut(Z,D)$, and corresponds in fact to the subgroup of elements of $\Aut(Z,D)$ which fix the point $\nu(C\cup E)$. Note that $\nu \Psi_1\nu^{-1}$ is an automorphism of $(Z,D)$ which sends the point $\nu(C\cup E)\in D$ onto the point $\nu(p)\in D$, where $p$ is the base-point of $\Psi_1$. Since the action of $\Aut(Z,D)$ on $D$ yields a surjective morphism $\rho:\Aut(Z,D)\rightarrow \Aut(D)\cong\PGL(2,\k)$, and because $\nu^{-1}\Aut(X,B)\nu$ contains the kernel of $\rho$ and its image by $\rho$ is a maximal group, then $\nu^{-1}\Aut(X,B)\nu$ and $\nu^{-1}\Psi_1\nu$ generate $\Aut(Z,D)$. This shows $7(a)$.

Assume now that $n\geq 3$, let $\Aut(X,B,\Psi_1)$ be the group of automorphisms of $(X,B)$ which fix the proper base-point $p$ of $\Psi_1$, and let $A_0=<\Aut(X,B,\Psi_1),\Psi_1>$. Then, clearly $A_0$ and $\Aut(X,B)$ (respectively $J$) generate $A$ (respectively $\Aut(S)$). Let us prove that we have an amalgamated free product in both cases. Let $g=a_n\circ j_n\circ...\circ a_1\circ j_1$, where each $a_i\in A_0\setminus \Aut(X,B)$ and $j_i\in \Aut(X,B)\setminus A_0$. Then, $a_i$ is a reversion centered at $p$ and $j_i$ is an automorphism of $(X,B)$ which moves $p$. Consequently, the decomposition $g=(a_n j_n)\circ ... \circ (a_1j_1)$ has no simplification and is minimal (Lemma~\ref{Lem:UnicityFacto}), so $g$ is not trivial. Assume now that each $a_i$ belongs to $A_0\setminus J$ and each $j_i$ belongs to $J\setminus A_0$. Once again, each $a_i$ is a reversion centered at $p$, and now $j_i$ is either an automorphism which moves $p$ or a fibered modification. We may group the $j_i$ which belongs to $\Aut(X,B)$ with $a_i$ and obtain a decomposition of $g$ of minimal length (applying once again Lemma~\ref{Lem:UnicityFacto}), so $g$ is not trivial. This yields $7(b)$. 

It remains to prove the explicit forms of $(8)$. Let $\psi=\eta_P\circ (\mu_P)^{-1}:(X,B)\dasharrow (\mathbb{F}_1, F\tr C)$, and recall that $\psi$ restricts to a birational morphism $S=X\setminus B\rightarrow \A^2=\mathbb{F}_1\setminus (F\cup C)$. According to Lemma~\ref{Lem:LiftTriangularC2}, $J=\psi^{-1} J'\psi$, where $J'$ is the group of eleements of $\dJo=\left\{(x,y)\mapsto (ax+P(y),by+c) \ |\ a,b\in\k^{*}, c\in \k, P\in\k[y])\right\}$ which preserve the points blown-up by $\psi^{-1}$ (or $\eta_P^{-1}$); furthermore, $\Aut(X,B)=\psi^{-1}(J'\cap \Aff) \psi$. The proper base-points of $\eta^{-1}$ are the points $(x_i,0)$ where $P(x_i)=0$. Furthermore, the other base-points lying on the transform of the line $L$ (which corresponds to $y=0$), $J'$ is the subgroup of elements of $\dJo$ which preserve the set of points of the form $(x_i,0)$ with $P(x_i)=0$. This means that $J'$ is generated by \[\begin{array}{rclcl}
H'&=&\{(x,y)\mapsto (x,a y)&|& a\in \k^{*}\},\\
Sp'&=&\{(x,y)\mapsto (ax+b,y)&|& a\in\k^{*}, b\in\k, P(ax+b)\mbox{ is a multiple of } P(x)\},\mbox{ and }\\
T'&=&\{(x,y)\mapsto (x+yQ(y),y)\!\!&|&Q\in \k[y]\}.\end{array}\]
The lift of these groups give respectively $H,Sp,T$, which generate $J$. Note that $J'\cap \Aff$ is generated by $H',Sp'$, and $T_0'=T'\cap \Aff$. The lift of these groups give $H,Sp,T_0$, which generate $\Aut(X,B)$. The proper base-point of $\Psi_1$ corresponds to $(0:1:0)\in\Pn$ (see Lemma~\ref{Lem:Rev01k}), which corresponds in $\A^2$ to the pencil of lines of the form $ax+b=0$. The group $\Aut(X,B,\Psi_1)$ is thus the lift of $<H',Sp'>$.  The remaining parts of $(8)$ follow directly.
\end{proof}

\subsection{$1$-standard pairs with a zigzag of type $(0,-1,-2,-3)$ or $(0,-1,-3,-2)$}
The surfaces with a zigzag of type $(0,-1,-n_1,-n_2)$ are the most simple immediately after the surfaces  described in the previous section. All these surfaces can give new examples of affine surfaces with unexpected properties. We give here the special case where the surface is smooth, the zigzag is of type $(0,-1,-2,-3)$ or $(0,-1,-3,-2)$, and where each component of the degenerate fibre is $\k$-rational. Properties distinct from the previous surfaces already show up in this simple example (Proposition~\ref{Prp:Graph23}). The general case will be treated in a forthcoming article.

Firstly, we describe a family of $1$-standard pairs (\S \ref{Exa:Pairs23}), and then prove that these are the only examples (Lemma \ref{Lem:IsoClass23}). We give the links between these maps by studying the possible reversions (Lemma~\ref{Lem:Rev23}), and then use this result to describe the properties of the $\mathbb{A}^1$-fibrations and of the automorphism group (Proposition~\ref{Prp:Graph23}).

\begin{enavant}\label{Exa:Pairs23}
We define here four families of $1$-standard pairs of surfaces $(X,B)$ of type $(0,-1,-2,-3)$ or $(0,-1,-3,-2)$, where $X\setminus B$ is affine and smooth. The map $\eta\colon X\to \mathbb{F}_1$ (as in \ref{Lem:GoingDownToF1}) is described here by its set of base-points which are in each case four points belonging, as proper or infinitely near points, to $L_0=L\setminus C\subset \mathbb{A}^2\subset \mathbb{F}_1$. Recall that $\mathbb{A}^2$ is viewed in $\mathbb{F}_1$ via the embedding $(x,y)\to ((x:y:1),(y:1))$, and that $L_0$ is the line of equation $y=0$ in $\mathbb{A}^2$ (see \ref{Explicitform}).

{\bf $\mathbf{I}$: Reduced case of type $\mathbf{(0,-1,-2,-3)}$}: there is only one surface here, called $(X_1,B_1)$. The map $\eta_1\colon X_1\to \mathbb{F}_1$ is the blow-up of $(0,0),(1,0)\in L_0$, and of the two points in the first neighbourhood of $(0,0)\in L_0$ corresponding to the two directions $x=0$ and $x=y$. $E_1\subset X_1$ is the proper transform of $L$ and $E_2\subset X_1$ is the curve obtained by blowing-up $(0,0)$. 
The following figure describes the morphism $(\mathbb{F}_1,F\tr C\tr L)\stackrel{\eta_1}{\leftarrow}(X_1,B_1)$.
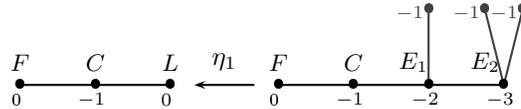
\begin{center}
\begin{pspicture}(1,0.7)(3,1.5)
\psline(1,1)(3,1)
\rput(1,1){\textbullet}\rput(1,1.3){{\small $F$}}\rput(0.95,0.8){{\scriptsize 0}}
\rput(2,1){\textbullet}\rput(2,1.3){{\small $C$}}\rput(1.95,0.8){{\scriptsize $-1$}}
\rput(3,1){\textbullet}\rput(3,1.3){{\small $L$}}\rput(2.95,0.8){{\scriptsize $0$}}
\end{pspicture}
\begin{pspicture}(0.4,0.7)(1.6,2)
\psline{<-}(0.6,1)(1.4,1)
\rput(1,1.3){$\eta_1$}
\end{pspicture}
\begin{pspicture}(1,0.7)(4,2)
{\darkgray 
\psline[linecolor=darkgray](3,1)(3,2)
\psline[linecolor=darkgray](4,1)(4.25,2)
\psline[linecolor=darkgray](4,1)(3.75,2)
\rput(3,2){\textbullet}
\rput(2.75,1.95){{\scriptsize $-1$}}
\rput(4.25,2){\textbullet}
\rput(4,1.95){{\scriptsize $-1$}}
\rput(3.75,2){\textbullet}
\rput(3.5,1.95){{\scriptsize $-1$}}
}
\psline(1,1)(4,1)
\rput(1,1){\textbullet}\rput(1,1.3){{\small $F$}}\rput(0.95,0.8){{\scriptsize 0}}
\rput(2,1){\textbullet}\rput(2,1.3){{\small $C$}}\rput(1.95,0.8){{\scriptsize $-1$}}
\rput(3,1){\textbullet}\rput(2.8,1.3){{\small $E_1$}}\rput(2.95,0.8){{\scriptsize $-2$}}
\rput(4,1){\textbullet}\rput(3.75,1.3){{\small $E_2$}}\rput(3.95,0.8){{\scriptsize $-3$}}
\end{pspicture}
\end{center}

{\bf $\mathbf{II}$: Reduced case of type $\mathbf{(0,-1,-3,-2)}$}: there is a family here, parametrised by a parameter $a\in \k\setminus \{0,1\}$. The pair is called $(X_{2,a},B_{2,a})$. The map $\eta_{2,a}\colon X_{2,a}\to \mathbb{F}_1$ is the blow-up of $(0,0),(1,0),(a,0)\in L_0$, and of the point in the first neighbourhood of $(0,0)\in L_0$ corresponding to the two direction $x=0$. $E_1\subset X_{2,a}$ is the proper transform of $L$ and $E_2\subset X_{2,a}$ is the curve obtained by blowing-up $(0,0)$. 
The following figure describes the morphism $(\mathbb{F}_1,F\tr C\tr L)\stackrel{\eta_{2,a}}{\leftarrow}(X_{2,a},B_{2,a})$.
\begin{center}
\begin{pspicture}(1,0.7)(3,1.5)
\psline(1,1)(3,1)
\rput(1,1){\textbullet}\rput(1,1.3){{\small $F$}}\rput(0.95,0.8){{\scriptsize 0}}
\rput(2,1){\textbullet}\rput(2,1.3){{\small $C$}}\rput(1.95,0.8){{\scriptsize $-1$}}
\rput(3,1){\textbullet}\rput(3,1.3){{\small $L$}}\rput(2.95,0.8){{\scriptsize $0$}}
\end{pspicture}
\begin{pspicture}(0.4,0.7)(1.6,2)
\psline{<-}(0.6,1)(1.4,1)
\rput(1,1.3){$\eta_{2,a}$}
\end{pspicture}
\begin{pspicture}(1,0.7)(4,2)
{\darkgray 
\psline[linecolor=darkgray](3,1)(3.25,2)
\psline[linecolor=darkgray](3,1)(2.75,2)
\psline[linecolor=darkgray](4,1)(4,2)
\rput(3.25,2){\textbullet}
\rput(3,1.95){{\scriptsize $-1$}}
\rput(2.75,2){\textbullet}
\rput(2.5,1.95){{\scriptsize $-1$}}
\rput(4,2){\textbullet}
\rput(3.75,1.95){{\scriptsize $-1$}}
}
\psline(1,1)(4,1)
\rput(1,1){\textbullet}\rput(1,1.3){{\small $F$}}\rput(0.95,0.8){{\scriptsize 0}}
\rput(2,1){\textbullet}\rput(2,1.3){{\small $C$}}\rput(1.95,0.8){{\scriptsize $-1$}}
\rput(3,1){\textbullet}\rput(2.75,1.3){{\small $E_1$}}\rput(2.95,0.8){{\scriptsize $-3$}}
\rput(4,1){\textbullet}\rput(3.8,1.3){{\small $E_2$}}\rput(3.95,0.8){{\scriptsize $-2$}}
\end{pspicture}
\end{center}

{\bf $\mathbf{III}$: Non-reduced case of type $\mathbf{(0,-1,-2,-3)}$}: there is a family here, parametrised by a parameter $a\in \k\setminus \{0,1\}$. The pair is called $(X_{3,a},B_{3,a})$. The map $\eta_{3,a}\colon X_{3,a}\to \mathbb{F}_1$ is the blow-up of $p_0=(0,0)\in L_0$, of the point $p_1$ in the first neighbourhood of $p_0$ corresponding to the two direction $x=0$, and of two points in the neighbourhood of $p_1$. In coordinates, $(u,v)\mapsto (u,u^2v)$ is the blow-up of $p_0$ and $p_1$, and the last two points correspond to $(u,v)=(0,1)$ and $(u,v)=(0,a)$. $E_1\subset X_{3,a}$ is the proper transform of $L$ and $E_2\subset X_{3,a}$ is the curve obtained by blowing-up $p_1$ (in the above coordinates it corresponds to $u=0$). 
The following figure describes the morphism $(\mathbb{F}_1,F\tr C\tr L)\stackrel{\eta_{3,a}}{\leftarrow}(X_{3,a},B_{3,a})$.

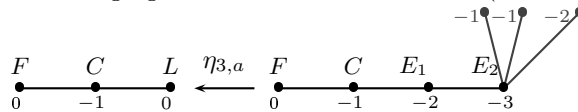
\begin{center}
\begin{pspicture}(1,0.7)(3,1.5)
\psline(1,1)(3,1)
\rput(1,1){\textbullet}\rput(1,1.3){{\small $F$}}\rput(0.95,0.8){{\scriptsize 0}}
\rput(2,1){\textbullet}\rput(2,1.3){{\small $C$}}\rput(1.95,0.8){{\scriptsize $-1$}}
\rput(3,1){\textbullet}\rput(3,1.3){{\small $L$}}\rput(2.95,0.8){{\scriptsize $0$}}
\end{pspicture}
\begin{pspicture}(0.4,0.7)(1.6,2)
\psline{<-}(0.6,1)(1.4,1)
\rput(1,1.3){$\eta_{3,a}$}
\end{pspicture}
\begin{pspicture}(1,0.7)(4,2)
{\darkgray 
\psline[linecolor=darkgray](4,1)(4.25,2)
\psline[linecolor=darkgray](4,1)(3.75,2)
\psline[linecolor=darkgray](4,1)(5,2)
\rput(4.25,2){\textbullet}
\rput(4,1.95){{\scriptsize $-1$}}
\rput(3.75,2){\textbullet}
\rput(3.5,1.95){{\scriptsize $-1$}}
\rput(5,2){\textbullet}
\rput(4.7,1.95){{\scriptsize $-2$}}
}
\psline(1,1)(4,1)
\rput(1,1){\textbullet}\rput(1,1.3){{\small $F$}}\rput(0.95,0.8){{\scriptsize 0}}
\rput(2,1){\textbullet}\rput(2,1.3){{\small $C$}}\rput(1.95,0.8){{\scriptsize $-1$}}
\rput(3,1){\textbullet}\rput(2.8,1.3){{\small $E_1$}}\rput(2.95,0.8){{\scriptsize $-2$}}
\rput(4,1){\textbullet}\rput(3.75,1.3){{\small $E_2$}}\rput(3.95,0.8){{\scriptsize $-3$}}
\end{pspicture}
\end{center}

{\bf $\mathbf{IV}$: Non-reduced case of type $\mathbf{(0,-1,-3,-2)}$}: there only one pair here, called $(X_4,B_4)$. The map $\eta_{4}\colon X_{4}\to \mathbb{F}_1$ is the blow-up of $(1,0)\in L_0$, of $p_0=(0,0)\in L_0$, of the point $p_1$ in the first neighbourhood of $p_0$ corresponding to the two direction $x=0$, and of one more point in the neighbourhood of $p_1$. In coordinates, $(u,v)\mapsto (u,u^2v)$ is the blow-up of $p_0$ and $p_1$, and the last point corresponds here to $(u,v)=(0,1)$. $E_1\subset X_{4}$ is the proper transform of $L$ and $E_2\subset X_{4}$ is the curve obtained by blowing-up $p_1$ (in the above coordinates it corresponds to $u=0$). 
The following figure describes the morphism $(\mathbb{F}_1,F\tr C\tr L)\stackrel{\eta_{4}}{\leftarrow}(X_{4},B_{4})$.

\begin{center}
\begin{pspicture}(1,0.7)(3,1.5)
\psline(1,1)(3,1)
\rput(1,1){\textbullet}\rput(1,1.3){{\small $F$}}\rput(0.95,0.8){{\scriptsize 0}}
\rput(2,1){\textbullet}\rput(2,1.3){{\small $C$}}\rput(1.95,0.8){{\scriptsize $-1$}}
\rput(3,1){\textbullet}\rput(3,1.3){{\small $L$}}\rput(2.95,0.8){{\scriptsize $0$}}
\end{pspicture}
\begin{pspicture}(0.4,0.7)(1.6,2)
\psline{<-}(0.6,1)(1.4,1)
\rput(1,1.3){$\eta_{4}$}
\end{pspicture}
\begin{pspicture}(1,0.7)(4,2)
{\darkgray 
\psline[linecolor=darkgray](3,1)(3,2)
\psline[linecolor=darkgray](4,1)(4,2)
\psline[linecolor=darkgray](4,1)(5,2)
\rput(3,2){\textbullet}
\rput(2.75,1.95){{\scriptsize $-1$}}
\rput(4,2){\textbullet}
\rput(3.75,1.95){{\scriptsize $-1$}}
\rput(5,2){\textbullet}
\rput(4.7,1.95){{\scriptsize $-2$}}
}
\psline(1,1)(4,1)
\rput(1,1){\textbullet}\rput(1,1.3){{\small $F$}}\rput(0.95,0.8){{\scriptsize 0}}
\rput(2,1){\textbullet}\rput(2,1.3){{\small $C$}}\rput(1.95,0.8){{\scriptsize $-1$}}
\rput(3,1){\textbullet}\rput(2.8,1.3){{\small $E_1$}}\rput(2.95,0.8){{\scriptsize $-2$}}
\rput(4,1){\textbullet}\rput(3.75,1.3){{\small $E_2$}}\rput(3.95,0.8){{\scriptsize $-3$}}
\end{pspicture}
\end{center}
\end{enavant}

\begin{lem}[Isomorphism classes of surfaces of type $(0,-1,-2,-3)$]\label{Lem:IsoClass23}
Let $(X,B=F\tr C\tr E,\overline{\pi})$ be a $1$-standard pair, such that $B$ is of type $(0,-1,-2,-3)$ or $(0,-1,-3,-2)$, such that $X\setminus B$ is smooth and affine and let $\eta:X\rightarrow \mathbb{F}_1$ be a birational morphism as in $\S \ref{Explicitform}$ above. Assuming that any component in the singular fibre of $\overline{\pi}$ is $\k$-rational, the following hold:

\begin{enumerate}
\item
there exist an automorphism $\alpha$ of $(\mathbb{F}_1,F\tr C\tr L)$ such that $\alpha\eta$ is equal to one of the morphisms $\eta_1,\eta_{2,a},\eta_{3,a},\eta_4$ defined in $\S\ref{Exa:Pairs23}$ above. In particular, $(X,B)$ is isomorphic to one of the pairs given in $\ref{Exa:Pairs23}$;
\item
any $1$-standard pair $(X',B',\overline{\pi}')$ such that $(X\setminus B,\pi)\cong(X'\setminus B',\pi')$ is isomorphic to $(X,B,\pi)$.
\end{enumerate}
\end{lem}
\begin{proof}
Since $S=X\setminus B$ is affine and smooth, only one fiber of $\overline{\pi}$ is singular and any irreducible component of this fibre touches $B$ or belongs to $B$. The self-intersections of the components in the boundary being given, $\eta\colon X\to \mathbb{F}_1$ is the blow-up of exactly four points. One checks that all possibilities are given in the four cases described in $\ref{Exa:Pairs23}$.

The second assertion can be checked directly, using the description of the base-points and applying Lemma~\ref{Lem:IsoFibTriangular}.\end{proof}

\begin{lem}\label{Lem:Rev23}
Let $(X,B)$ and $(X',B')$ be pairs described in $\ref{Exa:Pairs23}$, given with maps $\eta\colon X\to \mathbb{F}_1$ and $\eta'\colon X'\to \mathbb{F}_1$.
Suppose that there exists a reversion $\psi\colon(X,B)\dasharrow (X',B')$, centred at $p\in X$, whose inverse is centred at $p'\in X'$, where $p$ and $p'$ correspond respectively via $\eta$ and $\eta'$ to $((\lambda:1:0),(1:0)),((\lambda':1:0),(1:0))\in \mathbb{F}_1$ for some $\lambda,\lambda'\in \k$. 

Then, up to automorphisms of the pairs $(X,B)$, $(X',B')$, one of the following situations occurs for $\psi$ or its inverse, and every such situation can be realised:

$(1)$ $(X,B)=(X_1,B_1)$, $\lambda\in \k\setminus\{0,1\}$, $(X',B')=(X_{2,1-1/\lambda},B_{2,1-1/\lambda})$, $\lambda'\in \k^*$. 

$(2)$ $(X,B)=(X_1,B_1)$, $\lambda\in\{0,1\}$, $(X',B')=(X_{4},B_{4})$, $\lambda'\in \k$. 

$(3)$ $(X,B)=(X_{2,a},B_{2,a})$, $\lambda=0$, $(X',B')=(X_{3,a},B_{3,a})$, $\lambda'\in \k$.
\end{lem}
\begin{proof}
Let us fix some notation. We denote by $\epsilon_p\colon X_p\to X$ the blow-up of $p$, and by $\mathcal{E}_p\subset X_p$ the exceptional curve produced, and
write $B=F\tr C\tr E_1\tr E_2$, $B'=F'\tr C'\tr E_1'\tr E_2'$, and $\overline{\pi}'\colon X'\to \p^1$ the fibration associated to $F'$. We also denote by $\mathcal{R}'$ the set of components of the singular fibre of $\overline{\pi}'$ which intersect $B'$. There are $1$, $2$ or $3$ elements in $\mathcal{R}'$, depending in which family the pair $(X',B')$ is.
We will compute the number of elements of $\mathcal{R}'$ and their self-intersection using the information on $(X,B)$ and $\lambda$ to know in which of the four families the pair $(X',B')$ is.

Denote by $\mathcal{T}$ the set of curves of $X$ which are sent by $\psi$ on curves of $\mathcal{R}'$. It follows from the decomposition of $\psi$ given in \ref{DescriptionReversion1} (or from its resolution given in \ref{Exa:32resol}) that $\psi$ factors through $\epsilon_p$ and that $\psi_p=\psi\circ(\epsilon_p)^{-1}$ restricts to an isomorphism from $\mathcal{E}_p\backslash F$ to $E_2'\backslash E_1'$.
Moreover, if a curve of $\mathcal{R}'$ has self-intersection $-r$, the corresponding curve in $\mathcal{T}$  has self-intersection $-r+1$, and it intersects the boundary $B$ transversally and only at $p$. Since $r\in \{-1,-2\}$, the curve of $\mathcal{T}$ is the proper transform by $(\tau\eta)^{-1}$ of a line passing through $\tau\eta(p)=(\lambda:1:0)$ and through one or two points blown-up by $\eta$. 

We describe now the set of curves in $\mathcal{T}$ for each family and each $\lambda$.

(I) If $(X,B)=(X_{1},B_{1})$, the line of equation $z=x-\lambda y$ passes through $(\lambda:1:0)$ and through the point $(1:0:1)$ blown-up by $\eta_1$. Hence, its transform on $X$ gives an element of $\mathcal{T}$ of self-intersection $0$, and thus an element of $\mathcal{R}'$ of self-intersection $-1$. 
\begin{itemize}
\item[(Ia)] 
If $\lambda\notin\{0,1\}$, there is no other element of $\mathcal{T}$, hence $(X',B')$ is equal to $(X_{2,a},B_{2,a})$ for some $a\in \k\setminus\{0,1\}$.
\item[(Ib)] 
If $\lambda\in\{0,1\}$, the line $x=\lambda y$ passes through the point $(0:0:1)$, which is blown-up by $\eta_{1}$, \emph{and} by one of the two points in its neighbourhood which are also blown-up $\eta_{1}$. In this case, the transform of the line is an element of self-intersection $-1$ of $\mathcal{T}$, and gives and element of $\mathcal{R}'$ of self-intersection $-2$. In consequence, $(X',B')=(X_4,B_4)$.\end{itemize}

(II) If $(X,B)=(X_{2,a},B_{2,a})$ for some $a\in \k\setminus\{0,1\}$, the lines of equation $z=x-\lambda y$ and $az=x-\lambda y$ pass through $(\lambda:1:0)$ and respectively through $(1:0:1)$ and $(a:0:1)$. Hence, their transforms give two elements of $\mathcal{T}$ of self-intersection $0$ on $X$, and so two elements of $\mathcal{R}'$ of self-intersection $-1$. Moreover, the composition of the projection from $(\lambda:1:0)$ to the line $L\subset \mathbb{P}^2$ (of equation $y=0$) with $\epsilon_p$ gives rise to an isomorphism $\mathcal{E}_p\backslash F\to L\backslash F$, which induces with $\psi\circ (\eta_p)^{-1}$ an isomorphism $L\backslash F\to E_2'\backslash E_1'$. Call $D_0\subset X$ the proper transform of the line of $\p^2$ of equation $x=\lambda y$, which passes through $(\lambda:1:0)$ and the point $(0:0:1)$ blown-up by $\eta_{2,a}$.
\begin{itemize}
\item[(IIa)] 
If $\lambda\not=0$, $D_0$ has self-intersection $0$ in $X$ and intersects $E_2$. It does not belong to $\mathcal{T}$, so $(X',B')$ is equal to $(X_1,B_1)$. Moreover, $\psi(D_0)$ also has self-intersection $0$, and intersects $B'$ into two points, which are $p'=\psi(E_2)$ and $(\psi\circ (\eta_p)^{-1})(D_0\cap \mathcal{E}_p)\in E_2'$. In consequence, $\psi(D_0)$ is the lift by $(\eta_1)^{-1}$ of the line of equation $x=\lambda'y$.  The isomorphism $L\backslash F\to E_2'\backslash E_1'$ sends $(1:0:1)$ and $(a:0:1)$  onto the directions of the lines $x=y$ and $x=0$, each of these two curves passing through $(0:0:1)=\eta_1(E_2')$. Moreover, the point $(0:0:1)\in L$ is sent onto the direction of the line $x=\lambda'y$. Up to an exchange of the two directions $x=y$ and $x=0$ (which is induced by an automorphism of $(X_1,B_1)$), we obtain an automorphism of $\mathbb{A}^1$ which sends respectively $1,a,0$ onto $0,1,\lambda'$. This automorphism is $x\mapsto (x-1)/(a-1)$, so $\lambda'=1/(1-a)$.

\item[(IIb)] If $\lambda=0$, the curve $D_0$ belongs to $\mathcal{T}$, since the line of equation $x=\lambda y$ passes through the point $(0:0:1)$, \emph{and} by the point in its neighbourhood which is also blown-up $\eta_{2,a}$. In this case, $(D_0)^2=-1$ and $D_0$ corresponds to an element of $\mathcal{R}'$ of self-intersection $-2$. In consequence, $(X',B')=(X_{3,b},B_{3,b})$, for some $b\in \k\backslash \{0,1\}$. Let us prove that $b=a$. The isomorphism $\mathcal{E}_p\to E_2'$ sends respectively the direction of $z=x-\lambda y$, $az=x-\lambda y$ and $0=x-\lambda y$ onto the points corresponding to the three curves of $\mathcal{R}'$. Taking the coordinates $(u,v)$ as in the definition of family $\mathbf{III}$, the three points correspond respectively to $1$, $b$, $0$. We get an automorphism of $\mathbb{A}^1$ which sends respectively $1,a,0$ onto $1,b,0$. In consequence, $b=a$.  \end{itemize}

(III) If $(X,B)=(X_{3,a},B_{3,a})$ for some $a\in \k\setminus\{0,1\}$, the line of equation $x=\lambda y$ passes through $(\lambda:1:0)$ and through $(0:0:1)$, blown-up by $\eta_{3,a}$. Its transform is the unique element of $\mathcal{T}$, of self-intersection $0$ on $X$. Hence, $(X',B')=(X_{2,b},B_{2,b})$ for some $b\in \k\setminus\{0,1\}$. Moreover, $b=a$ since this reversion is the inverse of the one described in (IIa).

(IV) If $(X,B)=(X_{4},B_{4})$, the line of equation $z=x-\lambda y$ passes through $(\lambda:1:0)$ and $(1:0:1)$, blown-up by $\eta_{4}$. Its transform is the unique element of $\mathcal{T}$, and has self-intersection $0$ on $X$. Hence, $(X',B')=(X_{1},B_{1})$.

By the above list, there are three possible cases for $\psi$ or its inverse, which are $I\rightarrow II$, $I\rightarrow IV$ and $II\rightarrow III$. It remains to study each case and to give the values of the parameters associated to the surfaces or to the points.

$I\rightarrow II$. It follows from (Ia) and (IIa) that $(\lambda,\lambda')\in \k\backslash \{0,1\}\times\k^{*}$ and that each couple of this form is possible. We prove now that here the parameter of the surface $(X_{2,a},B_{2,a})$ is $a=1-1/\lambda$ (which proves in particular that any element of family $II$ can be obtained by a reversion on $(X_1,B_1)$). This link being the inverse of the one described in (IIa), the equality $\lambda= 1/(1-a)$, follows from the equality computed above.

$I\rightarrow IV$. The equality $\lambda=0$ follows from (Ib). Moreover, (IV) shows that  $\lambda'$ can take all possible values in $\k$.

$II\rightarrow III$. It follows from (IIb) that $\lambda=0$, and that the parameters of each pair are the same. Moreover, (III) shows that  $\lambda'$ can take all possible values in $\k$.
\end{proof}

\begin{prop}\label{Prp:Graph23}
All pairs described in $\ref{Exa:Pairs23}$ give the same affine surface $S$, up to isomorphism. Moreover, the graph $\mathcal{F}_S$ associated is the following:
\[\xymatrix@R=0.01cm@C=1cm{&&(X_{2,a},B_{2,a})&\ar@{<->}[l](X_{3,a},B_{3,a})\\
(X_4,B_4)\ar@{<->}[r]&(X_1,B_1)\ar@{<->}[ru]\ar@{<->}[rd]&\vdots&\vdots\\
&&(X_{2,b},B_{2,b})\ar@{<->}[r]&(X_{3,b},B_{3,b})}\]
where the $a,b$ correspond to all values in $\k\setminus\{0,1\}$, up to equivalence $a\sim a^{-1}$.

There are infinitely many equivalence classes of $\mathbb{A}^1$-fibrations on $S$ if and only if $k$ is infinite. Furthermore, $\Aut(S)$ is generated by the automorphisms of $\mathbb{A}^1$-fibration.
\end{prop}
\begin{rem}
The structure of $\Aut(S)$ can be described by this method; it is an amalgamated product of the group of automorphisms of $\mathbb{A}^1$-fibrations. 
\end{rem}
\begin{proof}
According to Lemma~\ref{Lem:Rev23}, we may obtain $(X_4,B_4)$ and any surface of type $(X_{2,a},B_{2,a})$ by applying a reversion on $(X_1,B_1)$. Applying a reversion on $(X_{2,a},B_{2,a})$, we get either $(X_1,B_1)$ or $(X_{3,a},B_{3,a})$.

Due to the descriptions of the families, $(X_{2,a},B_{2,a})$ is isomorphic to $(X_{2,b},B_{2,b})$ if and only if there exists an element of $\Aff$ which sends the points blown-up by $\eta_{2,a}$ onto points blown-up by $\eta_{2,b}$. This amounts to ask for the existence of an automorphism of the affine line $L_0\subset \mathbb{A}^2$ which fixes $(0,0)$ and sends $\{(1,0),(a,0)\}$ onto $\{(1,0),(b,0)\}$, and is thus equivalent to say that $a=b^{\pm 1}$. The case of family III is similar. Moreover, two pairs are isomorphic if and only if they induce the same affine fibred surfaces (Lemma~\ref{Lem:IsoClass23}). 

This gives the fact that all affine surfaces provided by the four families are isomorphic and also the description of the graph $\mathcal{F}_S$. We obtain the last assertion by applying  Proposition~\ref{Prp:FS}.
\end{proof}

\begin{rem}
In fact, taking $a,b\in \k^{*},c\in \k$, $a\not=b$, the following equations in $\mathbb{A}^4=\mathrm{Spec}(\k[w,x,y,z])$ define a smooth affine surface $S_{a,b,c}$, already studied in \cite{Dub2} (see also \cite{BandMakar}).
\[\begin{array}{rcl}
xz&=&y(y-a)(y-b)\\
yw&=&z(z-c)\\
xw&=&(y-a)(y-b)(z-c)\end{array}\]
The projection on the $x$-factor induces a $\mathbb{A}^1$-fibration which can be compactified by a pair of family $\mathbf{II}$ (\cite{Dub2}). In fact, we can check that the surface is $(X_{2,b/a},B_{2,b/a})$. The projection on the $w$-factor also gives an $\mathbb{A}^1$-fibration, and one can observe that this one belongs to family $\mathbf{I}$ if $c\not=0$ and to family $\mathbf{IV}$ otherwise. Proposition~\ref{Prp:Graph23} gives information on this affine surface and also shows that the isomorphism class does not depend of the parameters $(a,b,c)\in (\k^*)^2\times \k$.
\end{rem}
\bibliographystyle{amsplain}

\end{document}